\documentclass[reqno,11pt]{amsart}
\usepackage{amsmath,amsfonts,amssymb,amsxtra,latexsym,amscd,enumerate,amsthm,verbatim}
\usepackage{graphicx}

\usepackage[dvipsnames]{xcolor}
\usepackage{hyperref}
\usepackage{yhmath}
\hypersetup{colorlinks=true, pdfstartview=FitV, linkcolor=BrickRed,citecolor=BrickRed, urlcolor=BrickRed}

\usepackage{cancel,soul}

\newcommand{\f}{\frac}

\usepackage[margin=1.40in]{geometry}
\numberwithin{equation}{section}

\newcommand{\R}{\mathbb{R}}

\newcommand{\la}{\langle}
\newcommand{\ra}{\rangle}
\newcommand{\nn}{\nonumber}
\newcommand{\les}{\lesssim}

\numberwithin{equation}{section} %pour numeroter les equations par section

\newtheorem{theorem}{Theorem}[section]
\newtheorem{lemma}[theorem]{Lemma}

\newtheorem{corollary}[theorem]{Corollary}
\newtheorem{definition}[theorem]{Definition}
\newtheorem{remark}[theorem]{Remark}

\begin{document}

\title{Global Solutions for 5D Quadratic Fourth-Order Schr\"{o}dinger equations}

\author{Ebru Toprak}
\address{Yale University}
\email{ebru.toprak@yale.edu}
\author{Mengyi Xie}
\address{Yale University}
\email{mengyi.xie@yale.edu}

\maketitle

\begin{abstract}
We prove small data scattering for the fourth-order Schr\"odinger equation  with quadratic nonlinearity 
\begin{equation*}
	i\partial_t u+\Delta^2 u+\alpha u^2 + \beta \bar{u}^2=0\qquad\text{in }\mathbb{R}^5  
\end{equation*}
for $\alpha, \beta \in \R$. We extend the space-time resonance method, originally introduced by Germain, Masmoudi, and Shatah, to the setting involving the bilaplacian. We show that under a smallness condition on the initial data measured in a suitable norm, the solution satisfies $\|u\|_{L^{\infty}_x}\lesssim t^{-\frac{5}{4}} $ and scatters to the solution to the free equation. Although our work builds upon an established method, the fourth-order nature of the equation presents substantial challenges, requiring different techniques to overcome them.
\end{abstract}

\maketitle

\setcounter{tocdepth}{1}

\tableofcontents
%\textbf{Acknowledgement}

\section{Introduction}

Consider the nonlinear fourth order Schr\"odinger equation given by 
\begin{align} \label{nonf}
& iu_t= -\Delta^2 u + N(u,\bar{u}), \\ \nn 
& u(0,x)=u_0(x).  
\end{align}

Variants of this equation were introduced by Karpman \cite{Karp} and Karpman and Shagalov \cite{KarpSh}
to account for small fourth-order dispersion in the propagation of laser beams in a bulk
medium with Kerr nonlinearity. Additionally, in one dimension, this equation is related to the motion of an isolated vortex filament embedded in an inviscid, incompressible fluid, as proposed in \cite{Fukmof}. Beyond these applications, the equation provides a framework for modeling other "high-dispersion" phenomena.

In this paper, we investigate the small-data scattering problem for \eqref{nonf} in the case $d=5$ and $N(u) = \alpha u^2+ \beta \bar{u}^2$: 

\begin{align} \label{maineq}
& iu_t= -\Delta^2 u + \alpha u^2 + \beta \bar{u}^2 \\ \nn 
& u(1,x)= e^{-i\Delta^2}f_1
\end{align}
where $\alpha,\beta \in \R $. 

\subsection{Prior Work}
The equation \eqref{nonf} has been extensively studied.  As far as we have traced, the first work on the existence of solutions appears in \cite{FIP}, which initiated a line of research later expanded to address various dimensions; see, for example, \cite{GCui,HHWmult,Ben3,HHW,HJ,PNL, Aoki,HIT} and references therein. Among these, \cite{HHWmult} specifically considers nonlinearities that include quadratic terms of the type studied here. More broadly, the authors establish local well-posedness for nonlinearities of the form \( F(\partial^{\alpha}_{x} u, \partial^{\alpha}_{x} \bar{u}) \) for \( |\alpha| \leq 2 \) and \( d \geq 2 \). Their result, when applied to \eqref{maineq}, corresponds to the local well-posedness for initial data in \( H^{63/2}(\R^d)\cap H^{21} (\R^d; |x|^6 dx)\). In contrast, while our main result focuses on scattering, it also implies the global existence of solutions for small initial data in a function space with less regularity requirements (see \eqref{x_norm})

Before exploring the existing results on the scattering properties of \eqref{nonf}, it is helpful to first review its scaling properties. Let $N(u,\bar{u})$ be a homogeneous function with degree $p+1$. Then, \eqref{nonf} is invariant under the scaling: 
 $$
 u_{\lambda}(t,x)= \lambda^{\frac{4}{p}} u(\lambda^4 t, \lambda x), \,\,\ \lambda \in \R
 $$
Additionally, one may compute:
$$
\|u_{\lambda}\|_{\dot{H}^s}= \lambda^{\frac{4}{p} +s - \frac{d}{2}}\|u\|_{\dot{H}^s}, 
$$
which identifies the critical Sobolev space $\dot{H}^{s_c}$ with $s_c:= \frac{d}{2}-\frac{4}{p}$. In particular, the mass-critical nonlinearity power corresponds to $p=\frac{8}{d}$ and the energy-critical nonlinearity power corresponds to $p=\frac{8}{d-4}$. It is worth noting that the quadratic nonlinearity we consider lies in the mass-subcritical regime and necessitates a careful analysis of space-time resonances.

The scattering properties of \eqref{nonf} have been a topic of extensive investigation starting with \cite{seg1}, where the one-dimensional equation with cubic nonlinearity was analyzed. Over time, much of the subsequent work has centered on power-type nonlinearities, with the available methods largely enabling analysis in the mass supercritical regime.

In the defocusing case, where \( N(u,\bar{u}) = -|u|^p u \), the solutions to \eqref{nonf} are known to scatter in the \( H^2 \) space within the intercritical regime. Specifically, for dimensions \( 1 \leq d \leq 4 \), scattering has been established for \( \frac{8}{d} < p < \infty \), and for \( d \geq 5 \), the range is given by \( \frac{8}{d} \leq p < \frac{8}{d-4} \); see \cite{PX, GuaC, Saa, Dinh}. Furthermore, for the focusing case, scattering results are available in \cite{Pfoc, MXZ, PS, GQ}, where the energy-critical power in dimensions \( d \geq 5 \) and the intercritical power in every dimension are studied for radial initial data.

  On the other hand, in the energy supercritical regime, scattering has been demonstrated for \( p = 2 \) in dimensions \( 5 \leq d \leq 8 \) \cite{Ben3, MWZcubic}, and for \( d \geq 9 \), when \( s_c \geq 1 \), in \cite{MWZcubic}. Additionally, in \cite{RWZ}, small energy scattering is established in modulation spaces for \( d \geq 2 \) and for \( p  > \frac{8}{d} \) and \( p\) even. 

Among the existing works, only two studies have been able to address the mass subcritical regime. The first, in \cite{HMN}, explores the one-dimensional equation and establishes scattering for any \( p > 4 \). The second, in \cite{Wang}, studies the case for \( d \geq 2 \) with radial data, showing that scattering holds for \( s_c > - \frac{3n-1}{2n+1} \). Notably, for \( d = 5 \), this paper does not extend below \( p = \frac{14}{11} \). Furthermore, \cite{Wang} leverages the improved Strichartz estimate for the fourth-order equation in the radial setting—a property unavailable in non-radial configurations. Therefore, our work is the first to address the regime below the mass-critical point in multidimensional settings.

\subsection{Method and main result}\label{sec:method}

In the case of the nonlinear Schr\"odinger equation, it is well known that the scattering problem becomes increasingly challenging in lower dimensions when the nonlinearity has low homogeneity. Under such conditions, the most effective approach to establishing scattering is the space-time resonances method, which has proven pivotal since its introduction in \cite{GMS}.

%This is because most of the existing methods are effective only when the  homogeneity of the nonlinearity are sufficiently high to ensure strong dispersive effects.

Before presenting our main theorem, we wish to pay tribute to the elegant and profoundly impactful method of space-time resonances.  Consider the nonlinear Schr\"odinger equation in $\R^3$: 
 \begin{align} \label{nneq}
 i \psi_t = -\Delta \psi+ N(\psi, \bar{\psi}), \,\,\,  \psi(0,x)=\psi_0(x)
 \end{align}
whose solutions are represented by Duhamel's formula as 
  \begin{align*}
\psi(t,x) & = e^{-it\Delta} \psi_0(x) + \int_0^t e^{-i(t-s)\Delta} N(\psi(s))  ds \\
& = e^{-it\Delta} \Big[ \psi_0(x) + T(N(\psi))(x) \Big]
 \end{align*}
where $T(f)(x):= \int_0^{t} e^{is\Delta} f(s) ds$. 
Given the following well-known inequality
\begin{align}\label{lin:sch}
\|e^{-it\Delta}g\|_{L^{\infty}} \leq C t^{-\frac{3}{2}} \Big[ \| \hat{g}\|_{L^{\infty}} + t^{-\frac{1}{4}} \|x^2 g\|_{L^2} ]
\end{align}
the goal of small-data scattering is  to identify an appropriate norm $X$  such that 
$$
\| \widehat{T(N(\psi))}\|_{L^{\infty}_t L^{\infty}_x} +  \|t^{-\frac{1}{4}} x^2 T(N(\psi))\|_{L^{\infty}_t L^2} \leq C_1 \| T(N(\psi))\|_X \leq C_2 \|\psi\|^{p+1}_{X} 
$$
which would allow the application of the contraction mapping argument to establish the existence of solutions that exhibit the natural decay rate of the free equation, $t^{-\frac{3}{2}}$. The norm $X$ is typically expressed in terms of bounds on the profile of the solution, $f= e^{it\Delta} \psi$, which leads to the conclusion of scattering behavior.

  When dealing with $N(\psi)=\psi^2$ in $\R^3$, it becomes significantly more challenging to control $\|x^2 T(N(\psi))\|_{L^2}$, which corresponds to managing the second derivative of the following expression on the Fourier side:  
  $$
  \int_0^{t} \int_{\R^3}e^{is(|\xi|^2- |\xi-\eta|^2 -|\eta|^2)} \hat{f}(\xi) \hat{f}(\xi-\eta) d \eta ds.  
  $$

In particular, integration by parts is necessary to gain integrability in \( s \), as the nonlinearity alone is insufficient to fully recover the time integration. Moreover, it is necessary to eliminate the polynomial that depends on \( \xi \) and \( \eta \), which arises from differentiating the exponential term. However, the phase contains critical points: specifically, in the \( \eta \)-integral, a critical point arises when
\[
\partial_{\eta}\{ |\xi - \eta|^2 + |\eta|^2 \} = 0,
\]
which occurs at \( \xi = 2\eta \), referred to as the space resonance. Similarly, in the \( s \)-integral, a critical point appears when
\[
|\xi|^2 - |\xi - \eta|^2 - |\eta|^2 = 0,
\]
which occurs at \( \xi = \eta \), referred to as the time resonance.
Here, the space-time resonances method becomes essential. This technique takes advantage of the fact that the space resonances and time resonances intersect only at the origin, enabling simultaneous integration by parts in both \(s\) and \(\eta\).

In this paper, we extend this method to the fourth-order Schr\"odinger equation in $\R^5$. Let $f$ denote the profile of $u$ in \eqref{maineq}: $f=e^{it\Delta^2}u$. By Duhamel Formula 
\begin{align}\label{Duhamel}
	\hat{f}(\xi,t)&=\widehat{f_1}(\xi)+ \frac{ i \alpha}{(2 \pi)^{\f 52}} \int_1^t ie^{is|\xi|^4}\widehat{u^2}(\xi,s)ds  \nn \\ 
	&=\widehat{f_1}(\xi)+\frac{ i \alpha}{(2 \pi)^{\f 52}}  \int_1^t \int_{\mathbb{R}^5}e^{is\varphi(\xi,\eta)}\widehat{f}(\xi-\eta,s)\hat{f}(\eta,s)d\eta ds 
\end{align}
where  $\varphi=\varphi(\xi,\eta)=|\xi|^4-|\eta|^4-|\xi-\eta|^4$. Here, we consider when $\beta=0$. However, the method will also work when $\beta \neq 0$. 

If we denote space resonances with $\mathcal{S}$ and time resonances with $\mathcal{T}$, simple computations give 
\begin{equation*}
	\begin{split}
		\mathcal{S}=\{\partial_\eta \varphi=0  \}=\{(\xi,\eta) \ &|\ |\xi-\eta|^2(\xi -\eta)=|\eta|^2\eta \}=\{\xi=2\eta\}\\
		\mathcal{T}=\{\varphi=0\}=\{&(\xi,\eta) \ |\ |\xi|^4-|\eta|^4=|\xi-\eta|^4 \}\\
		\mathcal{R}&=\mathcal{S}\cap\mathcal{T}=\{(0,0)\}
	\end{split}
\end{equation*}
Thus, the space-time resonance method appears to be a suitable approach. In this paper, we prove 
\begin{theorem} \label{th:main} Let $f$ denote the profile of $u$ in \eqref{maineq} and define the following Banach space $X$ by its norm 
\begin{multline}\label{x_norm}
\|f\|_X=\|\widehat{f}\|_{L^\infty L^\infty}+ \|f\|_{L^\infty L^2}+\|xf\|_{L^\infty L^2}+\|\frac{x^2}{\log(t)}f\|_{L^\infty L^2} \\ +\|\frac{x^3}{t^{\alpha}}f\|_{L^\infty L^2}+\|t^{\frac{5}{4}}u\|_{L^\infty L^\infty}+\|f\|_{L^\infty H^{10}} 
\end{multline}
where $\alpha=\frac{1}{2}+\frac{1}{47}$. There exists a solution to \eqref{nonf} in $X$
for any $f_1$ such that  $\| f_1\|_{X} <\delta $, for some $\delta$ small enough. Furthermore, this solution scatters in $L^2(\R^5)$.
\end{theorem}
\begin{remark}\label{rmk_main}
	The number $\frac{1}{2}+\frac{1}{47}$ is not special. In fact, it can be shown that $\alpha$ can be pushed to $\frac{1}{2}+$.
\end{remark}
Although our work builds upon an established method, the fourth-order nature of the equation presents substantial challenges, requiring different techniques to overcome them. In the next section, as we discuss the setup and the proof outline, we address the terms appearing in \eqref{x_norm} and elaborate on these difficulties. We finish this section by giving the organization of the paper and the notations. 

\subsection{Organization of the paper} The paper is organized as follows. In Section~\ref{sec:setup}, we outline the proof. In Section~\ref{sec:technical}, we revisit fractional integration and Coifman-Meyer operators with flag singularities. This section also includes the proofs of technical lemmas crucial for the subsequent estimates. In Section~\ref{sec:dispersion}, we prove dispersive estimates for the solution to the linear equation.  In Sections~\ref{sec:gbounds} and~\ref{sec:hbounds}, we prove the \(L^2\) and weighted \(L^2\) bounds. Section~\ref{sec:hbounds}
also includes the proof of $\|t^{-\f54}u\|_{L_t^{\infty}L_x^{\infty}}$. Finally, in Section~\eqref{proof}, we close the bootstrap argument and complete the proof of Theorem~\ref{th:main}. 
\subsection{Notations} 
\begin{itemize}
\item  $ f(x) \les g(x)$ indicates  that there exists a uniform  constant $C>0$ so $ f(x) \leq Cg(x)$ for all $x$.
\item Let $\psi$ be a smooth nonnegative real radial function supported in $\{\frac{3}{4}<|\xi|<\frac{8}{3}\}$ such that 
$$\sum_j \psi^2\left(\frac{\xi}{2^j}\right)=1, \,\,\,\ \xi\neq 0 . $$
We define $\psi_j(\cdot):=\psi(\frac{\cdot}{2^j})$, $P_j f = \psi_j * f$ and  $P_{\leq j}:=\sum_{k\leq j}P_k$. 
\item Let $\varphi$ be a smooth nonnegative real radial function supported in $\{\frac{3}{4}<|\xi|<\frac{8}{3}\}$ such that $$\sum_j \varphi\left(\frac{\xi}{2^j}\right)=1, \,\,\,\ \xi\neq 0. $$
We define $\varphi_j(\cdot):=\varphi(\frac{\cdot}{2^j})$. With a slight abuse of notation, we also write \( P_j f \) when referring to  \( \varphi_j * f \).

\item We define $\Psi_1, \Psi_2$ as homogeneous functions of $\xi$ and $\eta $ with degree 0 which are $\mathcal{C}^\infty $ away from $(0,0)$ and satisfy 
\begin{equation}\label{cut_2}
	\begin{cases}
	\Psi_1(\xi,\eta)+\Psi_2(\xi,\eta)=1\\
	\text{on Supp$(\Psi_1)$   \ \ \ \   }|\eta|\leq 4|\xi-\eta|\\
	\text{on Supp$(\Psi_2)$   \ \ \ \   }|\xi-\eta|\leq 4|\eta|\\
\end{cases}
\end{equation}
\item Similarly, we define $\widetilde{\Psi_1}, \widetilde{\Psi_2}$ as homogeneous functions of $\xi$ and $\eta $ with degree 0 which are $\mathcal{C}^\infty $ away from $(0,0)$ and satisfy 
\begin{equation}\label{cut_2_2}
	\begin{cases}
	\widetilde{\Psi_1}(\xi,\eta)+\widetilde{\Psi_2}(\xi,\eta)=1\\
	\text{on Supp$(\widetilde{\Psi_1})$   \ \ \ \   }|\xi|\leq 4|\xi-\eta|\\
	\text{on Supp$(\widetilde{\Psi_2})$   \ \ \ \   }|\xi-\eta|\leq 4|\xi|\\
\end{cases}
\end{equation}
\item We define $\phi_1(\xi,\eta,\sigma), \phi_2(\xi,\eta,\sigma), \phi_3(\xi,\eta,\sigma)$ as homogeneous functions with degree 0 which are  $\mathcal{C}^\infty $ away from $(0,0)$ and satisfy
\begin{equation}\label{cut_3}
	\begin{cases}
	\phi_1(\xi,\eta,\sigma)+\phi_2(\xi,\eta,\sigma)+\phi_3(\xi,\eta,\sigma)=1\\
	\text{on Supp$(\phi_1)$   \ \ \ \   } |\eta-\sigma|,|\sigma|\leq 4|\xi-\eta| \\
	\text{on Supp$(\phi_2)$   \ \ \ \   }|\xi-\eta|,|\sigma|\leq 4|\eta-\sigma|\\
	\text{on Supp$(\phi_3)$ \ \ \ \ } |\xi-\eta|, |\eta-\sigma|\leq 4|\sigma|
\end{cases}
\end{equation}
\item We use $Q_k$ to denote homogeneous polynomials of degree $k$ in its variables
\item We will use both $\|\cdot\|_{L^p}$ and $\|\cdot\|_p$ to denote $L^p$ spaces
\end{itemize}

\section{Set-up and Outline of the proof}\label{sec:setup}

It is known that the evolution of the fourth-order Schr\"odinger equation satisfies \(\|e^{it\Delta^2}f\|_{L^{\infty}} \les t^{-\frac{d}{4}} \|f\|_{L^1}\). In \cite{BKS}, Ben-Artzi, Koch, and Saut demonstrated this result using Bessel's function expansion. However, their method is unsuitable for obtaining a linear estimate in the form of \eqref{lin:sch}. In Section~\ref{sec:dispersion}, we use multidimensional oscillatory analysis to prove (see Lemma~\ref{lem:linear}) 
\begin{align}\label{lin:fsch}
\| e^{it\Delta^2}f\|_{L^\infty} \les t^{-\frac{5}{4}} \|\hat{f}\|_{L^{\infty}} + t^{-\frac{11}{8}} \| \hat{f} \|_{H^3}.
\end{align}

As a result of this estimate, it is expected that all terms in the \(X\)-norm \eqref{x_norm} are included, except for the last one. It turns out that the last term is required for a technical reason and can be relaxed to \(L^{\infty} H^{5+}\), see Remark~\ref{middle_j} and the analysis of the terms in \eqref{H9}. 

 When \eqref{lin:fsch} is compared to \eqref{lin:sch}, two observations emerge immediately. First, the gain of decay in the second term is lower in \eqref{lin:fsch}, and second, the appearance of the \(\|x^3 f\|_{L^2_x}\) term. Both of these aspects introduce disadvantages. Luckily, the estimate \eqref{lin:fsch} is optimal for small frequencies. In particular, the group velocity of the fourth-order free wave is proportional to \(\xi|\xi|^2\), indicating that small frequencies travel slower, while large frequencies travel faster compared to the Schr\"odinger free particle, see Lemma~\ref{lem:ben}. Indeed, the fact that large frequencies are faster plays a crucial role in reducing \(\alpha\) as small as \(\frac{49}{94}\).

As explained in Section~\ref{sec:method}, we use the Banach fixed-point theorem, i.e., we aim to show that:
\begin{align}\label{contraction}
	\hat{f}(\xi,t)\rightarrow & \widehat{f_1}(\xi)+\frac{i \alpha}{(2\pi)^{\f 52}}\int_1^t\int e^{is\varphi}\widehat{f}(\xi-\eta,s)\hat{f}(\eta,s)d\eta ds \nn \\
    =:& \widehat{f_1}(\xi) + \frac{ i \alpha}{(2 \pi)^{\f 52}} \hat{B}(f,f)
\end{align} 
  is a contraction mapping on a small neighborhood of the origin with respect to the $X$ norm. More precisely, our goal is to show the following estimate:
\begin{equation}
	 \| B(f,f) \|_X\lesssim\|f\|_{X}^2.
\end{equation}

\subsection{Set-up for space-time resonance for bilinear terms}
To follow the space-time resonances strategy, we first look for a function \( P( \xi, \eta) \) such that \( Z := \varphi + P \cdot \varphi_\eta \) vanishes only at \(\eta = \xi = 0\), the set of space-time resonances. Specifically, we define
\begin{equation*}
	P=-\eta+\frac{\xi}{5}
\end{equation*}
and compute:
\begin{align}
	Z&=6|\eta|^4+\frac{4}{5}|\xi|^4+\frac{28}{5}(\eta\cdot \xi)^2-\frac{48}{5}|\eta|^2(\eta\cdot \xi)-\frac{12}{5}|\xi|^2(\eta\cdot \xi)+\frac{15}{4}|\eta|^2|\xi|^2 \nn \\
	&=|a\eta-n\xi|^4+\frac{6}{5}|\eta|^4+\frac{1}{2}|\xi|^4+\frac{4}{5}(\eta\cdot \xi)^2+\frac{2}{5}|\eta|^2|\xi|^2 \nn \\
	&\geq \frac{1}{2}(|\eta|^4+|\xi|^4)
	\gtrsim \max((|\xi-\eta|^4+|\xi|^4),|\xi-\eta|^4+|\eta|^4) \label{Zlowb}
\end{align}
with $a^2=\frac{12}{\sqrt{30}}$ and $b^2=\frac{3}{\sqrt{30}}$. This computation indicates that \( Z \) vanishes only at \((\xi, \eta) = (0, 0)\).   Noting that $\partial_s e^{is\varphi}=i\varphi e^{is\varphi}$ and $\partial_\eta xe^{is\varphi}=is\varphi_\eta e^{is\varphi}$, we have the following identity
\begin{equation}\label{key2}
	e^{is\varphi}=\frac{1}{\frac{1}{s}+iZ}(\partial_s+\frac{1}{s}+\frac{P}{s}\cdot \partial_\eta)e^{is\varphi}.
\end{equation}
\subsection{Decomposition of $f$} The space-time resonance set-up also naturally lends itself to a decomposition of $f$, allowing us to take advantage of the distinct behaviors of each component. Using \eqref{key2} in \eqref{contraction}, we may write
\begin{equation}\label{decompB}
\begin{split}
	\hat{B}(f,f)&=i\int_1^t \int_{\mathbb{R}^5} \frac{1}{\frac{1}{s}+iZ}(\partial_s+\frac{1}{s}+\frac{P}{s}\partial_\eta)e^{is\varphi}\widehat{f}(\xi-\eta,s)\hat{f}(\eta,s)d\eta ds\\
	&=:\widehat{f_*} +\hat{g}(\xi,t)+\hat{h}(\xi,t)
\end{split}	
\end{equation}
where
\begin{align}
 & \widehat{f_*}(\xi)=-i\int_{\mathbb{R}^5}\frac{1}{1+iZ}e^{i\varphi}\widehat{f}_1(\xi-\eta) \hat{f}_1(\eta) d\eta \nn \\
 & \hat{g}(\xi,t)=i\int_{\mathbb{R}^5}\frac{1}{\frac{1}{t}+iZ}e^{it\varphi}\widehat{f}(\xi-\eta,t)\hat{f}(\eta,t)d\eta \label{eqn_g}
\end{align}
and
\begin{equation}\label{eqn_h}
	\begin{split}
		\hat{h}(\xi,t)&=-i\int_1^t \frac{1}{s^2}\int_{\mathbb{R}^5}\frac{1}{(\frac{1}{s}+iZ)^2}e^{is\varphi}\widehat{f}(\xi-\eta,s)\hat{f}(\eta,s)d\eta ds\\
		&+2i\int_1^t \int_{\mathbb{R}^5} \frac{1}{\frac{1}{s}+iZ}e^{is\varphi}\widehat{f}(\xi-\eta,s)\partial_s\hat{f}(\eta,s)d\eta ds\\
		&+i\int_1^t \frac{1}{s}\int_{\mathbb{R}^5} \frac{1}{\frac{1}{s}+iZ}e^{is\varphi}\widehat{f}(\xi-\eta,s)\hat{f}(\eta,s)d\eta ds\\
		&+i\int_1^t \frac{1}{s}\int_{\mathbb{R}^5} \frac{P}{\frac{1}{s}+iZ}\partial_\eta(e^{is\varphi})\widehat{f}(\xi-\eta,s)\hat{f}(\eta,s)d\eta ds\\
		&=:\widehat{h_1}(\xi,t)+\widehat{h_2}(\xi,t)+\widehat{h_3}(\xi,t)+\widehat{h_4}(\xi,t)
	\end{split}
\end{equation}
 Intuitively, $e^{it\Delta^2}g$ has better pointwise decaying behavior, whereas $h$ is the collection of terms that are highly localized in the $L^2$-norm. In particular, the maximal growth \(t^{\alpha}\) is due to the \(g\) term, whereas it satisfies the following pointwise decay: \(e^{it\Delta^2}g = O(t^{-\frac{3}{2}+}) \|\hat{g}\|_{L^{\infty}}\). More precisely, we show the following estimates:
\begin{equation}\label{est_g}
	\begin{split}
	\|g\|_{L^2}\lesssim \|f\|_{X}^2, \qquad \|xg\|_{L^2}&\lesssim \|f\|_{X}^2, \qquad \|x^2g\|_{L^2}\lesssim\|f\|_X^2+\|f\|_X^3, \\
		\qquad \|x^3g\|_{L^2}\lesssim t^{\frac{1}{2}+\frac{1}{47}} \|f\|_{X}^2, &\qquad \|e^{-it\Delta^2}g\|_{L^\infty}\lesssim t^{-\frac{3}{2}+}\|f\|_{X}^2
	\end{split}
\end{equation}
\begin{equation}\label{est_h}
	\begin{split}
		\|h\|_{L^2}\lesssim \|f\|_{X}^2+\|f\|_X^3, \qquad \|xh\|_{L^2}&\lesssim \|f\|_{X}^2, \qquad \|x^2h\|_{L^2}\lesssim\log(t)(\|f\|_{X}^2+\|f\|_X^3), \\
		\qquad \|x^3h\|_{L^2}\lesssim t^{\frac{1}{24}} (\|f\|_{X}^2+\|f\|_X^3&+\|f\|_X^4), \qquad \|e^{-it\Delta^2}h\|_{L^\infty}\lesssim t^{-\frac{5}{4}}\|f\|_{X}^2
	\end{split}
\end{equation}
In particular, we will establish the weighted $L^2$-norms via a bootstrap argument. Therefore, the proof will rely on the following a priori estimates for $g$ and $h$: 

\begin{equation}\label{Boot_est_g}
	\begin{split}
	\|g\|_{L^2}\lesssim \log(t)\|f\|_{X}^2, &\qquad \|xg\|_{L^2}\lesssim \log(t) \|f\|_{X}^2\\ 
\|x^2g\|_{L^2}\lesssim \log(t) \|f\|_X^2,
		&\qquad \|x^3g\|_{L^2}\lesssim t^{\frac{1}{2}+\frac{1}{45}} \|f\|_{X}^2
	\end{split}
\end{equation}
\begin{equation}\label{Boot_est_h}
	\begin{split}
		\|h\|_{L^2}\lesssim \log(t)\|f\|_{X}^2,& \qquad \|xh\|_{L^2}\lesssim \log(t)\|f\|_{X}^2\\
		\|x^2h\|_{L^2}\lesssim t^{\varepsilon_0}\|f\|_{X}^2, 
		&\qquad \|x^3h\|_{L^2}\lesssim t^{\frac{1}{16}} \|f\|_X^2
	\end{split}
\end{equation}
where $\varepsilon_0$ is an arbitarily small number, like $\frac{1}{10000}$. We will prove \eqref{est_g} and \eqref{est_h} using these assumptions. As suggested in Remark \ref{rmk_main}, the numbers $\frac{1}{47}, \frac{1}{24}, \frac{1}{16}$ are not special, but just picked to close the bootstrap. 
%Therefore, we have:
%\begin{equation}\label{est_f}
%	\begin{split}
%		\|f\|_{L^2}\lesssim \|f\|_{X}^2, \qquad \|xf\|_{L^2}&\lesssim \|f\|_{X}^2, \qquad \|x^2f\|_{L^2}\lesssim\log(t)(\|f\|_X^2+\|f\|_X^3), \\
%		\qquad \|x^3f\|_{L^2}\lesssim t^{\frac{1}{2}+\frac{1}{47}} \|f\|_{X}^2, &\qquad \|e^{-it\Delta^2}f\|_{L^\infty}\lesssim t^{-\frac{5}{4}}\|f\|_{X}^2
	%\end{split}
%\end{equation}

This decomposition \eqref{decompB} plays a significant role in controlling the local decay of the terms in the bilinear operator of $B(f,f)$. In particular, the components involving $g$ enable us to transform the bilinear operators into trilinear operators, yielding an additional $t^{-\frac{1}{4}}$ decay. Meanwhile, the components involving only the $h$ term exhibit better local decay properties.  For example, if we let  $L(f,f)$ be a bilinear operator, we can decompose 
\begin{equation}\label{Decom_Bi}
	\begin{split}
		L(f,f)=L(f,f_1+f_*)+L(f,g)+L(f_1+f_*,h)+L(g,h)+L(h,h)
	\end{split}
\end{equation}
The terms involving $f_1+f_*$ are easier to estimate as \(\|f_1+f_*\|_{X} \) is small independent of $t$ . Since $h$ has better local behavior than $f$, $L(g,h)$ is easier to work with than $L(f,g)$. Thus, the major difficulties we have are with terms $L(f,g)$ and $L(h,h)$. . %, but they can be resolved using the previous observations. 

\subsection{Set-up for space-time resonance for trilinear terms.}
Controlling the third derivative of $ \hat{B}(f,f) $ necessitates the use of the space-time resonance strategy for the associated trilinear terms of the form: 
\begin{align*}
 \int_1^{\infty} \int_{\R^5}\int_{\R^5} \frac{Q_k(\xi,\eta)}{\f 1{s} +iZ }e^{it\psi(\xi,\eta,\sigma)} \hat{f}(\xi-\eta)  \hat{f}(\eta-\sigma)\hat{f}(\sigma) \: d\sigma \: d \eta \: ds
\end{align*}
where 
\begin{equation}\label{phase_3}
	\psi=\psi(\xi,\eta,\sigma)=|\xi|^4-|\xi-\eta|^4-|\eta-\sigma|^4-|\sigma|^4.
\end{equation}
We may compute
\begin{equation*}
	\begin{split}
		\mathcal{S_{(\xi,\eta,\sigma)}}=\{\partial_{\eta}\psi=0=\partial_\sigma \psi\}=\{ |\xi-\eta|^2(\xi-\eta)=&|\eta-\sigma|^2(\eta-\sigma)=|\sigma|^2\sigma \}=\{\xi=3\sigma=\frac{3}{2}\eta \}\\
		\mathcal{T_{(\xi,\eta,\sigma)}}=\{\psi=0 \}=\{|\xi|^4=|\xi-&\eta|^4+|\eta-\sigma|^4+|\sigma|^4  \}\\
		\mathcal{R_{(\xi,\eta,\sigma)}}=S_{(\xi,\eta,\sigma)}&\cap T_{(\xi,\eta,\sigma)}=\{0\}. 
	\end{split}
\end{equation*}
Moreover, we  pick two polynomials:
\begin{equation*}
	Q(\xi,\eta,\sigma)=2\xi-3\eta, \qquad S(\xi,\eta,\sigma)=\xi-3\sigma.
\end{equation*}
We then have
\begin{equation*}
	\begin{split}
		Y:&=\psi+Q\cdot \psi_\eta+S\cdot \psi_\sigma \geq 4(|\xi-\eta|^4+|\eta-\sigma|^4+|\sigma|^4)
	\end{split}
\end{equation*}
which only vanishes at zero. Similar to \eqref{key2}, we have:
\begin{equation}\label{key3}
	e^{is\psi}=\frac{1}{\frac{1}{s}+iY}\left(\frac{1}{s}+\partial_s+\frac{Q}{s}\cdot \partial_\eta+\frac{S}{s}\cdot \partial_\sigma \right)e^{is\psi}. 
\end{equation}

This identity allows us simultaneous integration by parts in \(\eta\), \(s\), and \(\sigma\). However, after performing integration by parts, the resulting expressions are Coifman–Meyer operators with flag singularities which we provide more details on in the next section.

\section{Fractional Integration and Pseudo-Product Operators} \label{sec:technical}
\subsection{Fractional Integration} As suggested in \eqref{key2} and \eqref{eqn_g}, we work  with multipliers with denominator like $\frac{1}{t}+iZ$. Thus, we need to study the fractional integration for $\Delta^2$:
\begin{equation}
	\nabla_t^\alpha:=(1/t+\Delta^2)^{\frac{\alpha}{4}}
\end{equation}
We have the following useful lemma  whose proof is adapted from the standard proof for the fractional integration for $\Delta$. 
\begin{lemma}\label{lem_frac_int}
	For $\alpha\geq 0$, $1\leq p,q<\infty$ and $0\leq \frac{1}{q}-\frac{1}{p}\leq\frac{\alpha}{5}$; or $1\leq p,q\leq\infty$ and $0\leq \frac{1}{q}-\frac{1}{p}<\frac{\alpha}{5}$, we have:
	\begin{equation}\label{frac_int}
		\|\nabla_t^{-\alpha}f\|_{L^p}\lesssim t^{\frac{\alpha}{4}+\frac{5}{4}(\frac{1}{p}-\frac{1}{q})}\|f\|_{L^q}
	\end{equation}
\end{lemma}
\begin{proof}
 For $1\leq p,q<\infty$, we have 
	\begin{equation*}
		\begin{split}
			\|\nabla_t^{-\alpha}f\|_{L^p}&=t^{\frac{\alpha}{4}-\frac{5}{4}+\frac{5}{4p}}\|\mathcal{F}^{-1}((1+|\eta|^4)^{-\frac{\alpha}{4}}\hat{g})\|_{L^p} \ \ \text{where }\hat{g}(\eta)=\hat{f}(t^{-\frac{1}{4}}\eta)\\
			&\approx t^{\frac{\alpha}{4}-\frac{5}{4}+\frac{5}{4p}}\|g\|_{W^{-\alpha,p}}\\
			&\lesssim t^{\frac{\alpha}{4}-\frac{5}{4}+\frac{5}{4p}} \|g\|_{L^q} \ \ \text{for }q\in[p^*,p]\ \text{with }\frac{1}{p}+\frac{\alpha}{5}=\frac{1}{p^*}\\
			&=t^{\frac{\alpha}{4}+\frac{5}{4}(\frac{1}{p}-\frac{1}{q})}\|f\|_{L^q}
		\end{split}
	\end{equation*}
Similarly for $p=\infty$ 
	\begin{equation*}
		\begin{split}
			\|\nabla_t^{-\alpha}f\|_{L^\infty}&=t^{\frac{\alpha}{4}-\frac{5}{4}}\|\mathcal{F}^{-1}((1+|\eta|^4)^{-\frac{\alpha}{4}}\hat{g})\|_{L^\infty} \ \ \text{where }\hat{g}(\eta)=\hat{f}(t^{-\frac{1}{4}}\eta)\\
			&\lesssim t^{\frac{\alpha}{4}-\frac{5}{4}}\|(1+|\eta|^4)^{-\frac{\alpha}{4}}\hat{g}  \|_{L^1}\\
			&\lesssim  t^{\frac{\alpha}{4}-\frac{5}{4}}\|(1+|\eta|^4)^{-\frac{\alpha}{4}}\|_{L^q}\|g\|_{L^q}  \ \ \text{for }q\in(p^*,p]\ \text{with }\frac{1}{p}+\frac{\alpha}{5}=\frac{1}{p^*}\\
			&=t^{\frac{\alpha}{4}-\frac{5}{4q}}\|f\|_{L^q}
		\end{split}
	\end{equation*}
\end{proof}

\subsection{Coifman-Meyer Operators} Most of the terms we will estimate are closely related to Coifman-Meyer operators that are either bilinear or trilinear. These operators are defined using a Fourier multiplier $m=m(\xi,\eta)$ or $m(\xi,\eta,\sigma)$, and we denote them by $T_m$:
\begin{equation*}
	\begin{split}
		T_m(f,g)&=\mathcal{F}^{-1}\int m(\xi,\eta)\hat{f}(\xi-\eta)\hat{g}(\eta)d\eta\\
		T_m(f,g,h)&=\mathcal{F}^{-1}\int m(\xi,\eta,\sigma)\hat{f}(\xi-\eta)\hat{g}(\eta-\sigma)\hat{h}(\sigma)d\sigma d\eta\\
	\end{split}
\end{equation*}
The first is for bilinear cases, while the second is trilinear. 
\begin{theorem}[Coifman-Meyer] Let $m\in L^\infty$ be a bounded function that is smooth away from the origin. Suppose that $m$ satisfies:
	\begin{equation}\label{CM_con}
		|\partial^\alpha m(\xi)|\leq \frac{C(\alpha,m)}{|\xi|^{|\alpha|}}
	\end{equation}
	for sufficiently many multiindices $\alpha$. In other words, we have:
	\begin{equation}\label{multiplier_norm}
		\begin{split}
			\|m\|_{CM}=\sup_{\xi, |\alpha|\leq N} |\xi|^{|\alpha|}|\partial^\alpha m(\xi)|
		\end{split}
	\end{equation}
	is finite with $N$ is sufficiently large.
	Then we have:
	\begin{equation*}
		T_m: L^{p_1}\times L^{p_2}\times \dots\times L^{p_n}\to L^{p}
	\end{equation*}
	is a bounded operator for $\frac{1}{p}=\sum_{j=1}^n\frac{1}{p_j}$, where $1<p_j\leq \infty $ and $0<p<\infty$. Moreover, the operator norm is less than a multiple of $\|m\|_{CM}$.
\end{theorem}
For our case, working with Coifman-Meyer operators of negative degrees is useful. More precisely, we are interested in operators with multipliers like:
\begin{equation}\label{m_neg_deg}
	m_t(\xi,\eta)=\frac{Q_k(\xi,\eta)}{(\frac{1}{t}+iZ)^n}
\end{equation}
where $Q_k$ is a homogeneous polynomial of order $k$ in $\xi$ and $\eta$ with $0\leq k\leq 4n$.
\begin{corollary}\label{CM_neg_deg}
	For multiplier $m_t$ defined as in \eqref{m_neg_deg}, we have:
	\begin{equation}\label{est_m_neg_deg}
		\|T_{m_t}(f,g)\|_{L^r}\lesssim \|\nabla_t^{k-4n}f\|_{L^{p_1}}\|g\|_{L^{q_1}}+\|f\|_{L^{p_2}}\|\nabla_t^{k-4n}g\|_{L^{q_2}}
	\end{equation}
	where $\frac{1}{r}=\frac{1}{p_1}+\frac{1}{q_1}=\frac{1}{p_2}+\frac{1}{q_2}$ with $1<p_1,p_2,q_1,q_2\leq \infty $, $0<r<\infty$.
\end{corollary}
\begin{proof}
We first decompose multiplier $m_t$ into $m_1$ and $m_2$:
\begin{equation*}
	\begin{split}
		m_t(\xi,\eta)&=m_t(\xi,\eta)\Psi_1(\xi,\eta)+m_t(\xi,\eta)\Psi_2(\xi,\eta)\\
		&=:m_1(\xi,\eta)+m_2(\xi,\eta)
	\end{split}
\end{equation*}
where $\Psi_j$'s are cut-off functions given in \eqref{cut_2}.
It suffices to show for $m_1$:
\begin{equation*}
	\begin{split}
		\mathcal{F}T_{m_1}(f,g)(\xi)&=\int m^*_1(\xi,\eta)\hat{f}(\xi-\eta)\hat{g}(\eta)d\eta\\
		&=\int \frac{Q_k(\xi,\eta)(\frac{1}{t}+|\xi-\eta|^4)^{n-\frac{k}{4}}}{(\frac{1}{t}+iZ)^n}\Psi_1(\xi,\eta)\frac{1}{(\frac{1}{t}+|\xi-\eta|^4)^{n-\frac{k}{4}}}\hat{f}(\xi-\eta)\hat{g}(\eta)d\eta\\
		&=\int \frac{Q_k(\xi,\eta)(\frac{1}{t}+|\xi-\eta|^4)^{n-\frac{k}{4}}}{(\frac{1}{t}+iZ)^n}\Psi_1(\xi,\eta)\widehat{\nabla_t^{k-4n} f}(\xi-\eta)\widehat{g}(\eta)d\eta\\
	\end{split}
\end{equation*}
Thus, it suffices to show that $\frac{Q_k(\xi,\eta)(\frac{1}{t}+|\xi-\eta|^4)^{n-\frac{k}{4}}}{(\frac{1}{t}+iZ)^n}\Psi_1(\xi,\eta)$ is a Coifman-Meyer multiplier with bounds independent from $t$. First, easily notice that $\frac{Q_k(\xi,\eta)(|\xi-\eta|^4)^{n-\frac{k}{4}}}{Z^n}\Psi_1(\xi,\eta)$ is homogeneous of degree 0 and thus automatically a Coifman-Meyer multiplier. It follows that the "smoothed" multiplier $\phi(\xi,\eta)=\frac{Q_k(\xi,\eta)(1+|\xi-\eta|^4)^{n-\frac{k}{4}}}{(1+iZ)^n}\Psi_1(\xi,\eta)$ also satisfies the Coifman-Meyer condition \eqref{CM_con}. Last, notice that 
\begin{equation*}
	\frac{Q_k(\xi,\eta)(\frac{1}{t}+|\xi-\eta|^4)^{n-\frac{k}{4}}}{(\frac{1}{t}+iZ)^n}\Psi_1(\xi,\eta)=\phi(t^{\frac{1}{4}}\xi,t^{\frac{1}{4}}\eta)
\end{equation*}
Thus, the multiplier norm of  $\frac{Q_k(\xi,\eta)(\frac{1}{t}+|\xi-\eta|^4)^{n-\frac{k}{4}}}{(\frac{1}{t}+iZ)^n}\Psi_1(\xi,\eta)$ defined in \eqref{multiplier_norm} is independent of $t$, and thus the corresponding operator's norm is independent of $t$. 
\end{proof}
\begin{remark}
%\begin{itemize}
%\item[i)] 	The cut-off functions $\Psi_1$, $\Psi_2$ ensure the new multipliers are smooth functions, especially when $k$ is odd. For example, note that the second derivative $|\xi-\eta|^3$ is not continuous. Thus, for cases where $k$ is even, we have:
	%\begin{equation*}
		%\|T_{m_t}(f,g)\|_{L^r}\lesssim \min(\|\nabla_t^{k-4n}f\|_{L^{p_1}}\|g\|_{L^{q_1}},\|f\|_{L^{p_2}}\|\nabla_t^{k-4n}g\|_{L^{q_2}})
	%\end{equation*}
   % \item[ii)] 
   In \eqref{m_neg_deg}, one can replace \(Q_k(\xi,\eta)\) with \(q_k(\xi,\eta,t)\), where 
\begin{align} \label{mdef}
q_k(\xi,\eta,t) = t^{-s} Q_{k+4l}(\xi,\eta) \left(\frac{1}{t} + iZ\right)^{-s-l}, \quad s, l \in \mathbb{N}.
\end{align}
In particular, 
\[
\frac{q_k(\xi,\eta)\left(\frac{1}{t} + |\xi-\eta|^4\right)^{n-\frac{k}{4}}}{\left(\frac{1}{t} + iZ\right)^n}\Psi_1(\xi,\eta) = \phi(t^{\frac{1}{4}}\xi, t^{\frac{1}{4}}\eta)
\]
satisfies the Coifman-Meyer condition and has a norm that is independent of \(t\). Hence, the proof of Corollary~\ref{CM_neg_deg} applies similarly. Throughout the analysis, we often use multipliers of the form $q_k $, which allow us to collectively represent several terms exhibiting the same structural features.

 %\end{itemize}
\end{remark}

\subsection{Coifman-Meyer Operators With Flag Singularities} In our estimates, we will also need trilinear terms where the multiplier has some property called the flag singularity:
\begin{definition}
	Consider multiplier $m^*=m^*(\xi,\eta,\sigma)$. We say that $m^*$ is a Coifman-Meyer multiplier with flag singularity if it admits the following type of decomposition:
	\begin{equation}\label{defm*}
		m^*(\xi,\eta,\sigma)=m_1(\xi,\eta,\sigma)m_2(\xi,\eta)m_3(\eta,\sigma)
	\end{equation}
	where all $m_j$'s are Coifman-Meyer multipliers defined in \eqref{CM_con}.
\end{definition}

Pseudo-product operators with flag singularities have been extensively studied by Muscalu \cite{Mus}, as well as Muscalu and Schlag \cite{MS}, where very general estimates were studied and proved. However, in order to discuss the terms of our concern, we need estimates that are different from their framework. Specifically, we need to apply the estimate proved in \cite{GMS2}, which is \textit{Theorem 3} in their paper.

\begin{theorem} \label{th:flagcf}
	The operator $T_{m^*}$ maps $L^{p_1}\times L^{p_2}\times L^{p_3}$ into $L^p$ bounded with operator norm depending on $\|m_j\|_{CM}$, for every $1<p_1,p_2,p_3<\infty$ with $0<p<\infty$ and  $\frac{1}{p_1}+\frac{1}{p_2}+\frac{1}{p_3}=\frac{1}{p}$. In addition, one of the $p_i$'s are allowed to be $\infty$ when $p=2$.
\end{theorem}
As for the regular pseudo-product operators, we will need flag pseudo-product operators with negative degrees. More specifically, the multiplier that we are concerned with are like:
\begin{equation}\label{m_flag_neg_deg}
	m^*_t(\xi,\eta,\sigma)= \frac{Q_k(\xi,\eta)}{(\frac{1}{t}+iZ)^n}\frac{Q_l(\eta,\sigma)}{(\frac{1}{t}+iX)^r }\frac{Q_w(\xi,\eta,\sigma)}{(\frac{1}{t}+iY)^s}
\end{equation}
with 
\begin{equation*}
	X:=\varphi(\eta,\sigma)+P(\eta,\sigma)\cdot \varphi_\sigma(\eta,\sigma)\geq \frac{1}{2}(|\eta|^4+|\sigma|^4)\gtrsim \max((|\eta-\sigma|^4+|\eta|^4), |\eta-\sigma|^4+|\sigma|^4)
\end{equation*}
\begin{corollary}\label{flag_neg_deg}
	For multiplier $m_t^*$ defined as in \eqref{m_flag_neg_deg}, we have:
	\begin{equation}
	\begin{split}
		\|T_{m_t^*}(f,g,h)\|_r\les t^{\frac{4n-k}{4}}&(\|\nabla_t^{w-4s}f\|_{p_1}\|g\|_{q_1}\|\nabla_t^{l-4r}h\|_{s_1}+\|\nabla_t^{w-4s}f\|_{p_2}\|\nabla_t^{l-4r}g\|_{q_2}\|h\|_{s_2}\\
		&+\|f\|_{p_3}\|\nabla_t^{l-4r+w-4s}g\|_{q_3}\|h\|_{s_3}+\|f\|_{p_4}\|g\|_{q_4}\|\nabla_t^{l-4r+w-4s}h\|_{s_4} )
	\end{split}
	\end{equation}
\end{corollary}
\begin{proof}
	For simplicity, we prove for the case where $n=r=s=1$ here. The other cases follow the same reasoning.
We decompose $m_t^*$ into 6 parts:
\begin{equation*}
\begin{split}
	m_t^*&=m_t^*\phi_1\Psi_1(\eta,\sigma)+m_t^*\phi_1\Psi_2(\eta,\sigma)+m_t^*\phi_2\widetilde{\Psi_1}(\xi,\eta)+m_t^*\phi_2\widetilde{\Psi_2}(\xi,\eta)+m_t^*\phi_3\widetilde{\Psi_1}(\xi,\eta)+m_t^*\phi_3\widetilde{\Psi_2}(\xi,\eta)\\&=:m_1+m_2+m_3+m_4+m_5+m_6
\end{split}
\end{equation*}
Here $\Psi_j$, $\tilde{\Psi_j}$ and $\phi_j$ are cut-off functions defined in \eqref{cut_2}, \eqref{cut_2_2} and \eqref{cut_3}.
\begin{equation*}
	\begin{split}
		&\mathcal{F}T_{m_1}(f,g,h)(\xi)\\
		&=\iint \frac{Q_k(\xi,\eta)}{\frac{1}{t}+iZ}\frac{Q_l(\eta,\sigma)}{\frac{1}{t}+iX }\frac{Q_w(\xi,\eta,\sigma)}{\frac{1}{t}+iY}\phi_1(\xi,\eta,\sigma)\Psi_1(\eta,\sigma)\widehat{f}(\xi-\eta)\widehat{g}(\eta-\sigma)\widehat{h}(\sigma)d\sigma d\eta\\
		&=\iint m^*(\xi,\eta,\sigma) \phi_1(\xi,\eta,\sigma)\Psi_1(\eta,\sigma)\widehat{\nabla_t^{k+w-8}f}(\xi-\eta)\widehat{\nabla_t^{l-4}g}(\eta-\sigma)\widehat{h}(\sigma)d\sigma d\eta \\
		&=\mathcal{F}T_{m^*}(\nabla_t^{k+w-8}f, \nabla_t^{l-4}g,h)(\xi)
	\end{split}
\end{equation*}
where 
$$
m^*(\xi,\eta,\sigma) = \frac{Q_k(\xi,\eta)(\frac{1}{t}+|\xi-\eta|^4)^{\frac{4-k}{4}}}{\frac{1}{t}+iZ}\frac{Q_l(\eta,\sigma)(\frac{1}{t}+|\eta-\sigma|^4)^{\frac{4-l}{4}}   }{\frac{1}{t}+iX }\frac{Q_w(\xi,\eta,\sigma)(\frac{1}{t}+|\xi-\eta|^4)^{\frac{4-w}{4}}  }{\frac{1}{t}+iY}.
$$

Using the same argument as for Corollary \ref{CM_neg_deg}, it is easy to verify that $m^*$ is a Coifman-Meyer multiplier with flag singularity with norm not depending on $t$. Similarly:
\begin{equation*}
	\begin{split}
		T_{m_2}(f,g,h)&=T_{m^*}(\nabla_t^{k+w-8}f,g, \nabla_t^{l-4}h)\\
		T_{m_3}(f,g,h)&=T_{m^*}(\nabla_t^{k-4}f,\nabla_t^{l+w-8}g, h)\\
		T_{m_4}(f,g,h)&=\nabla_t^{k-4}T_{m^*}(f,\nabla_t^{l+w-8}g, h)\\
		T_{m_5}(f,g,h)&=T_{m^*}(\nabla_t^{k-4}f,g, \nabla_t^{l+w-8}h)\\
		T_{m_6}(f,g,h)&=\nabla_t^{k-4}T_{m^*}(f,g, \nabla_t^{l+w-8}h)
	\end{split}
\end{equation*}
Using Lemma \ref{lem_frac_int}, we are done.
\end{proof}
\begin{remark}
	As in Corollary~\ref{CM_neg_deg}, we can replace the \( k \)-th degree polynomial \( Q_k(\xi,\eta,\sigma) \) with \( q_k(\xi,\eta,\sigma, s) \). Moreover, we define \( q_k(\xi,\eta,\sigma,s) \) in a similar manner as
\begin{align*} 
q_k(\xi,\eta,\sigma,t) = t^{-s} Q_{k+4l}(\xi,\eta,\sigma) \left(\frac{1}{t} + iY\right)^{-s-l}, \quad s, l \in \mathbb{N}.
\end{align*}
\end{remark}

\section{Dispersion of Fourth-order Schr\"{o}dinger equations} \label{sec:dispersion}
\subsection{Stationary-Phase Lemmata}
We start with the following linear estimate for the solution to fourth order linear Schr\"odinger equation. 

\begin{lemma} \label{lem:linear} One has in $\R^5$
\begin{align}
e^{it\Delta^2}f(x)= \rho(x,t) \hat{f}(\xi_*) + t^{-\frac{11}{8}} \| \hat{f}\|_{H^3}
\end{align}
where $\sup_{x}|\rho(x,t)| \lesssim t^{-\f 54} $ and 
\[\xi_* = \Big( - \frac{|x|}{ 4 t} \Big)^{\f13} \cdot \frac{x}{|x|}.\]
In particular, 
\[ \|e^{it\Delta^2} f\|_{L^{\infty}} \les t^{-\f54} \|\hat{f}\|_{L^{\infty}}+ t^{-\frac{11}{8}}\|\hat{f}\|_{H^3}
\]
\end{lemma}
\begin{proof}
We start writing 
\begin{align}\label{mainosc}
e^{it \Delta^2}f(x) = \int_{\R^5} e^{i(|\xi|^4 + \frac{\xi \cdot x}{t})} \hat{f}(\xi) d \xi 
\end{align}
Let $\phi_{t,x}(\xi):=|\xi|^4 + \frac{\xi \cdot x}{t}$. One can easily see that $\nabla \phi_{t,x} (\xi_*) = 0 $ and therefore the oscillatory integral in $\eqref{mainosc}$ has a  critical point at $\xi_*$. 
Moreover, for any $v\in \R^5$  the following lower bound holds:
$$ \la \bf{H}_{\phi_{t,x}}v,v \ra \geq |\xi|^2 \la v, v \ra $$
enabling us to derive 
\begin{align} \label{cpbound}
|\nabla \phi_{t,x}(\xi)| \gtrsim |\xi - \xi_*| |\xi|^2.
\end{align}
To estimate \eqref{mainosc}, we define $\chi(\xi)$ as a smooth cut-off function which is one when $|\xi|<1$ and zero outside of $|\xi|\geq 2 $ and decompose 
\eqref{mainosc} as follows:
\begin{align}
e^{it \Delta^2}f(x)= \int_{\R^5} e^{it \phi_{t,x}(\xi)}  \chi_{\xi_*}(\xi) \hat{f}(\xi) \: d \xi + \int_{\R^5} e^{i t \phi_{t,x}(\xi)} \tilde{\chi}_{\xi_*}(\xi) \hat{f}(\xi) \: d \xi = I_1(x,t)+ I_2(x,t)
\end{align}
where $\chi_{\xi_*}(\xi)= \chi \big(|\xi|/|\xi-\xi_*|)$.

We further decompose $I_1(x,t)$ as 
\begin{align} \label{mtdec}
I_1(x,t) & = \hat{f}(\xi_*) \int_{\R^5} e^{it \phi_{t,x}(\xi)}  \chi_{\xi_*}(\xi) \: d \xi + \int_{\R^5} e^{it \phi_{t,x}(\xi)}  \chi_{\xi_*}(\xi)[ \hat{f}(\xi) - \hat{f}(\xi_*)] \: d \xi \\ & =: \hat{f}(\xi_*) \rho(x,t) + I_3(x,t) \nn 
\end{align}
We first establish $|\rho(x,t)| \lesssim t^{-\f54}$. To do that we write 
\begin{align*}
\rho(x,t)= \int_{\R^5} e^{it \phi_{t,x}(\xi)}  \chi_{\xi_*}(\xi) \chi(t^{\f14} |\xi-\xi_*|)  \: d \xi +\int_{\R^5} e^{it \phi_{t,x}(\xi)}  \chi_{\xi_*}(\xi) \tilde{\chi}(t^{\f14} |\xi-\xi_*|) \: d \xi
\end{align*}
The first term in $\rho(x,t)$ is bounded by $t^{-\f54}$. Using 
\begin{align} \label{ibp}
e^{it \phi_{t,x}(\xi)} = \frac{\nabla( e^{it \phi_{t,x}(\xi)}) \cdot \nabla\phi_{t,x}(\xi)}{it |\nabla\phi_{t,x}(\xi)|^2}
\end{align}
we apply integration by parts twice, to bound the second term by 
\begin{align}\label{rho2}
t^{-2}\int_{\R^5}  \Bigg| \nabla \cdot \Bigg(\frac{\nabla\phi_{t,x}(\xi)}{|\nabla\phi_{t,x}(\xi)|^2} \nabla \cdot  \Big( \frac{ \nabla\phi_{t,x}(\xi) \chi_{\xi_*}(\xi) \tilde{\chi}(t^{\f14} |\xi-\xi_*|) }{|\nabla\phi_{t,x}(\xi)|^2 } \Big) \Bigg) \Bigg| d \xi 
\end{align}
Using the bounds on $| \nabla \phi_{t,x}(\xi)|$, we have the following estimates  
\begin{align}
 \bigg | \nabla \cdot \bigg( \frac{ \nabla \phi_{t,x}(\xi)}{|\nabla \phi_{t,x}(\xi)|^2} \bigg) \bigg| \les \frac{ | \Delta \phi_{t,x}(\xi)| } {|\nabla \phi_{t,x}(\xi)|^2} & \les \frac{1}{| \xi - \xi_*|^2 |\xi|^2} \label{ibp1} \\ 
  \bigg |\nabla \cdot \bigg( \frac{ \nabla \phi_{t,x}(\xi)}{|\nabla \phi_{t,x}(\xi)|^2} \nabla \cdot \bigg(\frac{ \nabla \phi_{t,x}(\xi)}{|\nabla \phi_{t,x}(\xi)|^2} \bigg) \bigg) \bigg | &\les \frac{ | \Delta \phi_{t,x}(\xi)|^2 } {|\nabla \phi_{t,x}(\xi)|^4} + \frac{ |\nabla \Delta \phi_{t,x}(\xi)| }{|\nabla \phi_{t,x}(\xi)|^3 } \nn \\ 
& \les \frac{1}{|\xi - \xi_*|^4 |\xi|^4} + \frac{1}{|\xi - \xi_*|^3 |\xi|^5}\label{ibp2}
\end{align}
Moreover, if the derivative acts on any of the cutoff functions, it introduces a factor of $|\xi-\xi|^{-1}$. Therefore, using the support of $\chi_{\xi_*}(\xi)$  we obtain
\begin{align*}
|\eqref{rho2}| \lesssim t^{-2} \int_{|\xi-\xi_*| \gtrsim t^{-\f14}} |\xi-\xi_*|^{-8} \: d\xi \lesssim t^{-\f54} \nn
\end{align*}

We next bound $I_3(x,t)$. Noting that we removed the boundary term, we apply one time integration by parts to find, 
\begin{align*}
I_3(x,t)= \frac{1}{it} \int_{\R^5}e^{it \phi_{t,x}(\xi)} F(\xi,x,t) \chi(t^{\f14} |\xi-\xi_*|) \:d\xi +\frac{1}{it}   \int_{\R^5}e^{it \phi_{t,x}(\xi)} F(\xi,x,t) \tilde{\chi}(t^{\f14} |\xi-\xi_*|) \:d\xi
\end{align*}
where $F(\xi,x,t) := \nabla \cdot \Big( \frac{\nabla \phi_{t,x}(\xi) \chi_{\xi_*}(\xi) [\hat{f}(\xi)- \hat{f}(\xi_*)]}{|\nabla \phi_{t,x}(\xi)|^2} \Big) $. Using, \eqref{ibp1}, we can estimate 
\begin{align}
|F(\xi,x,t)| \lesssim \frac{|\hat{f}(\xi)- \hat{f}(\xi_*)|}{|\xi- \xi_*|^4} + \frac{|\nabla \hat{f}|}{|\xi| |\xi- \xi_*|^2}. 
\end{align}
We have 
\begin{align}\label{I31a}
 \int_{\R^5}\frac{|\hat{f}(\xi)- \hat{f}(\xi_*)|}{|\xi- \xi_*|^4} \chi(t^{\f14} |\xi-\xi_*|) d\xi  \lesssim \sup_{\xi} \Bigg| \frac{\hat{f}(\xi)- \hat{f}(\xi_*)}{|\xi- \xi_*|^{\f12}} \Bigg| \int_{|\xi-\xi_*|\lesssim t^{-\f14}} \frac{d\xi}{|\xi - \xi_*|^{\f72}} \les t^{-\frac{3}{8}} \| \hat{f} \| _{H^3 }
\end{align}
where in the last inequality we used Sobolev embedding. 

On the other hand, 
\begin{align} \label{I31b}
 \int_{\R^5}\frac{|\nabla \hat{f}|}{ |\xi| |\xi- \xi_*|^2} \chi(t^{\f14} |\xi-\xi_*|)  d\xi  \lesssim \Bigg\| \frac{\nabla \hat{f}}{\cdot} \Bigg\|_{L^{\frac{10}{3}}} \|\chi(t^{\f14} |\cdot-\xi_*|) |\cdot -\xi_*|^{-2}\|_{L^{\frac{10}{7}}} \lesssim t^{-\frac{3}{8}} \| \hat{f} \| _{H^3}
\end{align}
where in the last inequality, we first used Hardy's inequality and then Sobolev embedding to bound 
\begin{align}\label{HS}
\Bigg\| \frac{\nabla \hat{f}}{\cdot} \Bigg\|_{L^{\frac{10}{3}}} \les \| \hat{f} \|_{W^{2,\frac{10}{3}}} \les \|\hat{f}\|_{H^3}.
\end{align}
Applying  another  integration by parts to the second term in  $I_3(x,t)$, we proceed to estimate the remaining expression:  
\begin{align}
\frac{1}{t^2} \int_{\R^5} \Bigg| \nabla \cdot \Bigg[ F(\xi,x,t) \tilde{\chi}(t^{\f14} |\xi-\xi_*|)    \frac{\nabla \phi_{t,x}(\xi)}{|\nabla \phi_{t,x}(\xi)|^2  }\Bigg] \Bigg| \: d\xi 
\end{align}

One can compute 
\begin{align}
\Bigg| \nabla \cdot \Bigg[ F(\xi,x,t) \tilde{\chi}(t^{\f14} |\xi-\xi_*|)    \frac{\nabla \phi_{t,x}(\xi)}{|\nabla_{\xi} \phi_{t,x}(\xi)|^2  }\Bigg] \Bigg| \les \frac{|\hat{f}(\xi)- \hat{f}(\xi_*)|}{|\xi- \xi_*|^8} + \frac{|\nabla \hat{f}|}{|\xi||\xi- \xi_*|^6} + \frac{|\Delta \hat{f}|}{ |\xi| |\xi- \xi_*|^5}
\end{align}
In a similar argument that is applied in \eqref{I31a} and \eqref{I31b}, we have 
\begin{multline} \label{I32a}
t^{-2} \int_{\R^5} \Bigg[\frac{|\hat{f}(\xi)- \hat{f}(\xi_*)|}{|\xi- \xi_*|^8} + \frac{|\nabla \hat{f}|}{|\xi||\xi- \xi_*|^6} \Bigg]\tilde{\chi}(t^{\f14} |\xi-\xi_*|) \: d \xi  \\ \les \frac{\|\hat{f}\|_{W^{3,H^32}}}{t^2} \Big(
\|\tilde{\chi}(t^{\f14} |\cdot-\xi_*|) |\cdot -\xi_*|^{-\frac{7}{2}}\|_{L^{1}}+\|\tilde{\chi}(t^{\f14} |\cdot-\xi_*|) |\cdot -\xi_*|^{-6}\|_{L^{\frac{10}{7}}}\Big)\\
\les 
t^{-\frac{11}{8}} \| \hat{f} \|_{H^3}. 
\end{multline}
Moreover, we have by Hardy's inequality
\begin{align}\label{I32b}
t^{-2} \int_{\R^5} \frac{|\Delta \hat{f}|\tilde{\chi}(t^{\f14} |\xi-\xi_*|)}{ |\xi| |\xi- \xi_*|^5} \les \Bigg\| \frac{\Delta \hat{f}}{\cdot} \Bigg\|_{L^2} \| \tilde{\chi}(t^{\f14}|\cdot-\xi_*|) |\cdot-\xi_*|^{-5} \|_{L^2} \les t^{-\frac{11}{8}} \| \hat{f} \|_{H^3}.  
\end{align}

Combining all the estimates from \eqref{I31a}, \eqref{I31b}, \eqref{I32a} and \eqref{I32b} we obtain the estimate $|I_3(x,t)| \les t^{-\frac{11}{8}} \| \hat{f} \|_{H^3} $. 

In the support of $\tilde{\chi}_{\xi_*}$ one has $|\xi| \les |\xi -\xi_*|$. Consequently, the estimates for $I_2(x,t)$ can be derived in a similar manner using the cut-off $\chi(t^{\f14}|\xi|)$ rather than $\chi(t^{\f14}|\xi-\xi_*|)$. 

\end{proof}

Note that the plane wave solutions to the fourth order free quation is $e^{i x \cdot( \xi |\xi|^2 + t)}$, indicating that particles with large energies travel so much faster than the particles with small energies.  This characteristic is, in fact, crucial for many of our estimates. The following lemma provides a quantitative expression of this behavior. We define
$$f_j(x):= [\varphi_j(\xi) \hat{f}(\xi)]^{\vee}(x). $$

\begin{lemma}\label{lem:ben} Let $t>1$ and $1\leq p \leq 2$ we have 
 \begin{align} \label{ineq:ben}
\| e^{it \Delta^2} f_j\|_{L^{\frac{p}{p-1}}} \les 2^{- \frac{5(2-p)j}{p} } t^{-\frac{5(2-p)}{2p}} \|f_j\|_{L^p}
\end{align}

\end{lemma}
 \begin{proof}
We have the following expression 
\begin{align}
e^{it \Delta^2}f_j= [e^{it |\xi|^4} \varphi_j(\xi) \hat{f}(\xi) ]^{\vee} =: [K \ast f](x)
\end{align}
where 
$$ K(x) = \int_{\R^5} e^{it \phi_{t,x}(\xi)} \varphi_{j}(\xi) \: d \xi. $$
Given these settings, to prove the statement, we need to derive the following bound:
\begin{align}
\text{sup}_{x} |K(x)| \les  2^{-5j } t^{- \frac{5}{2}} \|f_j\|_{L^1}. 
\end{align}
Interpolation this inequality with $L^2$ conservation conclude \eqref{ineq:ben}. 

Applying twice integration by parts to $K(x)$ and decomposing, we obtain 
\begin{multline}\label{boundK}
 K(x)= \frac{1}{t^2} \int_{\R^5}e^{it \phi_{t,x}(\xi)} K_1(\xi:x,t) \chi(t^{\f12}2^{j}|\xi-\xi_*|)\: d \xi \\
 + \frac{1}{t^2} \int_{\R^5}e^{it \phi_{t,x}(\xi)} K_1(\xi:x,t)  \tilde{\chi}(t^{\f12}2^{j}|\xi-\xi_*|) \: d \xi 
\end{multline}
where 
$$K_1(\xi:x,t) = \nabla \cdot 
 \Bigg[\frac{\nabla \phi_{t,x}(\xi) }{|\nabla \phi_{t,x}(\xi)|^2} \nabla \cdot \Big[ \frac{\nabla \phi_{t,x}(\xi) \varphi_j (\xi)}{|\nabla \phi_{t,x}(\xi)|^2} \Big] \Bigg] $$
 and using \eqref{ibp1}, \eqref{ibp2} has the following bounds 
\begin{align*}
&|K_1(\xi:x,t)| \les \frac{2^{-4j}}{|\xi-\xi_*|^{2}} \max\{ 2^{-2j}, |\xi-\xi_*|^{-2}\} \varphi_j(\xi), \\
& \Bigg|\nabla \cdot 
 \Bigg[\frac{\nabla \phi_{t,x}(\xi) K_1(\xi:x,t) }{|\nabla \phi_{t,x}(\xi)|^2}  \Bigg]\Bigg| \les \frac{2^{-9j}\varphi_j(\xi)}{|\xi-\xi_*|^3}+  \frac{2^{-6j}\varphi_j(\xi}{|\xi-\xi_*|^6} . 
\end{align*}
Using the first bound, we can estimate the first term in \eqref{boundK} as following
\begin{align*}
t^{-2}\Bigg|\int_{\R^5}e^{it \phi_{t,x}(\xi)} K_1(\xi:x,t) \chi(t^{\f12}2^{j}|\xi-\xi_*|)\: d \xi \Bigg| & \les 2^{-4j} \int_{|\xi-\xi_*|  \les t^{-\f12}2^{-j}}  |K_1(\xi:x,t)| \: d \xi \\ 
& \les t^{-2} 2^{-4j} \max\{2^{-5j}t^{-\f32}, 2^{-j}t^{-\f12} \} \les t^{-\f52} 2^{-5j}
\end{align*}
Moreover, applying another integration by parts to the second term in \eqref{boundK}, we have the estimate 
\begin{align*}
t^{-2}\Bigg|\int_{\R^5}e^{it \phi_{t,x}(\xi)} & K_1(\xi:x,t) \tilde{\chi}(t^{\f12}2^{j}|\xi-\xi_*|)\: d \xi \Bigg|  \\ 
&\les t^{-3} \int_{|\xi-\xi_*|  \gtrsim t^{-\f12}2^{-j}}\Bigg|\nabla \cdot 
 \Bigg[\frac{\nabla \phi_{t,x}(\xi) K_1(\xi:x,t) }{|\nabla \phi_{t,x}(\xi)|^2}  \Bigg]\Bigg|  \varphi_j(\xi) 
 \: d \xi \nn   \\
& \les t^{-3}2^{-6j} \int_{|\xi-\xi_*|  \gtrsim t^{-\f12}2^{-j}} \frac{d \xi}{|\xi-\xi_*|^{6}} + t^{-3}2^{-12j} \int_{|\xi|\sim 2^j } \varphi_j(\xi) \: d\xi \les t^{-\f52} 2^{-5j}
 \end{align*}
 where in the second inequality we used the inequality $ 2^{-3j} |\xi-\xi_*|^{-3} \les 2^{-6j} + |\xi-\xi_*|^{-6}$. 

  \end{proof}
\subsection{Estimate of \( \|\hat{B}(f,f)\|_{L^{\infty}_xL^{\infty}_t} \)}

We start with the following oscillatory lemma. 
\begin{lemma}\label{lemsupB} We have 
\begin{align}
\bigg| \int_{\R^5} e^{-it(|\eta|^4+|\xi-\eta|^4)} \hat{f}(\xi - \eta) \hat{f}(\eta) \: d \eta \bigg| \les t^{-\f 98} \|f\|_{X}^2. 
\end{align}
\end{lemma}
\begin{proof}
Let $\phi_{\xi}(\eta) = |\eta|^4+|\xi-\eta|^4$, we analyze the following oscillatory integral.  
\begin{align}\label{oscbff}
 \int_{\R^5} e^{-it\phi_{\xi}(\eta)} \hat{f}(\xi - \eta) \hat{f}(\eta) \: d \eta 
\end{align}
We note that $ \nabla\phi_{\xi}(\eta) = 4\eta|\eta|^2-4(\xi-\eta)|\xi-\eta|^2$, and hence the critical point of the integral arises at $ \eta = \xi/2$. Moreover, one has 
$$
\la \bf{H}_{\phi_{\xi}}v,v \ra \geq (|\eta|^2 +|\xi-\eta|^2) \la v, v \ra 
$$
allowing us to derive 
\begin{align} 
|\nabla \phi_{\xi}(\eta)| \gtrsim |\eta -\xi/2| (|\eta|^2 +|\xi-\eta|^2).
\end{align}
Note that, one has 
$$
|\eta|^2+|\xi-\eta|^2 - |\eta -\xi/2|^2 = |\eta|^2 - \xi \cdot \eta + \frac{3|\xi^2|}{4} = |\eta -\xi/2|^2 + \frac{|\xi^2|}{4} \geq 0 
$$
and therefore,  $|\eta|^2+|\xi-\eta|^2 \geq |\eta -\xi/2|^2$.  

We again start with decomposing the integral in \eqref{oscbff} as in the following 
\begin{align*}
\int_{\R^5} e^{-it\phi_{\xi}(\eta)} \hat{f}(\xi - \eta) \hat{f}(\eta) \chi(t^{\f18}|\eta-\xi/2) \: d \eta +  \int_{\R^5} e^{-it\phi_{\xi}(\eta)} \hat{f}(\xi - \eta) \hat{f}(\eta) \tilde{\chi}(t^{\f18} |\eta - \xi/2|) \: d\eta = I + II 
\end{align*}
We have  $|I| \les t^{-\f98} \| f\|_{X}^2$. We need to apply twice integration by parts to $II$ to find 
\begin{align} \label{IIbound}
|II|\leq \frac{1}{t^2} \int_{\R^5} \bigg| \nabla \cdot \bigg(
\frac{\nabla \phi_\xi (\eta)}{|\nabla \phi_\xi (\eta)|^2} \nabla \cdot \Big(\frac{\nabla \phi_\xi (\eta)}{|\nabla \phi_\xi (\eta)|^2} \hat{f}(\xi - \eta) \hat{f}(\eta) \tilde{\chi}(t^{\f18} |\eta - \xi/2|) \Big) \bigg) \bigg| \: d\eta.
\end{align}
Using \eqref{ibp1} and \eqref{ibp2}, the integrand in \eqref{IIbound} can be bounded by the sum of the following terms 
\begin{align}\label{intbounds}
&\frac{\hat{f}(\xi - \eta) \hat{f}(\eta)}{|\eta - \xi/2|^4(|\eta|^2 + |\eta - \xi|^2)^2} + \frac{\hat{f}(\xi - \eta) \hat{f}(\eta)}{|\eta - \xi/2|^2(|\eta|^2 + |\eta - \xi|^2)^3} \les \frac{\hat{f}(\xi - \eta) \hat{f}(\eta)}{|\eta - \xi/2|^8}, \\
&\frac{\nabla(\hat{f}(\xi - \eta) \hat{f}(\eta)) }{|\eta - \xi/2|^3(|\eta|^2 + |\eta - \xi|^2)^2 } \les  \frac{\hat{f}(\xi - \eta) \: \nabla \hat{f}(\eta)+\hat{f}(\eta) \: \nabla \hat{f}(\xi - \eta)}{|\eta - \xi/2|^5 |\eta - \xi| |\eta| } \nn \\
&\frac{ \Delta(\hat{f}(\xi - \eta)\hat{f}(\eta))}{|\eta - \xi/2|^2(|\eta|^2 + |\eta - \xi|^2)^2} \les \frac{\hat{f}(\xi - \eta) \: \Delta\hat{f}(\eta)}{|\eta - \xi/2|^5 |\eta-\xi|} + \frac{\hat{f}(\eta)\: \Delta\hat{f}(\eta- \xi)}{|\eta - \xi/2|^5 |\eta|} + \frac{\nabla \hat{f}(\eta) \:\nabla \hat{f}(\xi -\eta) }{|\eta - \xi/2|^5 |\eta- \xi| }. \nn 
\end{align}
The estimates now follow by by Cauchy-Schwarz. In particular, we have 
\begin{align}
t^{-2}\int_{\R^5} \frac{\hat{f}(\xi - \eta) \hat{f}(\eta)}{|\eta - \xi/2|^8} \tilde{\chi}(t^{\f18} |\eta - \xi/2|) \: d \eta \les t^{-2}\|\hat{f}\|^2_{L^{\infty}} \int_{|\eta - \xi/2|\gtrsim t^{-\f18}} |\eta - \xi/2|^{-8} \les  t^{-\frac{13}{8}}  \|\hat{f}\|^2_{X}. 
\end{align}
Moreover, using also Hardy's inequality we have
\begin{multline*}
t^{-2} \int_{\R^5} \tilde{\chi}(t^{\f18} |\eta-\xi/2|) \frac{\hat{f}(\xi - \eta) \: \nabla \hat{f}(\eta)+\hat{f}(\eta) \: \nabla \hat{f}(\xi - \eta)}{|\eta - \xi/2|^5 |\eta - \xi| |\eta| } \: d \eta \\\les t^{-2} \sup_{\eta} \big\{\chi(t^{\f18} |\eta-\xi/2|)|\eta-\xi/2|^{-5}\big\} \Big\| \frac{ \hat{f}}{\cdot}\Big\|_{L^2} \Big\| \frac{\nabla \hat{f}}{\cdot}\Big\|_{L^2} \les t^{-\f98} \|f\|^2_{X}
\end{multline*}
Similarly
\begin{multline*}
t^{-2} \int_{\R^5} \tilde{\chi}(t^{\f18} |\eta-\xi/2|) \Bigg( \frac{\hat{f}(\eta) \: \Delta\hat{f}(\eta)}{|\eta - \xi/2|^5 |\eta|} + \frac{\hat{f}(\xi-\eta)\: \Delta\hat{f}(\eta- \xi)}{|\eta - \xi/2|^5 |\eta- \xi|}\Bigg) \: d \eta \\ \les t^{-2} \sup_{\eta} \big\{\chi(t^{\f18} |\eta-\xi/2|)|\eta-\xi/2|^{-5}\big\} \| \nabla \hat{f}\|_{L^2} \| \Delta \hat{f}\|_{L^2} \les t^{-\f98} \|f\|^2_{X}.  
\end{multline*}
Finally we estimate contribution of the last term in \eqref{intbounds} to \eqref{IIbound} as 
\begin{align*}
t^{-2} \int_{\R^5} \tilde{\chi}(t^{\f18}|\eta - \xi/2|)  \frac{\nabla \hat{f}(\eta) \:\nabla \hat{f}(\xi -\eta) }{|\eta - \xi/2|^5 |\eta- \xi|  } \: d \eta \les t^{-\frac{11}{8}} \|\nabla{\hat{f}}\|_{L^2}  \|\Delta\hat{f}\|_{L^2} \les t^{-\f98} \|f\|^2_{X}.
\end{align*}

This establishes the statement. 

\end{proof}
The following Corollary gives the first estimate for the contraction in $X$-norm. 
\begin{corollary}\label{corsupB} We have \( \| \hat{B}(f,f)\|_{L^\infty_t L^{\infty}_x} \les \|f\|_X^2 \)
    \end{corollary}
\begin{proof} Using the equality \eqref{contraction}, we have 
\[  | \hat{B}(f,f)| \les \int_1^{t} s^{-\f 98} \|f\|^2_X \: ds \les \|f\|^2_X \] 
where we also used Lemma~\ref{lemsupB} in the first inequality. 
    \end{proof}

\section{Estimates on $g$}\label{sec:gbounds}
As introduced in Section~\ref{sec:setup}, the function \( g \) is defined as 
\[
\hat{g}(\xi,t) = i \int_{\mathbb{R}^5} \frac{1}{\frac{1}{t} + iZ} e^{it\varphi} \widehat{f}(\xi-\eta,t) \hat{f}(\eta,t) \, d\eta
\]
where we suppress the variables in $Z(\xi,\eta)=\varphi(\xi,\eta) + P(\xi,\eta) \cdot \partial_{\eta}\varphi(\xi,\eta)$. By letting \( A := A_t = \frac{1}{t} + iZ \), \( g(x,t) \) can be expressed as 
\begin{align}\label{gexp}
    g(x,t) = i e^{it\Delta^2} \big[T_{\frac{1}{A}}(u, u)(x,t)\big],
\end{align}
where \( T_{\frac{1}{A}}(u, u) \) denotes a bilinear Coifman-Meyer operator with multiplier $\frac{1}{A}$. We also observe that the form of $g$ naturally leads to a fast decay in $L^\infty$ space:
\begin{lemma}\label{lem:decay_g}
	We have
	\begin{equation}
		\|e^{-it\Delta^2}g\|_{L^\infty}\lesssim t^{-\frac{3}{2}}\|f\|_X^2
	\end{equation}
\end{lemma} 
\begin{proof}
Observe that we can write $\hat{g}$ by:
\begin{equation*}
	\widehat{\nabla_t^{4}g}(\xi,t)=i \int_{\mathbb{R}^5} \frac{\frac{1}{t}+|\xi|^4}{\frac{1}{t} + iZ} e^{it\varphi} \widehat{f}(\xi-\eta,t) \hat{f}(\eta,t) \, d\eta
\end{equation*}
By a similar argument as in Lemma \ref{CM_neg_deg}, we have that $\frac{\frac{1}{t}+|\xi|^4}{A}$ is a Coifman-Meyer multiplier with bound independent of $t$. 
	Applying Lemma \ref{lem_frac_int}, we have:
	\begin{equation*}
		\begin{split}
			\|e^{-it\Delta^2}g\|_{L^\infty}&=\|\nabla_t^{-4}e^{-it\Delta^2}\nabla_t^4g\|_{L^\infty}=\|\nabla_t^{-4}(T_{\frac{\frac{1}{t}+|\xi|^4}{A}}(u,u))\|_{L^\infty}\\
			&\lesssim t^{1-\frac{5}{4q}}\|T_{\frac{\frac{1}{t}+|\xi|^4}{A}}(u,u)\|_{L^q}\lesssim t^{ 1-\frac{5}{4q}  }\|u\|_{L^q}\|u\|_{L^\infty}\\
			&\lesssim t^{-\frac{3}{2}+\frac{5}{4q}}\|f\|_{X}^2
		\end{split}
	\end{equation*}
	for any $2<q<\infty$, which implies \( \|e^{-it\Delta^2}g\|_{L^\infty}\lesssim t^{-\frac{3}{2}+}\|f\|_{X}^2 \) .
\end{proof}
Despite this fast-decay property of $g$, it is not very 'localized'. More specifically, we will rigorously establish the required estimates for \( g(x,t) \) as listed in \eqref{est_g}, in this section, assuming the bootstrap conditions in \eqref{Boot_est_g} and \eqref{Boot_est_h}. We start with $L^2$ bound. 
%\subsection{Estimates of $\|g\|_{L^2}$}

\begin{lemma}\label{lem:L2_g}
Assuming \eqref{Boot_est_g} and \eqref{Boot_est_h}, we have
\begin{equation}\label{L2_g}
	\|g\|_{L^2_x}\lesssim \|f\|_{X}^2. 
\end{equation}
\end{lemma}
\begin{proof} Using Lemma~\ref{CM_neg_deg} alongside Lemma~\ref{lem_frac_int}, we have
\begin{align}\label{TAbound}
\|T_{\frac{1}{A}}(u,u) \|_{L^2}\lesssim\|\nabla_t^{-4}u\|_{L^2_x}\|u\|_{L^\infty_x}\lesssim t^{-\frac{1}{4}}\|u\|_{L^2_x}\|t^{\f54}u\|_{L^\infty_x} \les t^{-\f 14} \|f\|_{X}^2. 
	\end{align}

Therefore, using the representation in \eqref{gexp} we obtain the statement for $t>0$. 
\end{proof}
%\subsection{Estimates of $\|xg\|_{L^2}$}
\begin{lemma}
Assuming \eqref{Boot_est_g} and \eqref{Boot_est_h}, we have
	\begin{equation}\label{Decay_xg}
		\|xg\|_{L^2_x}\lesssim t^{-\f14} \|f\|_{X}^2
	\end{equation}
\end{lemma}
\begin{proof}
In the proof of this lemma, we derive several bounds in a general form so that these estimates can later be applied to control other terms in subsequent lemmas.

First, note that by applying Lemma~\ref{lem_frac_int} we obtain
\begin{align*}
\|\nabla_t^{-(4-k)} e^{-it\Delta^2} f\|_{L^2} \lesssim \|f\|_{L^{\frac{10}{13-2k}}} \lesssim \|\langle x \rangle^{(4-k)+} f\|_{L^2}, \quad 2 \le k \le 3,
\end{align*}
where in the first inequality we used the weighted estimate
\begin{align}\label{weighted}
\|f\|_{L^p} \lesssim \|\langle x \rangle^{\alpha} f\|_{L^2}, \quad 1 \le p \le 2,\quad \alpha > 5\left(\frac{1}{p} - \frac{1}{2}\right).
\end{align}
Consequently, by Corollary~\ref{CM_neg_deg}, we obtain for \(2 \le k \le 3\)
\begin{align}
\|T_{\frac{q_{k}}{A}}(u,u)\|_{L^2} \lesssim \|\langle x \rangle^{(4-k)+} f\|_{L^2}\|u\|_{L^{\infty}}   \lesssim t^{-\frac{5}{4}} \|\langle x \rangle^{(4-k)+} f\|_{L^2}\|f\|_{X}, \label{polyb}
\end{align}
where the second inequality follows from Lemma~\ref{CM_neg_deg} and  we used Lemma~\ref{lem_frac_int} along with the fact that  \(\|t^{\f54} u\|_{L^{\infty}_x} \les \|f\|_{X}\) in the third inequality.
Second, we have for any $0\leq k \leq 4$
\begin{align}
\| e^{it\Delta^2} T_{\frac{q_k}{A}}(e^{it\Delta^2}(xf),u) \|_{L^2_x} & \les \|\nabla_t^{-(4-k)}(e^{it\Delta^2}(xf))\|_{L^2_x}\|u\|_{L^\infty_x}  + \|(e^{it\Delta^2}(xf))\|_{L^2_x}\|\nabla_t^{-(4-k)}u\|_{L^\infty_x} \nn \\
& \les t^{\frac{4-k}{4}} \|xf\|_{L^2} \|u\|_{L^{\infty}}\les t^{-\f{k+1}4} \|f\|_{X}^2.\label{fracwuse}
\end{align}

Finally, we can bound $\|xg\|_{L^2_x}$. Note that using the identity $
\widehat{xg}(\xi,t) = -i\,\partial_{\xi}\hat{g}(\xi,t)$, 
we can express
\begin{align}\label{eq:xg}
\widehat{xg}(\xi,t) = t\int_{\R^5} \frac{q_3(\xi,\eta,t)}{A}\,e^{it\varphi}\,\hat{f}(\xi-\eta)\,\hat{f}(\eta)\,d\eta
+ \int_{\R^5} \frac{e^{it\varphi}}{A}\,\widehat{xf}(\xi-\eta)\,\hat{f}(\eta)\,d\eta.
\end{align}
where $q_k$ is defined as in \eqref{mdef}. Using \eqref{polyb} with $k=3$, and \eqref{fracwuse} for $k=0$, we obtain 

$$\|\eqref{eq:xg}\|_{L^2_x} \les t \| e^{it\Delta^2} T_{\frac{q_3}{A}}(u,u) \|_{L^2_x} + \| e^{it\Delta^2} T_{\frac{q_0}{A}}(e^{it\Delta^2}(xf),u) \|_{L^2_x} \les t^{-\f14} \|f\|_X^2
$$

\end{proof}

%\subsection{Estimates of $\|x^2g\|_{L^2}$}
\begin{lemma}
	Assuming \eqref{Boot_est_g} and \eqref{Boot_est_h}, we have
	\begin{equation}\label{Decay_x2g}
		\|x^2g\|_{L^2}\lesssim \|f\|_X^2+\|f\|_X^3
	\end{equation}
\end{lemma}
\begin{proof} We start computing 
\begin{align}
	\label{x2g_1}\widehat{x^2g}(\xi,t)&=t^2\int_{\R^5}\frac{Q_6}{A}e^{it\varphi}\hat{f}(\xi-\eta)\hat{f}(\eta)d\eta+t\int_{\R^5}\frac{q_2(\xi,\eta,t)}{A}e^{it\varphi}\widehat{f}(\xi-\eta)\hat{f}(\eta)d\eta \\ \label{x2g_2}
	&+ t \int_{\R^5}\frac{q_3(\xi,\eta,t)}{A}e^{it\varphi}\widehat{xf}(\xi-\eta)\hat{f}(\eta)d\eta +\int_{\R^5} \frac{1}{A}e^{it\varphi}\widehat{x^2f}(\xi-\eta)\hat{f}(\eta)d\eta. 
\end{align}
 Using \eqref{fracwuse} with $k=0$ and $k=3$ , we obtain $ \|\eqref{x2g_2}\|_{L^2} \les\|f\|_X^2$. Also \eqref{polyb} with $k=2$ gives the bound $t^{-\f14+} \|f\|_{X}^2 $ for the $L^2$ norm of the  second term in \eqref{x2g_1}. 

We therefore, focus  on the first term in \eqref{x2g_1} which requires a more detailed analysis due to the fast growth in time. To handle this term, we apply the strategy illustrated in \eqref{Decom_Bi}.
We first define the operator
\begin{align*}
H_{6}(\widehat{f},\widehat{g})(\xi,t):=\int_{\R^3}\frac{Q_6}{A}e^{it\varphi}\hat{f}(\xi-\eta)\hat{g}(\eta)d\eta.  
\end{align*}
Now, recall that by  \eqref{Decom_Bi} we have the decomposition: $f=f_1+f_*+g+h$. Therefore, to bound the term $\|\eqref{x2g_1}\|_2$, we study  
\begin{align} \label{terms6}
H_{6}(\widehat{f},\widehat{f_1} +\widehat{f_*})(\xi,t), \,\,\  H_{6}(\widehat{f},\widehat{g})(\xi,t), \,\,\,\ H_{6}(\widehat{h},\widehat{h})(\xi,t), \,\,\,\ H_{6}(\widehat{g},\widehat{h})(\xi,t).
\end{align}

In particular, by Plancherel, it is enough to bound the $L^2$ norm of the each term in  \eqref{terms6} by $t^{-2}$.  We note that since $\|f_1+f_*\|_X$ is small independently of $t$, the estimate of $\|H_{6}(\widehat{f},\hat{f_1} +\hat{f_*})(\xi,t)\|_2$ naturally follows from the strategies applied to estimate $\|H_{6}(\widehat{f},\widehat{g})(\xi,t)\|_2$. 

We start with estimating $\|H_{6}(\widehat{f},\widehat{g})(\xi,t)\|_2$. Consider $|\xi-\eta| \sim 2^{j}$ and $|\eta| \les 2^j$, then by \eqref{Zlowb}, we have 
\[ | Q_6|| A|^{-1} \les  2^{6j} \ t (1+ t |\xi-\eta|^4)^{-1}\les 2^{6j} \min{(2^{-4j},t)}. 
\]
Moreover, if we define the multiplier 
\begin{align}\label{m6t}
    m(\xi,\eta,t)= \frac{Q_6 \psi_j(\xi-\eta)\chi(2^{-j}|\eta|)}{A \: 2^{6j} \: \min(2^{-4j},t)}
\end{align}
 we see that $m$ defines a Coifman-Meyer multiplier satisfying \eqref{CM_con}. More specifically, using frequency localization, we can normalize the multiplier to generate $T_{m}$ whose operator norm is independent on $t$. Hence, by symmetry and Littlewood-Paley projection, $P_j$, we estimate
\begin{equation} \label{LP}
	|H_{6}(\widehat{f},\widehat{g})(\xi,t)| \les 
 \sum_{j\in\mathbb{Z}} 2^{6j}\min{(2^{-4j},t)}|T_m(P_je^{-it\Delta^2}h,P_{\leq j}e^{-it\Delta^2}h )|.  
\end{equation}
Thus, by combining \eqref{LP} with Lemma~\ref{CM_neg_deg}, we obtain the following estimate;
\begin{equation*}
	\begin{split}
		t^2 \|H_{6}&(\widehat{h},\widehat{h})\|_{L^2_x}\lesssim t^2 \sum2^{2j}\|T_m(P_je^{-it\Delta^2}h,P_{\leq j}e^{-it\Delta^2}h )\|_{L^2}\\
		&\lesssim t^2\sum \min \left(2^{2j} \|P_je^{-it\Delta^2}h\|_{L^{10}}  \|P_{\leq j}e^{-it\Delta^2}h\|_{L^{\frac{5}{2}}},\ 2^{2j} \| P_je^{-it\Delta^2}h\|_{L^\infty} \|  P_{\leq j}e^{-it\Delta^2}h )\|_{L^2}\right )\\
		&\lesssim t^2 \sum \min \left( t^{-2} 2^{-2j}\|h\|_{L^\frac{10}{9}} 2^{\frac{j}{2}} \|P_{\leq j} h\|_{L^2}, \ t^{-\frac{5}{2}}2^{-3j}\|h\|_{L^1} 2^{2j}\|h\|_{L^\frac{10}{9}} \right)\\
		&\lesssim \sum \min\left(2^{\frac{j}{2}} \|x^{2+}h\|_{L^2}^2, \ t^{-\frac{1}{2}}2^{-j}\|x^{2+}h\|_{L^2}\|x^{2.5+}\|_{L^2}  \right)\\
		&\lesssim \sum\min\left(2^{\frac{j}{2}}t^{0+}, \ 2^{-j}t^{-\frac{1}{2}+\frac{1}{32}+}  \right)\|f\|_{X}^2 
		\lesssim \|f\|_{X}^2
	\end{split}
\end{equation*}

In the third inequality, we used the estimate $\eqref{ineq:ben}$ along with the following consequence of Bernstein's inequality in $\R^5$:
\begin{align}\label{bern}
\|P_{\leq j} f\|_{L^q} \les 2^{5j(\frac{1}{p}-\frac{1}{q})}\|f\|_{L^p}, \,\,\,\ 1\leq p \leq q \leq \infty. 
\end{align}
Also, in the last inequality we used  the weighted estimate \eqref{weighted}
in conjunction with the growth condition in the $X$-norm. 
\begin{remark}
 In our estimation, we rely on \eqref{ineq:ben} to mitigate time  and frequency growth. However, using this estimate alone is not ideal, as it may adversely affects the summation over \( j \) for small values of \( j \). 
 To address this, we balance it with \eqref{bern}, which allows us to gain $2^j$ for small $j$.We always conclude our estimates by relating the \( L^p \) norms to weighted \( L^2 \) norms with \eqref{weighted}, thereby leveraging the growth conditions in the \( X \)-norm. Nonetheless, it is important to note that we can not use \eqref{ineq:ben} to gain frequency decay when estimating the term localized with \( P_{\leq j} \). 
 While this technique is not the only one we will employ, it is a crucial component to follow the computations in the remainder of the paper.
\end{remark}

%\begin{align*}
%2^{-2j}\| P_{\leq j} h\|_{L^2} \les \|h\|_{L^{\frac{10}{9}}} \les \| \la x\ra^{2+}h\|_{L^2} \les t^{0+} \|f\|_X, \,\,\ \|h\|_{L^{1}} \les \| \la %x\ra^{2,5}h\|_{L^2} \les t^{\frac{1}{32}} \|f\|_X. 
%\end{align*}

 We next consider $H_{6}(\widehat{f}\:,\widehat{g})(\xi,t)$. Using \eqref{eqn_g}, the definition of $\widehat{g}$ we observe that $H_{6}(\widehat{f}\:,\widehat{g})(\xi,t)$ can be written as a trilinear term, i.e.,
\begin{align*}
	H_{6}(\widehat{f}\:,\widehat{g})(\xi,t)=t^2\int_{\R^{5}} \int_{\R^{5}} \frac{Q_6(\xi,\eta)}{A}\frac{1}{C}e^{it\psi}\widehat{f}(\xi-\eta)\widehat{f}(\eta-\sigma)\widehat{f}(\sigma)d\sigma d\eta 
	\end{align*}
where $C:=C_t=\frac{1}{t}+iZ(\eta,\sigma)$ and $\psi=|\xi|^4-|\xi-\eta|^4-|\eta-\sigma|^4-|\sigma|^4$.
We consider cut-off functions $\Psi_j$ as in \eqref{cut_2} and further decompose $F(\xi,t)$ into a sum of the following terms:
\begin{equation}\label{5.12_1_1}
	H_{6,i,\ell}(\widehat{f},\widehat{f},\widehat{f})(\xi,t):= \iint \frac{Q_6(\xi,\eta)}{A \: C} \Psi_{i}(\xi,\eta)\Psi_{\ell}(\eta,\sigma)  e^{it\psi}\widehat{f}(\xi-\eta)\widehat{f}(\eta-\sigma)\widehat{f}(\sigma)d\sigma d\eta
\end{equation}
For each of these terms, we use Littlewood-Paley decomposition. In particular, we represent $H_{6,1,1}(\widehat{f},\widehat{f},\widehat{f})(\xi,t)$ as 
\begin{align*}
	 \sum_j\sum_k\sum_{l\leq k+10} \iint \frac{Q_6(\xi,\eta)}{A}\frac{1}{C}\Psi_1(\xi,\eta)\Psi_1(\eta,\sigma) \psi_j^2(\xi-\eta) \psi_k^2(\eta-\sigma)\varphi_l(\sigma) e^{it\psi}\widehat{f}(\xi-\eta)\widehat{f}(\eta-\sigma)\widehat{f}(\sigma)   d\sigma d\eta
\end{align*}
As we did in \( H_{6}(\widehat{h}, \widehat{h})(\xi,t) \), our goal is to normalize the multiplier to a Coifman–Meyer multiplier, now with a flag singularity. Letting

\begin{align}\label{m611t}
m_2(\xi,\eta) = \frac{Q_6(\xi,\eta)\Psi_1(\xi,\eta)\psi_j(\xi-\eta)} {A \: 2^{6j} \: \min(t,2^{-4j})},\,\,\ m_3(\eta,\sigma) = \frac{\Psi_1(\eta,\sigma)\psi_k (\eta-\sigma)}{C \ \min(t,2^{-4k})}
\end{align}
We observe that \( m_*(\xi, \eta, \sigma) = m_2(\xi, \eta) \, m_3(\eta, \sigma) \) defines a Coifman–Meyer multiplier with flag singularity. Moreover, the operator norm of the operator associated with \( m_* \) is independent of \( t \).  Consequently, we have the following bound
\begin{equation*}
	t^2 \|H_{6,1,1}(\widehat{f},\widehat{f},\widehat{f})\|_{L^2}\leq t^2\sum_j\sum_k\sum_{l\leq k+10}2^{6j}\min(t,2^{-4j})\min(t,2^{-4k})\|T_{m^*}(P_ju, P_ku, P_lu)\|_2.
\end{equation*}
Let us first consider the case when $t^{-1}\leq 2^j\leq t$. Using Theorem~\ref{th:flagcf} together with the inequalities \eqref{ineq:ben} and \eqref{bern} we further can bound $t^2 \|H_{6,1,1}(\widehat{f},\widehat{f},\widehat{f})\|_2$ by: 

\begin{align*}
	  t^\frac{5}{2}\sum_j\sum_k&\sum_{l\leq k+10} \min\Big(2^{4j}2^{-4k}\|P_ju\|_{10}\|P_ku\|_\frac{10}{3}\|P_{\leq l}u\|_{10}, 2^{4j}2^{-4k}\|P_ju\|_{10}\|P_ku\|_{10}\|P_{\leq l}u\|_\frac{10}{3}\Big) \\
		&\lesssim t^\frac{1}{2}\sum_j\sum_k\min\Bigg(\sum_{l\leq k+10} 2^{4l}2^{-k}\|f\|_\frac{10}{9}^3,\sum_{l\leq k+10}t^{-2}2^{-8k}2^{3l} \|f\|_\frac{10}{9}^3\Bigg)\\
		&\lesssim t^\frac{1}{2}\log(t) \|f\|_X^3\sum_k\min\left(2^{3k}t^{0+} , \ t^{-2+}2^{-5k} \right) \lesssim  \|f\|_X^3. 
\end{align*}
In the last inequality, we utilize the condition $2^j \leq t$ to sum in $j$. We give Remark~\ref{middle_j} before we consider the other cases: $2^j\leq t^{-1}$ and $2^j \geq t$. 

\begin{remark}\label{middle_j}
	For estimating series $\sum_j |a_j|$, the ’determining‘ part is the middle-frequency case, \textit{i.e.,} when $t^{-1} <2^j<t $. More specifically, if $ 2^{j} <t^{-1}$ then for any $p,q <\infty$ with $\f1{p} + \f1{q} + \f1{s} =\f12$ we have 
    \begin{align}
        t^{k} \sum_{2^j \leq t^{-1}} 2^{nj} \|T_m(P_j(e^{it\Delta^2}f_1), &e^{it\Delta^2}f_2,e^{it\Delta^2}f_3) \|_2 \les t^{(k-n)}  \sum_{2^j \leq t^{-1}} \| P_j(e^{it\Delta^2}f_1)\|_{s} \| e^{it\Delta^2}f_2\|_p \| e^{it\Delta^2}f_3\|_q \nn   \\ & \les  t^{(k-n-\f54 + \frac{5}{2s})} \sum_{2^j \leq t^{-1}} 2^{ (\f52 - \frac{5}{s}) j} \|f_1\|_2 \|f_2\|_{\frac{q}{q-1}}\|f_3\|_{\frac{p}{p-1}}  \\ 
        &\les  t^{(k-n-\f{15}4 + \frac{5}{s})} \|f_1\|_2 \|f_2\|_{\frac{q}{q-1}}\|f_3\|_{\frac{p}{p-1}}. \label{lowerb} 
    \end{align}   
Moreover, for any $n < 5 -\frac{10}{r}$ and $\frac{1}{p}+\frac{1}{q} + \frac{1}{r} =\f12$, we have 
 \begin{align}
        t^{k}\sum_{2^{-j} \leq t^{-1}} 2^{nj} \|T_m(P_j(e^{it\Delta^2}f_1), &e^{it\Delta^2}f_2,e^{it\Delta^2}f_3) \|_2 \nn \\ &\les t^{k-\f52(1-\f1{p}-\f1{q})} \sum_{2^{-j} \leq t^{-1}}2^{nj} \| P_j(e^{it\Delta^2}f_1)\|_r \|f_2\|_{\frac{q}{q-1}}\|f_3\|_{\frac{p}{p-1}} \nn \\
        & \les t^{k-\f{15}{4} + \frac{5}{2r}}  \sum_{2^{-j} \leq t^{-1}} 2^{j(n-5+\frac{10}{r})} \|\la x \ra^{\frac{5(r-2)}{2r}+} f_1\|_2 \|f_2\|_{\frac{q}{q-1}}\|h\|_{\frac{p}{p-1}} \nn \\
        & \les t^{k-\f{15}{4} + \frac{5}{2r}-\alpha+}\|\la x \ra^{\frac{5(r-2)}{2r}+} f_1\|_2 \|f_2\|_{\frac{q}{q-1}}\|f_3\|_{\frac{p}{p-1}}  \label{upperb} 
    \end{align}
where $\alpha = 5- n+\frac{10}{r}>0 $. Furthermore, if  $T_m$ is a bilinear operator, the estimates take the form 
\begin{align}
 & t^{k}\sum_{2^j \leq t^{-1}} 2^{nj} \|T_m(P_j(e^{it\Delta^2}f_1), e^{it\Delta^2}f_2) \|_2 \les t^{k-n-\f52} \|f_1\|_2 \|f_2\|_2 \label{blowerb}\\ 
  &t^{k}\sum_{2^{-j} \leq t^{-1}} 2^{nj} \|T_m(P_j(e^{it\Delta^2}f_1), e^{it\Delta^2}f_2) \|_2 \les t^{k-\f{5}{2} + \frac{5}{2r}-\alpha+}\|\la x \ra^{\frac{5(r-2)}{2r}+} f_1\|_2 \|f_2\|_{\frac{2r}{r+2}}\label{bupperb}
\end{align}

In all these estimates if $e^{it\Delta^2}f_i =u$ then the related norms on the right side of the inequalities can be exchanged by $\|f_i\|_X$. 
%if we can prove a uniform bound on $a_j$: $|a_j|\leq C$, then 
	%\begin{equation*}
	%	\sum_{t^\alpha<2^j<t^\beta}|a_j|\leq C\sum_{t^\alpha<2^j<t^\beta} 1\approx C\int_{\alpha \log t}^{\beta \log t}dx=C(\beta-\alpha)\log t
	%\end{equation*}
	%Thus, we can manipulate $\alpha$ and $\beta$ however we want to make the low-frequency and high-frequency sum good without influencing the middle-frequency sum significantly. 
\end{remark}

Having Remark~\ref{middle_j}, we continue with the cases:  $2^j\leq t^{-1}$ and $2^j \geq t$. First, note that we have 
\begin{align*}
    t^2\sum_j \|T_{m^*}(P_ju, P_ku, P_lu)\|_2 \les t^3\sum_j\sum_k\sum_{l\leq k+10}2^{2j}\|T_{m^*}(P_ju, P_ku, P_lu)\|_2. 
\end{align*}
Moreover, by Young's inequality we have $\|P_j(f_1)\|_p \les \|f\|_p$. Therefore, using \eqref{lowerb} and \eqref{upperb} with $r=\infty$ we obtain $t^2\|H_{6,1,1}(\widehat{f},\widehat{f},\widehat{f})\|_2\leq \|f\|_X^3$. 

Note that $H_{6,1,2}(\widehat{f}\:,\widehat{g})$ can be addressed in a similar manner by expressing  $H_{6,1,2}(\widehat{f},\widehat{f},\widehat{f})(\xi,t)$ in the following way
\begin{equation*}
	\sum_j\sum_k\sum_{l\leq k+10} \iint \frac{Q_6(\xi,\eta)}{A \: C}\Psi_1(\xi,\eta)\Psi_2(\eta,\sigma)   \psi_j^2(\xi-\eta)\psi_k^2(\sigma)\varphi_l(\eta-\sigma) e^{it\psi}\widehat{f}(\xi-\eta)\widehat{f}(\eta-\sigma)\widehat{f}(\sigma)   d\sigma d\eta. 
\end{equation*}
and using the Coifmann-Meyer multiplier  $m^*=m^*_2m^*_3$
\begin{align} \label{m612t}
m^*_2(\xi,\eta) = \frac{Q_6(\xi,\eta)\Psi_1(\xi,\eta)\psi_j(\xi-\eta)} {A \: 2^{6j} \: \min(t,2^{-4j})} \,\,\ m^*_3(\eta,\sigma) = \frac{\Psi_2(\eta,\sigma)\psi_k (\sigma)}{C \ \min(t,2^{-4k})}. 
\end{align}

We next focus on $H_{6,2,1}(\widehat{f}\:,\widehat{g})$, where we localize in $\eta-\sigma$ and $\eta$ to express $H_{6,2,1}(\widehat{f},\widehat{f},\widehat{f})$ as 
\begin{align}\label{m621}
	\sum_k\sum_{j\leq k+10}\iint \frac{Q_6(\xi,\eta)}{A \:C}\Psi_2(\xi,\eta)\Psi_1(\eta,\sigma) \varphi_j(\eta)\psi_k^2(\eta-\sigma)e^{it\psi}\widehat{f}(\xi-\eta)\widehat{f}(\eta-\sigma)\widehat{f}(\sigma)   d\sigma d\eta
\end{align}
Similar to $H_{6,1,1}(\widehat{f},\widehat{f},\widehat{f})$, we let 
\begin{align}\label{m621t}
m_2(\xi,\eta) = \frac{Q_6(\xi,\eta) \Psi_2(\xi,\eta)\varphi_j(\eta)}{ A \ 2^{6j}\min(t,2^{-4j})}, \,\,\, m_3(\eta,\sigma)= \frac{\Psi_1(\eta,\sigma) \psi_k(\eta-\sigma)}{C \ \min(t,2^{-4k})}
\end{align}
and notice that $m^*= m_2m_3$ defines  a Coifman-Meyer multiplier with flag singularity, that does not depend on time. 
This enables the following estimate:
\begin{equation*}
	\begin{split}
		t^2\|H_{6,1,2}(\widehat{f},\widehat{f},\widehat{f})\|_2&\les t^2\left\|\sum_k\sum_{j\leq k+10} 2^{6j}\min(t,2^{-4j})\min(t,2^{-4k})e^{-it\Delta^2}T_{m^*}(u,P_ku,u) \right\|\\
		&\lesssim t^2\sum_k \sum_{j\leq k+10}2^{2j}2^{-4k}\|u\|_{10}^2\|P_ku\|_\frac{10}{3} \\
		&\lesssim \sum_k2^{-2k}\min(2^{3k}, t^{-1}2^{-k} )\|f\|_\frac{10}{9}\|f\|_X^2 \lesssim \|f\|_X^3
	\end{split}
\end{equation*}
We note that $H_{6,2,2}(\widehat{f},\widehat{f},\widehat{f})$ can be  handled similarly by localizing in $\eta$ and $\sigma$, i.e :
\begin{equation*}
	\sum_k\sum_{j\leq k+10}\iint \frac{Q_6(\xi,\eta)}{A \: C}\Psi_2(\xi,\eta)\Psi_2-(\eta,\sigma) \varphi_j(\eta)\psi_k^2(\sigma)e^{it\psi}\widehat{f}(\xi-\eta)\widehat{f}(\eta-\sigma)\widehat{f}(\sigma)   d\sigma d\eta
\end{equation*}
and using the multiplier $m^*=m^*_2 m^*_3$ where 
\begin{align} \label{m622t}
m_2(\xi,\eta) = \frac{Q_6(\xi,\eta) \Psi_2(\xi,\eta)\varphi_j(\eta)}{ A \ 2^{6j}\min(t,2^{-4j})}, \,\,\, m_3(\eta,\sigma)= \frac{\Psi_2(\eta,\sigma) \psi_k(\sigma)}{C \ \min(t,2^{-4k})}
\end{align}

Thus, we have completed the estimation of \( H_{6}(\widehat{f}, \widehat{g}) \). Note that \( H_{6}(\widehat{g}, \widehat{h}) \) can be estimated in a similar manner, since the time growth of the weighted \( L^2 \) norms of \( h \) is slower than that of \( f \), as suggested by \eqref{Boot_est_h}. Therefore, we omit the details for this term. This completes the proof of the statement. 
\end{proof}
%\subsection{Estimates of $\|x^3g\|_{L^2}$}
\begin{lemma}
	Assuming \eqref{Boot_est_g} and \eqref{Boot_est_h}, we have
	\begin{equation*}
		\|x^3g\|_{L^2_x}\les t^{\frac{1}{2}+\frac{1}{47}}\|f\|_X^2
	\end{equation*}
\end{lemma}
\begin{proof}
	Applying $\partial_\xi^3$ to $\widehat{g}$, we can express
	\begin{align}
		\partial_\xi^3\widehat{g}(\xi,t)&\label{x3g_1} =t^3\int \frac{Q_9}{A}e^{it\varphi}\widehat{f}(\xi-\eta,s)\widehat{f}(\eta,s)d\eta\\
		\label{x3g_2} &+t^2\int \frac{Q_6}{A} e^{it\varphi}\widehat{xf}(\xi-\eta,s)\widehat{f}(\eta,s)d\eta\\
		\label{x3g_3} &+t^2\int \frac{q_5}{A}e^{it\varphi}\widehat{f}(\xi-\eta,s)\widehat{f}(\eta,s)d\eta  \\
            \label{x3g_5}&+t\int \frac{q_3}{A}e^{it\varphi}\widehat{x^2f}(\xi-\eta,s)\widehat{f}(\eta,s)d\eta +\int\frac{1}{A}e^{it\varphi}\widehat{x^3f}(\xi-\eta,s)\widehat{f}(\eta,s)d\eta\\
		\label{x3g_6}&+ t\int  \frac{q_2}{A}e^{it\varphi}\widehat{xf}(\xi-\eta,s)\widehat{f}(\eta,s)d\eta+ t\int \frac{q_1}{A}e^{it\varphi}\widehat{f}(\xi-\eta,s)\widehat{f}(\eta,s)d\eta.
	\end{align}
 The terms in \eqref{x3g_5} and the first term in \eqref{x3g_6}  can be bounded by $t^{\f14} \|f\|_{X}^2$ using \eqref{fracwuse}. For the second term in \eqref{x3g_6}, we first observe that
\begin{align}\label{m1b}
\|T_{\frac{q_1}{A}}(u,u)\| &\les \| \nabla_{t}^{-3} f\|_{L^2} \|u\|_{L^{\infty}} \les t^{\frac{3}{4} - \frac{5}{4}( \frac{1}{2} - \alpha)} \|f\|_{L^{\frac{1}{1-\alpha}}} \|u\|_{L^{\infty}} \nn \\ 
& \les t^{-\frac{9}{8}+ \frac{5 \alpha}{4}} \| \la x \ra^{\frac{5-5\alpha}{2} +} f\|_{L^2} \|f\|_{X}, \,\,\ 0 \leq \alpha \leq \f12. 
\end{align}
Choosing \(\alpha = \frac{1}{10}\), we estimate the $L_2$ norm of this term by $t^{0+}\|f\|_{X}^2$. To estimate \(\|\eqref{x3g_3}\|_{L^2}\), we normalize the multiplier similar to $H_6(\widehat{h}, \widehat{h})$, see \eqref{LP}, to obtain 
\begin{align*}
\|\eqref{x3g_3}\|_{L^2} \les t^2\sum_j 2^j \|P_ju\|_{10}\|P_{\leq j}u\|_\frac{5}{2} 
 \les  t^2\sum_j \min(t^{-2}2^{-\frac{j}{2}} , 2^{\frac{15}{2}j})\|f\|_{\frac{10}{9}} \les t^\frac{1}{8}\|f\|_X^2.
 \end{align*}
Here, in  the third inequality we used the estimate 
 $\|P_j u\|_{10} \les \min( 2^{-4j} t^{-2}, 2^{4j}) \|f\|_{\frac{10}{9}}$.

Next, we focus on \eqref{x3g_2}. Observe that $\eqref{x3g_2}=H_6(\widehat{f},\widehat{f})(\xi,t)$ and 
therefore similar to \eqref{LP}, we have 
$$
 \|H_6(\widehat{xf},\widehat{f})\|_{L^2} \les \sum_{j} 2^{2j} \Big( \|T_m(P_je^{-it\Delta^2}xf, P_{\leq j}u)\|_2	+ \|T_m(P_{\leq j}e^{-it\Delta^2}xf, P_{j}u)\|_2	\Big). 
$$
We consider  when $t^{-1} \leq 2^j\leq t$ and obtain 
	\begin{align}
   & t^2\sum_{t^{-1}\leq 2^j \leq t} 2^{2j} \Big( \|T_m(P_je^{-it\Delta^2}xf, P_{\leq j}u)\|_2	+ \|T_m(P_{\leq j}e^{-it\Delta^2}xf, P_{j}u)\|_2	\Big) \nn \\
      & \les t^2 \sum_{t^{-1}\leq 2^j \leq t}  2^{2j} \Big(\|P_je^{-it\Delta^2}xf\|_\frac{10}{3}\|u\|_5 + \|P_ju\|_{10}\|P_{\leq j}e^{-it\Delta^2}xf\|_\frac{5}{2}\Big) \nn \\ 
      & \les t^\frac{1}{4} \|xf\|_\frac{10}{7}\|f\|_{\f54} + \log(t)\|f\|_\frac{10}{9}\|xf\|_\frac{5}{4} \les t^{\frac{1}{4}+\frac{1}{80}+}\|f\|_X^2 \label{x3g_2est}
     \end{align}
     where we use \eqref{weighted} in the last inequality. 
Having Remark~\ref{middle_j}, we have by \eqref{blowerb} 
\begin{align}
 t^2\sum_{2^{j} \leq t^{-1}} 2^{2j} \Big( \|T_m(P_je^{-it\Delta^2}xf, P_{\leq j}u)\|_2	+ \|T_m(P_{\leq j}e^{-it\Delta^2}xf, P_{j}u)\|_2	\Big) \les  \|f\|_X^2.  \label{x3g_2small}
\end{align}
Moreover, letting $r=\frac{10}{3}-$ and $r=\infty$ in \eqref{bupperb}, we have 
\begin{multline}
 t^2\sum_{2^j \leq t^{-1}} 2^{2j} \Big( \|T_m(P_je^{-it\Delta^2}xf, P_{\leq j}u)\|_2	+ \|T_m(P_{\leq j}e^{-it\Delta^2}xf, P_{j}u)\|_2	\Big)  \\ 
 \les t^{\f14} \| \la x \ra^{2+} f\|_2 \|f\|_X + t^{-\f12}\|xf\|_2\|f\|_X \label{x3g_2large}
\end{multline}
 
Combining \eqref{x3g_2est}, \eqref{x3g_2small} and \eqref{x3g_2large}, we obtain 
\begin{align*}
t^2\|\ H_6(\widehat{xf},\widehat{f})\|_{L^2} \les t^{\frac{21}{80}} \|f\|_{X}^2 
\end{align*}

%When $t\leq 2^j$, we use  $\|xf\|_{\frac{5}{4}} \les \| \la x \ra^{\f32+} xf\|_2$ to estimate 
%\begin{align}
%\|\eqref{x3g_2_1}\|_2\les t^2\sum_j2^{2j} \|P_je^{-it\Delta^2}xf\|_5\|u\|_{\frac{10}{3}}  \les t^{-\f12}\sum_j2^{-2j} \|xf\|_{\frac{5}{4}}\|f\|_{\frac{10}{9}} %\les \|f\|_X^2
%\end{align}

To estimate $\|\eqref{x3g_1}\|_2$, we employ a method similar to the one used for estimating \eqref{x2g_1}. Specifically, we define the operator 
$$
H_9(\widehat{f},\widehat{g})(\xi,t):=\int \frac{Q_9}{A}e^{it\varphi}\widehat{f}(\xi-\eta,t)\widehat{f}(\eta,t)d\eta
$$
and use the decomposition \eqref{Decom_Bi} to reduce the problem of estimating $\|\eqref{x3g_1}\|_2$ to study 
\begin{align}  \label{H9}
H_9(\widehat{h},\widehat{h})(\xi,t), \,\,\,\ H_9(\widehat{f},\widehat{g})(\xi,t), \,\,\,\ H_9(\widehat{g},\widehat{h})(\xi,t) 
\end{align}

We start with the term $H_9(\widehat{h},\widehat{h})(\xi,t)$. First not that  when $ 2^j\leq t^2$:
\begin{align*}
\|H_9(\widehat{h},\widehat{h})(\xi,t)|_2 & \les  t^3\sum_{2^j \leq t^2} 2^{5j} \|T_m(P_je^{-it\Delta^2}h, P_{\leq j}e^{-it\Delta^2}h)\|_2 \les t^3 \sum_{2^j \leq t^2} 2^{5j} \|P_je^{-it\Delta^2}h\|_\infty\|e^{-it\Delta^2}h\|_2 \\
& \les t^3 \sum_{2^j \leq t^2} 2^{5j} \min( t^{-\f54} \|f\|_X\|h\|_2 , 2^{-5j}  t^{-\f52}\| h\|_1 \|h\|_{2})\les t^{\frac{1}{2}+\frac{1}{48}+}\|f\|_X^2.
\end{align*}
Here $m$ is the same multiplier as in \eqref{m6t} where $Q_6$ is exchanged with $Q_9$. For high frequencies that is when $2^j \geq t^2$, we need to study $\|h\|_{H^6}$. In particular, \eqref{upperb} in Remark~\ref{middle_j} is not available when $n\geq 5$.  Multiplying $h$ with $Q_6$, we obtain 
	\begin{multline*}
			\widehat{\partial_x^6h}(\xi,t)=\int_1^t s^{-1}\int\frac{q_6(\xi,\eta,s) }{A}e^{is\varphi}\widehat{f}(\xi-\eta,s)\widehat{f}(\eta,s)d\eta ds \\ + \int_1^t \int \frac{Q_6}{A} e^{is\varphi}\widehat{f}(\xi-\eta,s)\partial_s\widehat{f}(\eta,s)d\eta ds
			 +\int_1^t s^{-1}\int \frac{Q_7}{A}e^{is\varphi}\widehat{xf}(\xi-\eta,s)\widehat{f}(\eta,s)d\eta ds \nn
	\end{multline*}
One can easily see that $L^2_\xi$ norm of the first term can be bounded by \(\|f\|_X^2\) using the bound on $\|H_6(\widehat{f}, \widehat{f})\|_2$. For the second term a in $\widehat{\partial_x^6h}(\xi,t)$, we have 
\begin{align*}
	\int_1^t \|H_6(\widehat{f}, \widehat{\partial_sf}) (\xi,s)\|_2 \: ds \les  \int_1^t \sum_j2^{2j} \Big( \|T_m(P_ju,P_{\leq j}u^2)\|_2 + \|T_m(P_{\leq j}u,P_{j}u^2)\|_2 \Big) ds
    \end{align*} 
In a similar manner to \eqref{x3g_2small}, \eqref{x3g_2est}, and \eqref{x3g_2large} we see that 
$\|H_6(\widehat{f}, \widehat{\partial_sf}) (\xi,s)\|_2 \les s^{-\frac{63}{32}}$. 
Moreover, for the third term in $\widehat{\partial_x^6h}(\xi,t)$, we have 
    
\begin{multline*}
   \|H_7(\widehat{xf}, \widehat{f}) (\xi,s)\|_2 \: \les 
\sum_j 2^{3j} \Big( \|T_m(P_je^{-is\Delta^2}(xf), P_{\leq j}u)\|_2 + \|T_m(P_{\leq j}e^{-is\Delta^2}(xf), P_{j}u)\|_2 \Big)  
\end{multline*}
Using $r=5+$ and $r=\infty$ in \eqref{blowerb} in a similar way to \eqref{x3g_2large}  together with $\eqref{bupperb}$,  it is enough to only consider

\begin{align*}
& \int_1^t s^{-1}\sum_{s^{-1} \leq 2^j \leq s} 2^{3j} ( \|P_j(e^{-is\Delta^2}(xf))\|_5\|P_{\leq j}u\|_\frac{10}{3}ds + \|P_{\leq j}(e^{-is\Delta^2}(xf))\|_\frac{5}{2}\|P_{ j}u\|_{10} \Big)ds \\
&\les \int_1^t \sum_{s^{-1} \leq 2^j \leq s}\big( s^{-\f52} 2^j \|xf\|_{\f54} \| u\|_{2} + s^{-3} \| xf\|_{2} \|f\|_{\frac{10}{9}} \big) \les \|f\|_X^2.   \end{align*}
which completes the estimate \(\|h\|_{H^6}\les \|f\|_X^2\). 

We now return to estimating $\|H_9(\widehat{h},\widehat{h})\|_2$ for $2^j\geq t^2$. We compute
	\begin{equation*}
		\begin{split}
			t^3\sum_{2^j \geq t^2} 2^{5j}\|T_m(P_je^{-it\Delta^2}h, P_{\leq j}e^{-it\Delta^2}h)\|_2&\les t^3\sum_j 2^{5j}\|P_je^{-it\Delta^2}h\|_2\|e^{-it\Delta^2}h\|_\infty\\
			&\les t^{\frac{7}{4}}\sum_j2^{-j}\|h\|_{H^6}\|f\|_X\les\|f\|_X^3
		\end{split}
	\end{equation*}
which gives \(\|H_9(\widehat{h},\widehat{h})\|_2\les t^{\frac{1}{2}+\frac{1}{48}+}\|f\|_X^2 \).

We continue with estimating $\|H_9(\widehat{f}, \widehat{g})\|_2$. We use \eqref{eqn_g} to represent $H_9(\widehat{f}, \widehat{g})$  as a trilinear term: 
	\begin{equation*}
		H_9(\widehat{f}, \widehat{f},\widehat{f})(\xi,t)=\iint \frac{Q_9}{A}\frac{1}{C}e^{is\psi}\widehat{f}(\xi-\eta,t)\widehat{f}(\eta-\sigma,t)\widehat{f}(\sigma,t)d\sigma d\eta
	\end{equation*}
Using cut-off functions we focus on $H_{9,i,\ell}(\widehat{f}, \widehat{g})(\xi,t)$ for $i,\ell=1,2$. Employing the multipliers given in \eqref{m611t} and \eqref{m612t}, with \( Q_6 \) replaced by \( Q_9 \) for the cases \( \ell = 1 \) and \( \ell = 2 \), respectively, we obtain:
\begin{align*}
  \|H_{9,1,\ell}(\widehat{f},\widehat{f},\widehat{f})\|_2 \les \sum_j\sum_k 2^{5j}\min(t,2^{-4k})\|T_{m^*}(P_ju,P_ku,u)\|_2
\end{align*}
We  consider the case where $t^{-1}\leq 2^k \leq t$ first to see 
	\begin{equation*}
		\begin{split}
			t^3  \|H_{9,1,\ell}(\widehat{f},\widehat{g})\|_2&\les t^4\log(t)\sum_j2^{5j}\min(\|u\|_6^3 ,\|P_ju\|_\infty\|u\|_4^2, \|P_ju\|_\frac{5}{2}\|u\|_{20}^2)\\
			&\les t^4\log(t)\sum_j2^{5j}\min(t^{-\frac{5}{2}}\|f\|_X^3, t^{-\frac{15}{4}}2^{-5j}\|f\|_1\|f\|_X^2, t^{-\frac{9}{4}}2^{-\frac{9}{2}j}\|f\|_{H^{10}} \|f\|_X^2 )\\
			&\les t^{\frac{1}{2}+\frac{1}{80}+\varepsilon}\|f\|_X^3
		\end{split}
	\end{equation*}
For the other cases note that $ \sum_{j} \| 2^{5j} P_j(f)\|_2 = \|f\|_{H^5} \les \|f\|_X$. Hence,  we can use Remark~\ref{middle_j} in summing $k$ with $r=2$ to complete the estimate on $t^3  \|H_{9,1,\ell}(\widehat{f},\widehat{g})\|_2$. 

Furthermore, using the multipliers given in \eqref{m621t} and \eqref{m622t}, with \( Q_6 \) replaced by \( Q_9 \) for the cases \( \ell = 1 \) and \( \ell = 2 \), respectively, we have 
	\begin{align*}
	\|H_{9,2,\ell}(\widehat{f},\widehat{f},\widehat{f})\|_2\les \sum_k\sum_{j\leq k+10}\sum_{l\leq j+10}2^{5j}	\min(t,2^{-4k})\|T_{m^*}(P_lu,P_ku,u)\|_2
	\end{align*}
Noting that $2^{5j} \min(t,2^{-4k}) \les 2^{k}$ in the domain of the sum and by Remark~\ref{middle_j} we only consider when $t^{-1}\leq 2^k \leq t$ to estimate 
	\begin{multline*}
			 t^3\sum_{t^{-1}\leq 2^k \leq t}\sum_{j\leq k+10}\sum_{l\leq j+10}2^{-4k}2^{5j}\|P_lu\|_{\frac{10}{3}}\|P_ku\|_{10}\|u\|_{10}\\
			\les \sum_{t^{-1}\leq 2^k \leq t} 2^{-8k}\sum_{j\leq k+10}\sum_{l\leq j+10}2^{5j}2^{3l}\|f\|_\frac{10}{9}^2\|f\|_X \les t^\varepsilon\|f\|_X^3
\end{multline*}
	This completes the estimate $t^3\|H_{9}(\widehat{f},\widehat{g})\|_2 \les t^{\f12 + \frac{1}{47}} \|f\|_X^2$. %For \eqref{x3g_1_3}, by change of variable, we have:
	%\begin{equation*}
	%	\eqref{x3g_1_3}=t^3\int\frac{P_9}{A^\prime}e^{it\varphi}\widehat{h}(\xi-\eta,t)\widehat{g}(\eta,t)d\eta
	%\end{equation*}
	As the weighted bounds that $f$ satisfies are also satisfied by $h$, $t^3\|H_{9}(\widehat{g},\widehat{h})\|_2 $ can be estimated in the same way after a change of variable. Therefore, we complete the proof.  
	\end{proof}

\section{Estimates on $h$}\label{sec:hbounds}
Recall by \eqref{eqn_h} \(h\) is the sum of \(h_1,h_2,h_3,h_4\), where
\begin{align*}
&h_1(x,t)=-i\int_1^t \frac{1}{s^2} e^{is\Delta^2}\big[T_{\frac{1}{A^2}}(u,u)\big](x,s) ds,\,\,\, h_2(x,t)= 2i \int_1^t  e^{is\Delta^2}\big[T_{\frac{1}{A}}(u,u^2)\big](x,s) ds \\ 
& h_3(x,t)= i \int_1^t \frac{1}{s} e^{is\Delta^2}\big[T_{\frac{1}{A}}(u,u)\big](x,s) ds, \,\,\,\, h_4(x,t)= i\int_1^t  e^{is\Delta^2}\big[T_{\frac{P \varphi_{\eta}}{A}}(u,u)\big](x,s) ds
\end{align*}
where we used 
\begin{equation}\label{dt_hat_f}
	\partial_t \hat{f}(\xi,t)=ie^{it|\xi|^4}\widehat{u^2}(\xi,t)
\end{equation}
 expressing $h_2$, and used $\partial_\eta( e^{is \varphi}) = i(\partial_{\eta}\varphi) e^{is \varphi} s =:\varphi_{\eta} e^{is \varphi} s  $ expressing $h_4$. In this section, we prove \eqref{est_h} for all these terms under the bootstrap assumption \eqref{Boot_est_g} and \eqref{Boot_est_h}. 

Unlike $g$, it is not straightforward to obtain the natural time decay $t^{-\f54}$ on $\|e^{it \Delta^2}f\|_{\infty}$  using the Coifman-Meyer theory. Instead, we will prove it using the dispersive decaying estimate in Lemma \ref{lem:linear}, combined with the bootstrap result.

%\textcolor{red}{-I think with these representations we can shorten up until $x^2 h$ bound using the relations with g. $h_3$ definitely can be shorten, but maybe even $h_1$, the following might be true but needs to be checked 
%\[\|x^{\alpha} h_2\|_{L^2}\les \int_1^t\|u\|_{L^{\infty}} \|x^{\alpha} T_{\frac{1}{A}}(u,u)  \|_{L^2} ds\] 
%}

\begin{lemma}
Assuming \eqref{Boot_est_g} and \eqref{Boot_est_h}, we have
	\begin{equation}\label{Decay_h}
		\|h\|_{L^2}\lesssim \|f\|_X^2+\|f\|_{X}^3
	\end{equation}
\end{lemma}
\begin{proof}
Using the uniform bound $|A_s^{-1}|\les s$ and  \eqref{TAbound}, we have  
    \begin{equation}
	\|h_1\|_{L^2}+ 	\|h_2\|_{L^2}+ \|h_3\|_{L^2}+  \les \int_1^{t}  (s^{-1} +\|u\|_{L^{\infty}})\|u\|_{L^2} \|u\|_{L^{\infty}} ds \les t^{-\f14}(\|f\|_X^2+ \|f\|_X^3)
	\end{equation} 
Similarly, using $\|T_{P\varphi_{\eta}}(u,u)\|_{L^2} \les \|f\|_{L^2} \|u\|_{L^{\infty}}$, we have  $\|h_4\|_{L^2} 
   \les t^{-\f14}\|f\|_{X}^2$
\end{proof}

%\subsection{Estimates of $\|xh\|_{L^2}$}
\begin{lemma}
Assuming \eqref{Boot_est_g} and \eqref{Boot_est_h}, we have
	\begin{equation}\label{Decay_xh}
		\|xh\|_{L^2}\lesssim \|f\|_X^2+\|f\|_{X}^3
	\end{equation}
\end{lemma}
\begin{proof} Note that using \eqref{Decay_xg}, we can easily see $\|xh_3\|_{L^2} \les \|f\|_X$. Moreover, applying $\partial_\xi$ to $h_1(\xi,t)$, we obtain 
\begin{align*}
\partial_{\xi}(h_1)= \int_1^t \int e^{is\varphi} \Big[ \frac{q_3(\xi,\eta,s)}{A} 	\widehat{f}(\xi-\eta,s)+ q_0(\xi,\eta,s) \widehat{xf}(\xi-\eta,s)\Big] \widehat{f}(\eta,s)d\eta ds 
\end{align*}
 Hence, using \eqref{polyb} with $k=3$, and \eqref{fracwuse} with $k=0$, we have 
\begin{align}\label{xh1}
\| xh_1\|_{L^2} \les \int_1^{t}  s^{-\f54}  \|f\|_X^2 \les \|f\|_X^2. 
\end{align}

Next, we compute 
\begin{align*}
\partial_{\xi}(h_2)= \int_1^t \int e^{is\varphi} \Big[\frac{  q_3(\xi,\eta,s) s }{A} 	\widehat{f}(\xi-\eta,s) +  \frac{1}{A} \widehat{xf}(\xi-\eta,s) \Big] e^{i|\eta|^4} \widehat{u^2}(\eta) d\eta ds
\end{align*}
Hence, using the norm estimates as in \eqref{TAbound}, we have 
\begin{align}
\|xh_2\|_{L^2} \les \int_1^{t} (s \|\nabla_t^{-1} f\|_{L^2}\|u\|^2_{L^{\infty}} +t \|xf\|_L^2 \|u\|^2_{L^\infty}) ds \les  \|f\|_X^3 
\end{align}
Finally, applying $\partial_{\xi}$ to $h_4(\xi,t)$ in \eqref{eqn_h}
\begin{align}
\partial_{\xi}(h_4) &= \int_1^{t} s^{-1} \int \partial_\eta(e^{is\varphi}) \Big[ \frac{q_0(\xi,\eta,s)}{A} \widehat{f}(\xi-\eta,s)+\frac{q_1(\xi,\eta,s)}{A} \widehat{xf}(\xi-\eta,s)\Big] \widehat{f}(\eta,s)d\eta ds   \\ 
\label{xh_7}&+\int_1^t s^{-1}\int \frac{q_1(\xi,\eta,s)}{A}\partial_\xi \partial_\eta(e^{is\varphi})\widehat{f}(\xi-\eta,s)\widehat{f}(\eta,s)d\eta ds
\end{align}
Expanding $\partial_\eta(e^{is\varphi})$, bounding the $L^2$ norm of the first term reduces to \eqref{xh1}. Moreover,  we compute 
\begin{equation*}
	\begin{split}
\eqref{xh_7}=\int_1^t \int\frac{q_3}{A}e^{is\varphi}\widehat{f}(\xi-\eta,s)\widehat{f}(\eta,s)d\eta ds+\int_1^t \int\frac{q_4}{A}\partial_\eta(e^{is\varphi})\widehat{f}(\xi-\eta,s)\widehat{f}(\eta,s)d\eta ds
	\end{split}
\end{equation*}
Whereas the first term appeared in $\partial_{\xi}(h_1)$, by integration by parts the second term reduces to again the terms in $\partial_{\xi}(h_1)$.

\end{proof}

Before, we continue to our estimates we prove  Lemma~\ref{lem_xf_4} which will be useful for further computations. 

%\subsection{Estimates of $\|x^2h\|_{L^2}$}

\begin{lemma}\label{lem_xf_4}
	Assuming \eqref{Boot_est_g} and \eqref{Boot_est_h}, we have
\begin{align}
& \label{decay_xf_4}\|e^{-it\Delta^2}(xf)\|_{4}^2 \lesssim  t^{-\frac{69}{64}+}\|f\|_X^2
\\ \label{Decay_xh_4}
& \|e^{-it\Delta^2}(xh)\|_4^2\lesssim t^{-\frac{223}{168}+}\|f\|_X^2 
\end{align}

\end{lemma}
\begin{proof}
Using Plancherel, we have:
	\begin{align}
			\|e^{-it\Delta^2}(xf)\|_{4}^2&=\left\|\int e^{it\phi}\nabla\widehat{f}(\xi-\eta)\nabla\widehat{\overline{f}}(\eta)d\eta  \right\|_2 \\
			&\lesssim \left\|\sum_{ 2^j\leq t^{-\beta}} \int e^{it\phi}\nabla \widehat{f}(\xi-\eta)\nabla\widehat{\overline{f}}(\eta)\varphi_j(\xi-\eta)\varphi_{\leq j}(\eta)d\eta \right\|_2 \nn \\
			&\ \ \ +\left\|\sum_{ 2^j> t^{-\beta}} \int e^{it\phi}\nabla \widehat{f}(\xi-\eta)\nabla\widehat{\overline{f}}(\eta)\varphi_j(\xi-\eta)\varphi_{\leq j}(\eta)d\eta \right\|_2=:I(xf)+II(xf) \nn 
	\end{align}
where $\phi=\phi(\xi,\eta)=|\eta|^4-|\xi-\eta|^4$.  We can estimate $I(xf)$  using \eqref{bern} as 
	\begin{equation*}
		\begin{split}
			I(xf)&\leq \sum_{2^j\leq t^{-\beta}} \|P_je^{-it\Delta^2}(xf)\|_\frac{10}{3}\|P_{\leq j}e^{-it\Delta^2}(xf)\|_5 \\ 
			&\lesssim \sum_{2^j\leq t^{-\beta}} t^{-1}2^{-2j}2^{\frac{5}{2}j}\|xf\|_\frac{10}{7}^2 \les t^{-1-\frac{\beta}{2}+}\|f\|_X^2. 
		\end{split}
	\end{equation*}
	For $II$, we apply integration by parts:
	\begin{align}
		II(xf)&\label{6.15}\leq \left\|\sum_{2^j\geq t^{-\beta}}\int e^{it\phi}\widehat{f}(\xi-\eta)\Delta \widehat{\overline{f}}(\eta)\varphi_j(\xi-\eta)\varphi_{\leq j}(\eta)d\eta  \right\|_2\\
			&\label{6.16_n}+	\left\|\sum_{2^j\geq t^{-\beta}} it \int \phi_\eta e^{it\phi}\widehat{f}(\xi-\eta)\nabla \widehat{\overline{f}}(\eta)\varphi_j(\xi-\eta)\varphi_{\leq j}(\eta)d\eta  \right\|_2\\
			&\label{6.17_n}+\left\|\sum_{2^j\geq t^{-\beta}} \int e^{it\phi}\widehat{f}(\xi-\eta)\nabla \widehat{\overline{f}}(\eta)\partial_\eta(\varphi_j(\xi-\eta)\varphi_{\leq j}(\eta))d\eta  \right\|_2
	\end{align}
Using Lemma \ref{lem:ben} we have the following estimates 
	\begin{equation*}
		\begin{split}
			\eqref{6.15} +\eqref{6.17_n} & \leq \sum_{2^j\geq t^{-\beta}}\|P_j e^{-it\Delta^2}f\|_{10} \Big( \|P_{\leq j}e^{-it\Delta^2}(x^2f)\|_\frac{5}{2} + 2^{-j}\|P_{\leq j}e^{-it\Delta^2}(xf)\|_\frac{5}{2}\Big) \\ &\lesssim t^{-2+}\sum_{2^j\geq t^{-\beta}} 2^{-\frac{7}{2}j}\|f\|_X^2\lesssim t^{-2+\frac{7}{2}\beta+}\|f\|_X^2, 
		\end{split}
	\end{equation*}
	\begin{equation*}
		\begin{split}
			\eqref{6.16_n}\leq t  \sum_{2^j\geq t^{-\beta}} 2^{3j}\|P_j e^{-it\Delta^2}f\|_{\infty}\|P_{\leq j}e^{-it\Delta^2}(xf)\|_2\lesssim t^{-\frac{3}{2}+\frac{1}{4}+\frac{1}{64}+} \sum_{2^j\geq t^{-\beta}}  2^{-j}\|f\|_X^2\lesssim t^{-\frac{3}{2}+\frac{1}{4}+\frac{1}{64}+\beta+}\|f\|_X^2. 
		\end{split}
	\end{equation*}
	Letting $\beta=\frac{5}{32}$, we obtain $
		\|e^{-it\Delta^2}(xf)\|_{4}^2\lesssim t^{-\frac{69}{64}+}\|f\|_X^2$. 
If one uses $h$	instead of $f$ and $\beta = \f14$, then one obtains similarly 
\begin{align*}
 &I(xh)\leq \sum_{2^j\leq t^{-\frac{1}{4}}} \|P_je^{-it\Delta^2}(xh)\|_\frac{10}{3}\|P_{\leq j}e^{-it\Delta^2}(xh)\|_5 \lesssim \sum_{2^j\leq t^{-\frac{1}{4}}} 2^{\frac{11}{2}j}\|xh\|_\frac{5}{4}^2 \les t^{-\frac{11}{8}+\frac{1}{21}+}\|f\|_X^2 \\
 & II(xh) \leq \sum_{2^j> t^{-\frac{1}{4}}} \|P_je^{-it\Delta^2}(xh)\|_5\|P_{\leq j}e^{-it\Delta^2}(xh)\|_\frac{10}{3} 
			\les \sum_{2^j> t^{-\frac{1}{4}}} t^{-\frac{3}{2}}2^{-\frac{j}{2}}\|xh\|_\frac{5}{4}^2
			\les t^{-\frac{3}{2}+\frac{1}{8}+\frac{1}{21}+}\|f\|_X^2.
\end{align*}	
which gives \eqref{Decay_xh_4}. 

\end{proof}

\begin{lemma}
Assuming \eqref{Boot_est_g} and \eqref{Boot_est_h}, we have	\begin{equation}\label{Decay_x2h}
		\|x^2h\|_{L^2}\lesssim \log(t)(\|f\|_X^2+\|f\|_{X}^3+\|f\|_X^4)
	\end{equation}
\end{lemma}
\begin{proof}
Using \eqref{Decay_x2g}, we can straightforwardly estimate $\|x^2 h_3\|_{L^2} \lesssim \log(t)\big(\|f\|_X^2 + \|f\|_X^3\big)$. Moreover, we compute
	\begin{multline}
\partial_\xi^2(\widehat{h_1}(\xi,t))=\int_1^t \int \Big( \frac{q_2(\xi,\eta,s)}{A}\widehat{f}(\xi-\eta,s)+ \frac{q_3(\xi,\eta,s)}{A}\widehat{xf}(\xi-\eta,s) \Big)e^{is\varphi} \widehat{f}(\eta,s)d\eta ds \\
		+\int_1^t q_0(\xi,\eta) e^{is\varphi}\widehat{x^2f}(\xi-\eta,s)\widehat{f}(\eta,s)d\eta ds. \label{2h1}
	\end{multline}
Using \eqref{polyb} for the first term and \eqref{fracwuse} for the last two terms, we have 
\begin{align*}
\|x^2 h_1(x,t)\|_{L^2}\les \int_1^{t}  ( s^{-\f54} \|\la x \ra^{2+}f\|_{L^2} + s^{-1} \|xf\|_{L^2} + s^{-\f54} \|x^2f\|_{L^2} ) \|f\|_{X} \: ds  \les \log(t) \|f\|_{X}. 
\end{align*}

Next, we consider  $\partial_\xi^2(\widehat{h_2})$ to obtain the following terms 
	\begin{align}
		\label{x2h_2_p6} & \int_1^t\int s^2\frac{Q_6(\xi,\eta)}{A}e^{is\varphi}\widehat{f}(\xi-\eta,s)\partial_s\widehat{f}(\eta,s)d\eta \: ds\\
		\label{x2h_2_f}&\int_1^t \int s e^{is\varphi} \Big( \frac{q_2(\xi,\eta,s)}{A}\widehat{f}(\xi-\eta,s)+ \frac{q_3(\xi,\eta,s)}{A}\widehat{xf}(\xi-\eta,s) \Big) \partial_s\widehat{f}(\eta,s)d\eta \:  ds\\
		\label{x2h_2_x2f}&\int_1^t \int \frac{Q_0}{A}e^{is\varphi}\widehat{x^2f}(\xi-\eta,s)\partial_s\widehat{f}(\eta,s)d\eta ds
	\end{align}

Using similar norm estimates as in \eqref{fracwuse}, we obtain  \(\|\eqref{x2h_2_f}\|_{L^2} +\|\eqref{x2h_2_x2f}\|_{L^2}\les \|f\|_{X}^2\) . 
	 Therefore, we focus on \eqref{x2h_2_p6} where we first use \eqref{dt_hat_f} and  then  apply \eqref{key3} to obtain
	\begin{equation*}
		\begin{split}
			\eqref{x2h_2_p6}&=\int_1^ts^2\iint\frac{Q_6 (\xi,\eta)}{A}e^{is\psi}\widehat{f}(\xi-\eta,s)\widehat{f}(\eta-\sigma,s)\widehat{f}(\sigma,s)d\sigma d\eta ds
			\\&=\int_1^t s^2 \iint \frac{Q_6(\xi,\eta)}{A\: B}\left(\frac{1}{s}+\partial_s+\frac{Q}{s}\cdot \partial_\eta+\frac{S}{s}\cdot \partial_\sigma \right)e^{is\psi}\widehat{f}(\xi-\eta,s)\widehat{f}(\eta-\sigma,s)\widehat{f}(\sigma,s)d\sigma d\eta ds 
		\end{split}
	\end{equation*}
In the equality above, we have suppressed the variables in $\psi$, $Q$, $S$ and $B$ for brevity. Recall from \eqref{phase_3} that $\psi := \psi(\xi, \eta, \sigma)$, $Q(\xi,\eta) = 2 \xi- 3 \eta$, $S(\xi,\sigma) = \xi-3 \sigma$. Furthermore, by \eqref{key3}, we have 
\(
B = B_s := \frac{1}{s} + iY,
\)
where $Y = Y(\xi, \eta, \sigma)$ satisfies the lower bound 
\(
|Y| \gtrsim |\xi - \eta|^4 + |\eta - \sigma|^4 + |\sigma|^4.
\)

We decompose \eqref{x2h_2_p6} into:
\begin{align}
\eqref{x2h_2_p6}& = t^2\iint \frac{Q_6}{A \: B }e^{it\psi}\widehat{f}(\xi-\eta,t)\widehat{f}(\eta-\sigma,t)\widehat{f}(\sigma,t)d\sigma d\eta \label{6.18}  \\ 
& - \iint \frac{Q_6}{A \: B }e^{i\psi}\widehat{f}_1(\xi-\eta)\widehat{f}_1(\eta-\sigma)\widehat{f}_1(\sigma)d\sigma d\eta \\ 
 \label{6.19}&+\int_1^t \iint \frac{Q_6(s+q_0(\xi,\eta,\sigma,s))}{A \: B }e^{is\psi}\widehat{f}(\xi-\eta,s)\widehat{f}(\eta-\sigma,s)\widehat{f}(\sigma,s)d\sigma d\eta ds \\
 &+ \int_1^t s^2\iint \frac{Q_6}{A \: B } e^{is\psi} \partial_{s} \{\widehat{f}(\xi-\eta,s)\widehat{f}(\eta-\sigma,s)\widehat{f}(\sigma,s) \} d\sigma d\eta \: ds \label{6.20}\\
 & + \int_1^t s\iint \frac{Q_6 \: Q }{A \: B } e^{is\psi}\partial_{\eta} \{\widehat{f}(\xi-\eta,s)\widehat{f}(\eta-\sigma,s) \}\widehat{f}(\sigma,s) d\sigma d\eta \:  d s \label{6.21} \\ 
 &+ \int_1^t s\iint \frac{Q_6 \: S }{A \: B } e^{is\psi} \partial_{\sigma} \{\widehat{f}(\eta-\sigma,s) \widehat{f}(\sigma,s)\} \widehat{f}(\xi-\eta,s)d\sigma d\eta \:  d s \label{6.22}
\end{align}

 As the second terms in time independent its bound arises trivially by the fact that the smallness of the initial condition.  For the first and third term, we use Lemma \ref{flag_neg_deg} for $k=3,\omega=3,l=0$ to estimate 
\[ \| [T_{\frac{Q_6}{A\: B}}(u,u,u)](t,x) \|_{L_x^2} \les t^{\f14}\|\nabla_t^{-1}u\|_{\frac{5}{2}}\|u\|_{20}^2 \les t^{-2}\|f\|_X^2 \|f\|_\frac{5}{3}\lesssim t^{-2} \|f\|_X^3
  \] 
where we used $\|\nabla_t^{-1}u\|_{\frac{5}{2}} \les \|\nabla_t^{-\f12}f\|_{2} \les \|f\|_\frac{5}{3}$  in the third inequality. Using this estimate we immediatly obtain
\[  \|\eqref{6.18}\|_2+ \|\eqref{6.19}\|_2 \les \log(t) \|f\|_X^3. \] 
Using Lemma~\ref{flag_neg_deg} with the same $k,\omega,l$ together with \eqref{dt_hat_f}, we can also estimate 
\[
\|\eqref{6.20}\|_2\lesssim \int_1^t s^\frac{5}{2} \|u\|_8^2\|u^2\|_4ds\lesssim \int_1^t s^{-\frac{15}{4}+\frac{5}{2}}ds\|f\|_X^4\lesssim \|f\|_X^4 .\]

Finally, we focus on the terms in \eqref{6.21} and \eqref{6.22}. Here, we use Lemma~\ref{flag_neg_deg} with $k=3,\omega=4,l=0$ to see 
\[
\|\eqref{6.21}\|_2+\|\eqref{6.22}\|_2\les  \int_1^t s^\frac{5}{4}\|u\|_8^2\|e^{-it\Delta^2}(xf)\|_4ds\lesssim \int_1^t s^{-\frac{183}{160}+\varepsilon}ds\|f\|_X^3\lesssim \|f\|_X^3.
\]
This establishes the esimate $\|\partial_\xi^2(\widehat{h_2})\|_2\lesssim \log(t)(\|f\|_X^3+\|f\|_X^4)$. 

Lastly, we turn our attention to \(h_4\). If no \(\xi\)-derivative falls on \(\partial_{\eta}(e^{it\varphi})\) in \(h_4\), the resulting terms are similar to those in \(\partial_\xi^2(h_1)\). However, if at least one derivative acts on \(\partial_{\eta}(e^{it\varphi})\), the resulting terms take the following form:

\begin{align}
 \label{eth} & \int_1^t \int \frac{q_3(\xi,\eta,s) }{A}e^{is\varphi}\widehat{f}(\xi-\eta,s)\widehat{f}(\eta,s)d\eta ds\\ 
		\label{x2h_4_1}&\int_1^t s^2\int \frac{Q_1(\xi,\eta)\cdot \varphi_\eta |\varphi_\xi|^2}{A}e^{is\varphi}\widehat{f}(\xi-\eta,s)\widehat{f}(\eta,s)d\eta ds\\
		\label{x2h_4_2}&\int_1^t s\int \frac{Q_1(\xi,\eta)\cdot \varphi_\eta \varphi_\xi}{A}e^{is\varphi}\widehat{xf}(\xi-\eta,s)\widehat{f}(\eta,s)d\eta ds
	\end{align}

A term similar to \eqref{eth} appeared in \(\partial_\xi(h_1)\), and its \(L^2\) norm is bounded by \(\log(t) \|f\|_X^2\). We simplify \eqref{x2h_4_1} by integration by parts in $\eta$ as 
	\begin{equation*}
		\begin{split}
			\eqref{x2h_4_1}&=\int_1^t s\int \frac{q_6(\xi,\eta,s)}{A}e^{is\varphi}\widehat{f}(\xi-\eta,s)\widehat{f}(\eta,s)d\eta ds \\
			&+\int_1^t s\int\frac{Q_7(\xi,\eta)}{A}e^{is\varphi} \partial_{\eta}\{ \widehat{f}(\xi-\eta,s)\widehat{f}(\eta,s)\} d\eta ds
		\end{split}
	\end{equation*}
Recall that \(T_{\frac{Q_6}{A}}(u,u)\) previously appeared in \eqref{x2g_1}, where we established \(\|T_{\frac{Q_6}{A}}(u,u)(x,t)\|_{L^2_x} \lesssim t^{-2} \|f\|_X^2\). With the same argument, we can bound the $L^2$ norm of the first term in \eqref{x2h_4_1} by $\log(t)\|f\|_X^2$.On the other hand, noting that $\varphi_\eta \varphi_\xi = Q_6(\xi,\eta)$ the second term is the same as \eqref{x2h_4_2}. 

Thus, we finish the proof estimating $\|\eqref{x2h_4_2}\|_2$. To achieve this, we apply \eqref{key2} once more. After this application, the resulting terms are either of the form of the first term in \eqref{x2g_2}, the terms in \eqref{2h1}, the second term in \eqref{x2h_2_f}, or one of the following:
\begin{align}
		&\label{x2h_4_4}\int_1^t s\int \frac{q_3(\xi,\eta,s)}{A} e^{is\varphi}\partial_s(\widehat{xf}(\xi-\eta,s))\widehat{f}(\eta,s)d\eta ds\\
		&\label{x2h_4_5}\int_1^t \int q_0(\xi,\eta,s) e^{is\varphi}\widehat{xf}(\xi-\eta,s)\widehat{xf}(\eta,s)d\eta ds
	\end{align}
Since the \(L^2\) norms of the three previously obtained terms are bounded by \(\log(t) \|f\|_X^2\), we now focus on estimating the new terms listed above.

Having \eqref{decay_xf_4}, we immediately estimate  $\|\eqref{x2h_4_5}\|_2 \les \|f\|_X^2$. Moreover,  \eqref{x2h_4_4} can be written as 
	\begin{equation*}
		\int_1^t s\int \frac{q_3(\xi,\eta,s)}{A}e^{is\varphi}\widehat{f}(\xi-\eta,s)\partial_s\partial_\eta(\widehat{f}(\eta,s))d\eta ds 
	\end{equation*}
    Integration by parts in $\eta$, we either obtain a term as in  \eqref{x2h_2_f} or 
	\begin{equation*}
		\begin{split}
			\int_1^t s^2\int \frac{q_3(\xi,\eta,s) \varphi_{\eta}}{A}e^{is\varphi}\widehat{f}(\xi-\eta,s)\partial_s\widehat{f}(\eta,s)d\eta ds
		\end{split}
	\end{equation*}

Noting that $\varphi_\eta=Q_3$, this term the same as \eqref{x2h_2_p6}, and its \(L^2\) norm is bounded by \(\log(t)(\|f\|_X^2 + \|f\|_X^3 + \|f\|_X^4)\). This completes the analysis of all terms in \(h\).

\end{proof}

%\subsection{Estimates of $\|x^3h\|_{L^2}$}
\begin{lemma}
	Assuming \eqref{Boot_est_g} and \eqref{Boot_est_h}, we have
	\begin{equation}\label{Decay_x3h1}
		\|x^3h_1\|_{L^2}\lesssim t^{\frac{1}{24}}\|f\|_X^3
	\end{equation}
\end{lemma}
\begin{proof}
To study $\partial_\xi^3(\widehat{h_1})$, we apply the decomposition illustrated in \eqref{Decom_Bi}: 
 $f= f_1 + f_*+ g + h$, and reduce the problem to estimating the following terms:
\begin{align}
	\label{h1_h_h}&\partial_\xi^3 \int_1^t s^{-2}\int\frac{1}{A^2}e^{is\varphi}\widehat{h}(\xi-\eta,s)\widehat{h}(\eta,s)d\eta ds\\
	\label{h1_f_g}&\partial_\xi^3 \int_1^t s^{-2}\int\frac{1}{A^2}e^{is\varphi}\widehat{f}(\xi-\eta,s)\widehat{g}(\eta,s)d\eta ds\\
	\label{h1_g_h}&\partial_\xi^3 \int_1^t s^{-2}\int\frac{1}{A^2}e^{is\varphi}\widehat{g}(\xi-\eta,s)\widehat{h}(\eta,s)d\eta ds
\end{align}
We begin with \eqref{h1_h_h}. Upon computing the third derivative, we obtain

\begin{align}
\eqref{h1_h_h}&= \int_1^t \int \Big(  \frac{q_0}{A \: s} \widehat{x^3h}(\xi-\eta,s)+ \frac{q_3}{A} \widehat{x^2h}(\xi-\eta,s) + \frac{q_1}{A} \widehat{h}(\xi-\eta,s) \Big) e^{is\varphi}\widehat{h}(\eta,s)d\eta ds \label{h1_1_1} \\ 
& + \int_1^t  \int \frac{q_2}{A}  e^{is\varphi} \widehat{xh}(\xi-\eta,s)\widehat{h}(\eta,s)d\eta ds +  \int_1^t s\int \frac{q_5}{A} e^{is\varphi} \widehat{h}(\xi-\eta,s)\widehat{h}(\eta,s)d\eta ds \label{h1_1_5}
\end{align} 
Using \eqref{m1b} with $ \alpha=\f 1{10}$, and \eqref{fracwuse} we have $\|\eqref{h1_1_1} \|_2 \les t^{0+}\|f\|_X^2$. We use normalization similar to \eqref{m6t} to estimate the $L^2$ norm of the second term in \eqref{h1_1_5} by 
\begin{equation*}
	\begin{split}
	& \int_1^t s\sum_j 2^j \min\left(\|P_j(e^{-is\Delta^2}h)\|_{10}\|P_\leq j(e^{-is\Delta^2}h)\|_\frac{5}{2},\|P_j(e^{-is\Delta^2}h)\|_{\infty}\|P_\leq j(e^{-is\Delta^2}h)\|_2 \right)ds\\
		&\lesssim \int_1^t s\sum_j 2^j\min\left(s^{-2}2^{-j}\|P_jh\|_\frac{10}{9}, s^{-\frac{5}{2}}2^{-3j}\|h\|_\frac{10}{9} \right)\|h\|_1 ds\\
		&\lesssim \int_1^t \sum_j \min(s^{-1+\frac{1}{16}+}2^{\frac{j}{2}}, s^{-\frac{3}{2}+\frac{1}{32}+}2^{-2j})ds\|f\|_X^2 \les \|f\|_X^2.
	\end{split}
\end{equation*}
We  decompose the second term in \eqref{h1_1_5} as 
\begin{align*}
I + II := \int_1^t \int \frac{q_2}{A}e^{is\varphi} \widehat{xh}(\xi-\eta,s)\widehat{h}(\eta,s) \Psi_1(\xi,\eta) d\eta ds + \int_1^t \int \frac{q_2}{A}e^{is\varphi} \widehat{xh}(\xi-\eta,s)\widehat{h}(\eta,s) \Psi_2(\xi,\eta) d\eta ds
\end{align*}
Similar to the proof of Corollary~\ref{CM_neg_deg}, we let $m_t(\xi,\eta)= A^{-1} m_2(\xi,\eta,s) (\frac{1}{s}+|\xi-\eta|^4)^\frac{1}{2}\Psi_1(\xi,\eta) $ and  estimate 
\begin{equation*}
    \|I\|_{2} 
        \leq\int_1^t \left\|T_m(\nabla_s^{-2}e^{-is\Delta^2}(xh),e^{-is\Delta^2}h )\right\|ds\lesssim \int_1^t s^{-1}\|xh\|_\frac{10}{7}\|f\|_Xds\lesssim t^{0+}\|f\|_X^2.
    \end{equation*}
We normalize the multiplier in $II$ as 
$$
m(\xi,\eta)= \frac{q_2(\xi,\eta) \Psi_2(\xi,\eta)\psi_j(\xi-\eta)}{A \: \min(s,2^{-4j})2^{2j}}
$$
and use $ \|T_m(P_{\leq j}e^{-is\Delta^2}(xh),P_j e^{-is\Delta^2}h)\|_2 \les \|P_{\leq j}e^{-is\Delta^2}(xh)\|_\frac{5}{2}\|P_j e^{-is\Delta^2}h\|_{10} $ to estimate 
\begin{align*}
     \| II\|_2&
       \lesssim \int_1^t \sum_j 2^{\frac{3j}{2}} \min(s, 2^{-4j})s^{-2} \|x^{2+}h\|_2 \|h\|_2 \\
       & \les \int_1^t \sum_j\min(s^{-1+}2^{\frac{3}{2}j}, s^{-2+}2^{-\frac{5}{2}j})\|f\|_X^2ds\lesssim t^{-\frac{3}{8}+}\|f\|_X^2
\end{align*}
where we used  $\|P_{\leq j}e^{-is\Delta^2}(xh)\|_{\frac{5}{2}} \les  2^{\frac{3j}{2}}\|x^2h\|_2$ and $\|P_j e^{-is\Delta^2}h\|_{10} \les 2^{-2j}s^{-2}\|h\|_2 $ in the first inequality.

Next, we focus on  \eqref{h1_f_g}. To begin, we express \eqref{h1_f_g} in a trilinear form using \eqref{eqn_g}. Then we apply \( \partial_\xi^3 \) to compute 

\begin{align}
\eqref{h1_f_g}=\label{h1_2_1}&\int_1^t s\iint \frac{q_5(\xi,\eta)}{AC}e^{is\psi} \widehat{f}(\xi-\eta,s)\widehat{f}(\eta-\sigma,s)\widehat{f}(\sigma,s)d\sigma d\eta ds\\
\label{h1_2_11}+&\int_1^t \iint \frac{q_1(\xi,\eta)}{AC}e^{is\psi} \widehat{f}(\xi-\eta,s)\widehat{f}(\eta-\sigma,s)\widehat{f}(\sigma,s)d\sigma d\eta ds\\
	+&\label{h1_2_2}\int_1^t  \iint s\frac{q_2(\xi,\eta)}{C}e^{is\psi} \widehat{xf}(\xi-\eta,s)\widehat{f}(\eta-\sigma,s)\widehat{f}(\sigma,s)d\sigma d\eta ds \\
	+&\label{h1_2_4}\int_1^t s\iint\frac{q_3(\xi,\eta)}{C}e^{is\psi} \widehat{x^2f}(\xi-\eta,s)\widehat{f}(\eta-\sigma,s)\widehat{f}(\sigma,s)d\sigma d\eta ds \\
	+&\label{h1_2_5}\int_1^t s^{-1} \iint \frac{q_0(\xi,\eta)}{ A C}e^{is\psi} \widehat{x^3f}(\xi-\eta,s)\widehat{f}(\eta-\sigma,s)\widehat{f}(\sigma,s)d\sigma d\eta ds
\end{align}
where $ C:=C(\eta,\sigma) = \frac{1}{s} + Z(\eta,\sigma)$. 

We use \eqref{key3} to reduce estimating $\|\eqref{h1_2_1}\|_2$ to estimating the $L^2$ norm of the following type of terms
\begin{align}
   & \label{h1_2_1_1}\int_1^t  \iint \frac{q_5(\xi,\eta)q_0(\xi,\eta,\sigma)}{ACB}     e^{is\psi} \widehat{f}(\xi-\eta,s)\widehat{f}(\eta-\sigma,s)\widehat{f}(\sigma,s)d\sigma d\eta ds \\
    &\label{h1_2_1_2}t \iint \frac{q_5(\xi,\eta)}{ACB}     e^{it\psi} \widehat{f}(\xi-\eta,t)\widehat{f}(\eta-\sigma,t)\widehat{f}(\sigma,t)d\sigma d\eta\\
    &\label{h1_2_1_3}\int_1^t s\iint \frac{q_5(\xi,\eta)}{ACB}     e^{is\psi}\partial_s\{ \widehat{f}(\xi-\eta,s)\widehat{f}(\eta-\sigma,s)\widehat{f}(\sigma,s)\}d\sigma d\eta ds\\
    &\label{h1_2_1_4}\int_1^t \iint \frac{q_2(\xi,\eta)}{A}\frac{1}{C}\frac{q_4(\xi,\eta,\sigma)}{B} e^{is\psi} \widehat{f}(\xi-\eta,s)\widehat{xf}(\eta-\sigma,s)\widehat{f}(\sigma,s)d\sigma d\eta ds + \text{similar terms}
\end{align}
Here, “similar terms” refers to expressions where the spatial variable \( x \) appears in product with other  \( f \). Using Corollary  \ref{flag_neg_deg} we have the following estimate  
\begin{align*}
\left\| \iint \frac{q_2}{A}\frac{1}{C}\frac{q_3}{B} e^{is\psi} \widehat{f}(\xi-\eta,s)\widehat{f}(\eta-\sigma,s)\widehat{f}(\sigma,s)d\sigma d\eta \right\|_2 
&\lesssim  s^{\frac{5}{4} }\|u\|_{20}^2\|\nabla_s^{-2}u\|_\frac{5}{2} \les s^{-1} \|f\|_\frac{5}{4}\|f\|_X^2. 
\end{align*}
Therefore, $\|\eqref{h1_2_1_1}\|_2 + \|\eqref{h1_2_1_2}\|_2 \les t^{0+}\|f\|_X^2 $. The same estimate also bounds $\| \eqref{h1_2_11}\|_2$ taking $k=1$, $l=0$ in Corollary~\ref{flag_neg_deg}. Moreover, $\|\eqref{h1_2_1_3}\|_2$ can be estimated as
\begin{equation} 
	\|\eqref{h1_2_1_3}\|_2\lesssim \int_1^t s^{\frac{11}{4}}\|u\|_8^2\|u^2\|_4ds\lesssim \|f\|_X^4\int_1^t s^{-1} ds\lesssim \log(t)\|f\|_X^4.\label{forh}
\end{equation}
To estimate $\eqref{h1_2_1_4}$, we first note that for $\alpha=1,2$ we have 
\begin{align}
		\|e^{-is\Delta^2}&(x^{\alpha}f)\|_\frac{40}{19}\leq \sum_{2^j<s^{-1/3}} \|P_j e^{-is\Delta^2}(x^{\alpha} f)\|_\frac{40}{19}+\sum_{2^j\geq s^{-\frac{1}{3}}}\|P_je^{-is\Delta^2}(x^{\alpha}f)\|_\frac{40}{19} \nn \\
		&\lesssim\sum_{2^j<s^{-1/3}}2^{\frac{j}{8}}\|x^{\alpha}f\|_2+\sum_{2^j\geq s^{-\frac{1}{3}}}s^{-\frac{1}{8}}2^{-\frac{j}{4}}\|x^{\alpha}f\|_\frac{40}{21} 
\lesssim s^{-\frac{1}{24}+ (\alpha -1) \frac{47}{720}+} \|f\|_X \label{Decay_xf}
\end{align}
Therefore, we have 
\begin{align}
	\|\eqref{h1_2_1_4}\|_2 & \lesssim \int_1^t s^\frac{3}{2}\|u\|_{80}^2\|e^{-is\Delta^2}(xf)\|_\frac{40}{19}
	+s\|u\|_{10}^2\|\nabla_s^{-2} e^{-is\Delta^2}(xf)\|_{\frac{10}{3}}ds   \nn \\ & \lesssim\int_1^t s^{\frac{3}{2}-\frac{1}{24}-\frac{39}{16}+}\|f\|_X^3+ s^{-1}\|xf\|_{\frac{10}{7}}\|f\|_X^2ds\lesssim t^{\frac{1}{48}+}\|f\|_X^3\label{6.44est}
	\end{align}
which finalizes the estimate $\|\eqref{h1_2_1}\|_2\lesssim t^{\frac{1}{48}+}\|f\|_X^3$. We note that $\|\eqref{h1_2_2}\|_2$ can be also estimated using $\eqref{Decay_xf}$ for $\alpha=1$ and \eqref{6.44est}.

To bound $\|\eqref{h1_2_4}\|_2$, we  use \eqref{Decay_xf} with $\alpha=2$ and have 
\begin{align*}
	\|\eqref{h1_2_4}\|_2\les \int_1^t s^{\frac{5}{4}}\|e^{-is\Delta^2}(x^2f)\|_\frac{40}{19}\|u\|_{80}^2ds\les \|f\|_X^3 \int_1^t s^{\frac{5}{4}-\frac{39}{16}+\frac{17}{720}}ds\lesssim t^{\frac{7}{720}}\|f\|_X^3.
\end{align*}
Next, we focus on the final term \eqref{h1_2_5}. Since Theorem~\ref{th:flagcf} is not applicable for \( p_1 = 2 \), we employ the same strategy used for \( H_{j}(\widehat{f},\widehat{g})(\xi,t) \) in \eqref{x2g_1} to estimate this term. Specifically, using cut-off functions, we decompose \eqref{h1_2_5} into the following components: 
\begin{align*}
F_{0,i}(\widehat{x^3f},\widehat{f},\widehat{f})(\xi,t)=\int_1^t  \iint \frac{q_0(\xi,\eta,s)}{C}\phi_i(\xi,\eta,\sigma) e^{is\psi} \widehat{x^3f}(\xi-\eta,s)\widehat{f}(\eta-\sigma,s)\widehat{f}(\sigma,s)d\sigma d\eta ds
\end{align*}
for $i=1,2,3$. We let 
$$
m^*(\xi,\eta,\sigma)=m_1(\xi,\eta) m_2(\eta,\sigma) = q_0(\xi,\eta,s) \frac{\phi_1(\xi, \eta,\sigma) \varphi_{j}(\xi-\eta) }{C\min(s,2^{-4j})} 
$$
and observe that $m^*(\xi,\eta,\sigma)$ defines a Coifman-Meyer operator with flag singularity whose operator norm is independent of $s$. Therefore, we have the estimate
\begin{equation*}
	\begin{split}
		\|F_{0,1}(\widehat{x^3f},\widehat{f},\widehat{f})\|_2&\leq \int_1^t \left\|\sum_j \min(s,2^{-4j})T_{m^*}(P_j(e^{-is\Delta^2}(x^3f), P_ku, u)\right\|_2ds\\ 
	&\les  \int_1^t \sum_j  \min(s, 2^{-4j}) \|P_j (e^{-is\Delta^2}(x^3f)\|_\frac{40}{19}\|u\|_{80}^2 \: ds \\ 
		&\lesssim \int_1^t \sum_j 
		s^{-\frac{39}{16}} 2^{\frac{j}{8}} \min (s, 2^{-4j})ds \|f\|_X^3\lesssim \|f\|_X^3
	\end{split}
\end{equation*}
We next estimate $\|F_{0,2}(\widehat{x^3f},\widehat{f},\widehat{f})\|_{2}$.  We let
\begin{align}
m^*(\xi,\eta,\sigma)=m_1^*(\xi,\eta) m^*_2(\eta,\sigma) = q_0(\xi,\eta,s) \: \frac{\phi_2(\xi, \eta,\sigma) \psi_{k}(\eta-\sigma) }{C\min(s,2^{-4k})} \label{m*F21}
\end{align}
and localize with $\varphi_j(\xi-\eta)\psi_k(\eta-\sigma) $ to see 
\begin{align*}
F_{0,2}(\widehat{x^3f},\widehat{f})(\xi,t)= 			\int_1^t \sum_k \sum_{j\leq k+10} \min(s, 2^{-4k}) T_{m^*}(P_j (e^{-is\Delta^2}(x^3f), P_ku,  u)ds.
\end{align*}
We have $\|T_{m^*}(P_j (e^{-is\Delta^2}(x^3f), P_ku,  u)\| \les \|P_j (e^{-is\Delta^2}(x^3f)\|_\frac{80}{39}\|u\|_{160}^2 \les 2^{\frac{j}{16}} \|x^3f\|_2 \|u\|_{160}^2$, and therefore, 
\begin{equation*}
		\|F_{0,2}(\widehat{x^3f},\widehat{f},\widehat{f})\|_2\lesssim\int_1^t s^{-\frac{79}{32}}\sum_k \min(s, 2^{-4k})2^{\frac{k}{16}+}\|f\|_X^3 ds \les t^{\frac{1}{24}}\|f\|_X^3 
	\end{equation*}
    
To estimate $\|F_{0,3}(\widehat{x^3f},\widehat{f},\widehat{f})\|_2$, we localize with $\varphi_j(\xi-\eta)\psi_k(\sigma) $ and use $m^*$ as in \eqref{m*F21} with $\psi_l(\eta-\sigma)$ replaced by \( \psi_l(\sigma) \). Then the bound $t^{\frac{1}{24}}\|f\|_X^3$ follows identically to $\|F_{0,2}(\widehat{x^3f},\widehat{f},\widehat{f})\|_2$.

We next focus on estimating \( \|\eqref{h1_g_h}\|_2 \). We note that the time growth of the weighted norms of \( h(x,t) \) is slower than that of the weighted norms of \( f(x,t) \). Therefore, much of the argument follows analogously to the estimate of \( \|\eqref{h1_f_g}\|_2 \) after a change of variables and by Lemma~\ref{dsh}, see Remark~\ref{rmk:dsh}. 

\begin{remark}\label{rmk:dsh}
 Estimating the term \( \|\eqref{h1_g_h}\|_2\), one sees that the proof only breaks down at the point where we apply \eqref{key2}, since in this case we obtain terms in which the \( s \)-derivative falls on the first  \( \widehat{f} \) in \eqref{h1_2_1_3}. For \( f \), we handled this by expressing \( \partial_s \widehat{f} \) using \eqref{dt_hat_f}. This representation, however, is not valid for \(\partial_s h \). To overcome this, we prove $\| (e^{-is|\xi|^2}(\partial_s \widehat{h}(\xi,s)))^\vee\|_p=\|e^{is\Delta^2}(\partial_s h)\|_p$ satisfies the same bound that $\|u^2\|_p$ satisfies in Lemma~\ref{dsh}. The estimate in \eqref{forh} then completes the proof of the estimate \eqref{Decay_x3h1}. 
\end{remark}

\end{proof}

\begin{lemma}
Assuming \eqref{Boot_est_g} and \eqref{Boot_est_h}, we have
	\begin{equation}
		\|x^3h_2\|_2\lesssim t^{\frac{1}{31}}(\|f\|_X^3+\|f\|_X^4)
	\end{equation} \label{Decay_x3h2} 
\end{lemma}
\begin{proof}
	Applying $\partial_\xi^3$ to $\widehat{h_2}$, we obtain the terms: 

	\begin{align}
		\label{h2_1}&\int_1^t s^3 \int \frac{Q_9}{A}e^{is\varphi} \widehat{f}(\xi-\eta,s)\partial_s\widehat{f}(\eta,s)d\eta ds\\
		\label{h2_2}+&\int_1^t s^2 \int \frac{Q_6}{A}e^{is\varphi} \widehat{xf}(\xi-\eta,s)\partial_s\widehat{f}(\eta,s)d\eta ds\\
		\label{h2_3}+&\int_1^t s^2 \int \frac{q_5}{A}e^{is\varphi} \widehat{f}(\xi-\eta,s)\partial_s\widehat{f}(\eta,s)d\eta ds \\
		\label{h2_4}+&\int_1^t s \int e^{is\varphi} \Big( \frac{q_3}{A} \widehat{x^2f}(\xi-\eta,s)+  \frac{q_2}{A}\widehat{xf}(\xi-\eta,s) \Big) \partial_s\widehat{f}(\eta,s)d\eta ds \\
		\label{h2_6}+&\int_1^t e^{is\varphi}  \int \Big( \frac{s\: q_1}{A}\widehat{f}(\xi-\eta,s)+  \frac{q_0}{A}\widehat{x^3f}(\xi-\eta,s)\Big) \partial_s\widehat{f}(\eta,s)d\eta ds
	\end{align}
Using \eqref{fracwuse} and \eqref{dt_hat_f}, we can easily see that $\|\eqref{h2_4}\|_2 \les \log(t)\|f\|_X^3 $. Moreover, using \eqref{m1b} for $\alpha=\frac{1}{10}$ and \eqref{fracwuse} together with \eqref{dt_hat_f}, we obtain $ \|\eqref{h2_6}\|_2 \les t^{\f1{47}}\|f\|_X^3$. 

We write  \eqref{h2_3} as a trilinear term using \eqref{eqn_g} and decomposing as 
	%\begin{equation*}
	%\begin{split}
	%	\eqref{h2_3}&=\int_1^t s^2 \iint \frac{q_5}{A}e^{is\psi} \widehat{f}(\xi-\eta,s)\widehat{f}(\eta-\sigma,s)\widehat{f}(\sigma,s)d\sigma d\eta ds
	%\end{split}
	%\end{equation*}

 %   To estimate this, we decompose $\eqref{h2_3}$ into 
    
\begin{align}
& F_{5,1}(\widehat{f},\widehat{f},\widehat{f}):= \int_1^t s^2 \iint \frac{q_5}{A}\Psi_1(\xi,\eta)  e^{is\psi} \widehat{f}(\xi-\eta,s)\widehat{f}(\eta-\sigma,s)\widehat{f}(\sigma,s)d\sigma d\eta ds \label{defF1}\\
& F_{5,2,j}(\widehat{f},\widehat{f},\widehat{f}):=\int_1^t s^2 \iint \frac{q_5}{A}\Psi_2(\xi,\eta)\Psi_j(\eta,\sigma)e^{is\psi} \widehat{f}(\xi-\eta,s)\widehat{f}(\eta-\sigma,s)\widehat{f}(\sigma,s)d\sigma d\eta ds\label{defF2}
\end{align}
    \begin{remark}
        Here, we observe that when considering \( \Psi_1(\xi,\eta) \), there is no need for an additional cut-off function \( \Psi_j(\eta, \sigma) \). Specifically, since the multiplier depends only on \( \xi \) and \( \eta \) within the support where \( |\eta| \lesssim |\eta-\xi| \), localizing in \( \eta-\xi \) is sufficient.

On the other hand, when considering the cut-off \( \Psi_2(\xi,\eta) \), which is supported in \( |\eta-\xi| \lesssim |\eta| \), localization in \( \eta \) alone is insufficient. In this case, we additionally introduce \( \Psi_j(\eta, \sigma) \).This gives us the control \( |\eta| \lesssim |\eta -\sigma| \) when \( j=1 \) and \( |\eta| \lesssim |\sigma| \) when \( j=2 \).
 
\end{remark}

To use Theorem~\ref{th:flagcf}, we let 
$$m^*(\xi,\eta,s)=m_1(\xi,\eta)= \frac{q_5(\xi,\eta,s)}{A\:2^{5j}\min(s, 2^{-4j})}\psi_j(\xi-\eta)\Psi_1(\xi,\eta)$$
and notice that $m^*(\xi,\eta,s)$ is a Coifman-Meyer multiplier with flag singularity as defined in \eqref{defm*}, and its norm is independent of $s$. 
Therefore, using $$\|T_{m^*}(P_ju, u,u)\|_2 \les \min(\|P_ju\|_\frac{10}{3}\|u\|_{10}^2, \|P_ju\|_{10}\|u\|_5^2)$$ we estimate
\begin{align} 
\|F_{5,1}(\widehat{f},\widehat{f},\widehat{f})\|_2&\lesssim \int_1^t s^2\sum_j 2^j \|T_m^*(P_ju,u,u)\|_2\nn  \\  
		&\lesssim \int_1^t \sum_j \min(2^{4j}\|f\|_\frac{10}{9}, s^{-1}2^{-j}\|P_jf\|_\frac{10}{7}, s^{-\frac{3}{2}}2^{-3j}\|f\|_\frac{10}{9}  )\|f\|_X^2ds \nn  \\ 
		&\lesssim\int_1^t \sum_j \min(2^{4j}, s^{-1}, s^{-\frac{3}{2}}2^{-3j}  )\|x^{2+}f\|_2\|f\|_X^2ds \lesssim t^{0+}\|f\|_X^3 \label{s22j}
	\end{align}
We localize $F_{5,2,1}$ with $\varphi_j(\eta)\varphi_k(\eta-\sigma) $ and let 
\begin{align} m^*(\xi,\eta,s) = m_1(\xi,\eta)=\frac{q_5(\xi,\eta,s)}{A\:2^{5j}\min(s, 2^{-4j})}\varphi_j(\eta)\Psi_2(\xi,\eta) \label{m*foreta}
\end{align}
to obtain  
\begin{align}
		\|F_{5,2,1}(\widehat{f},\widehat{f},\widehat{f})\|_2&\leq\int_1^t s^2\sum_{j\leq k+10}\sum_k 2^{5j}\min(s, 2^{-4j})\|T_{m^*}(u, P_ku,u)\|_2ds \lesssim t^{0+}\|f\|_X^3. \label{s22k}
\end{align}
we used the same norm estimates as in \eqref{s22j} in the last inequality. Moreover, $\|F_{5,2,2}(\widehat{f},\widehat{f},\widehat{f})\|_2$ can be estimated similarly localizing around $\eta$ and $\sigma$ with $\varphi_j(\eta)\varphi_k(\sigma)$, and   using $m^*$ as in \eqref{m*foreta}. 

We estimate \( \|\eqref{h2_2}\|_2 \) using the same decomposition as in \eqref{defF1} and \eqref{defF2}.  Consequently, we need to study $F_{6,1}(\widehat{xf},\widehat{f}, \widehat{f})$ and $F_{6,2,j}(\widehat{xf},\widehat{f}, \widehat{f})$ for $j=1,2$. 
Localizing with $\psi_j(\xi-\eta)$, we see that 
\begin{align}
\|F_{6,1}&(\widehat{xf},\widehat{f}, \widehat{f})\|_{2}  \les \int_1^t s^2\sum_j 2^{6j}\min(s,2^{-4j})\|T_{m^*}(P_je^{-is\Delta^2}(xf),u,u)\|_2ds \nn \\ 
&\lesssim \int_1^t s^2\sum_j 2^{2j}\min(\|P_je^{-is\Delta^2}(xf)\|_\frac{10}{3}\|u\|_{10}^2, \|P_je^{-is\Delta^2}(xf)\|_5\|u\|_{10}\|u\|_5 )ds \nn \\
&\lesssim \int_1^t \sum_j(2^{4j}\|xf\|_\frac{10}{7}, s^{-1}\|xf\|_\frac{10}{7}, s^{-\frac{5}{4}}2^{-j}\|xf\|_{\frac{5}{4}})\|f\|_X^2ds \lesssim t^{0+} \|f\|_X^3 \label{F61}
\end{align}
Localizing with $\varphi(\xi-\eta) \varphi_k(\eta) \psi_l^2(\eta-\sigma)$, we obtain 

\begin{align}
\|F_{6,2,1}&(\widehat{xf},\widehat{f}, \widehat{f})\|_{2}   \les \int_1^t s^2\sum_l\sum_{k\leq l+10}\sum_{j\leq k+10}2^{6k}\min(s,2^{-4k})\|T_{m^*}(P_je^{-is\Delta^2}(xf),P_lu,u)\|_2ds \nn \\
& \lesssim\int_1^t s^2\sum_l\sum_{k\leq l+10}\sum_{j\leq k+10} 2^{2k}\|P_je^{-is\Delta^2}(xf)\|_\frac{10}{3}\|P_lu\|_{10}\|u\|_{10}ds \nn \\
		&\lesssim \int_1^t \sum_l \min(s^{-1}2^{-l}\|f\|_\frac{10}{9}\|xf\|_2, s^{-1}\|f\|_\frac{10}{9}\|xf\|_\frac{10}{7}, s2^{9l} \|f\|_\frac{10}{9}\|xf\|_\frac{10}{7})\|f\|_Xds \nn\\ 
        &\lesssim t^{0+} \|f\|_X^3 \label{F621}
\end{align}
Finally, localizing with $\varphi(\xi-\eta) \varphi_k(\eta) \psi_l^2(\sigma)$, we can obtain $ \|F_{6,2,2}(\widehat{xf},\widehat{f}, \widehat{f})\|_{2} \les t^{0+} \|f\|_X^3$ analogously to  $\|F_{6,2,1}(\widehat{xf},\widehat{f}, \widehat{f})\|_{2}$. 

To estimate the final term $\|\eqref{h2_1}\|_{2}$, we apply \eqref{key3} and reduce the problem estimating the following terms:
\begin{align}
	\label{h2_1_1}&\int_1^t s^2\iint \frac{Q_5(\xi,\eta) q_0(\xi,\eta,\sigma,s)}{A}e^{is\psi}\widehat{f}(\xi-\eta,s)\widehat{f}(\eta-\sigma,s)\widehat{f}(\sigma,s)d\sigma d\eta ds\\
	\label{h2_1_2}&t^3\iint \frac{Q_5q_0(\xi,\eta,\sigma,s)}{A}e^{it\psi}\widehat{f}(\xi-\eta,t)\widehat{f}(\eta-\sigma,t)\widehat{f}(\sigma,t)d\sigma d\eta\\
	\label{h2_1_3}&\int_1^t s\iint \Big( \frac{Q_9}{A^2B} + \frac{Q_9}{AB^2} \Big)e^{is\psi}\widehat{f}(\xi-\eta,s)\widehat{f}(\eta-\sigma,s)\widehat{f}(\sigma,s)d\sigma d\eta ds\\
	\label{h2_1_4}&\int_1^t s^3\iint \frac{Q_9}{AB}e^{is\psi} \partial_s\{\widehat{f}(\xi-\eta,s)\widehat{f}(\eta-\sigma,s)\widehat{f}(\sigma,s)\}d\sigma d\eta ds\\ 
	\label{h2_1_5}&\int_1^t s^2\iint \frac{Q_{9}Q}{AB}e^{is\psi}\widehat{xf}(\xi-\eta,s)\widehat{f}(\eta-\sigma,s)\widehat{f}(\sigma,s)d\sigma d\eta ds + \text{similar terms} 
	\end{align}
Here by similar terms we mean $x$ can be multiplied by any other $f$ in the integral. Notice that letting $k=6$, $w=3$ in Corollary~\ref{flag_neg_deg}, we can estimate \eqref{h2_1_3} directly.  %\textcolor{red}{ I cant follow what is going on here, I dont think we need this term as it boils to \eqref{h2_1_1}, please check that}
%\hrule
%\textcolor{blue}{I agree!}
%\hrule
%\begin{equation*}
%	\begin{split}
	%\|\eqref{h2_1_3}\|_2\lesssim \int_1^t s^{\frac{3}{2}}\|u\|_{160}^2\|\nabla_s^{-1}e^{-is\Delta^2}(xf)\|_\frac{80}{39}ds\lesssim \int_1^t s^{-\frac{31}{32}+\varepsilon}ds\|f\|_X^3\lesssim t^{\frac{1}{31}}\|f\|_X^3
	%\end{split}
%\end{equation*}

 Noting the similarity with the term \eqref{h2_3}, we can decompose \eqref{h2_1_1} to 
\begin{align*}
&F_{5,1}(\widehat{f},\widehat{f},\widehat{f})= \int_1^t s^2\iint \frac{Q_5q_0}{A} \Psi_1(\xi,\eta)  e^{is\psi}\widehat{f}(\xi-\eta,s)\widehat{f}(\eta-\sigma,s)\widehat{f}(\sigma,s)d\sigma d\eta ds \\
&F_{5,2,i}(\widehat{f},\widehat{f},\widehat{f})=\int_1^t s^2\iint \frac{Q_5q_0}{A} \Psi_2(\xi,\eta )  \Psi_i(\eta,\sigma)e^{is\psi}\widehat{f}(\xi-\eta,s)\widehat{f}(\eta-\sigma,s)\widehat{f}(\sigma,s)d\sigma d\eta ds
\end{align*}
for $i=1,2$, where we have suppressed the variables in $Q$ and $q_0$ to shorten the notation. However, their dependence can be traced back using \eqref{h2_1_1}.
We consider the Coifman-Meyer multiplier with flag singularity 
\begin{align} \label{m*9}
 m^*(\xi,\eta,\sigma,s)=m_1(\xi,\eta) m_3(\xi,\eta,\sigma) = \frac{Q_5\Psi_1(\xi,\eta) \psi_j(\xi-\eta) }{A\: \min(s,2^{-4j}) 2^{5j}} \: q_0(\xi,\eta,\sigma, s)
     \end{align}
and localize with $\psi_j(\xi-\eta)$ to estimate $\|F_{5,1}(\widehat{f},\widehat{f},\widehat{f})\|_2 \leq t^{0+} \|f\|_X^3$  similar to \eqref{s22j}. 
To estimate $\|F_{5,2,j}(\widehat{f},\widehat{f},\widehat{f})\|_2$ for $j=1,2$, we consider 
\begin{align}
& m^*(\xi,\eta,\sigma,s)=m_1(\xi,\eta) m_3(\xi,\eta,\sigma)= \frac{Q_5\psi_j(\eta)\Psi_2(\xi,\eta)}{A\:2^{5j} \:\min(s,2^{-4j})) } \: q_0 \Psi_1(\eta, \sigma) \label{m*91}\\
& m^*(\xi,\eta,\sigma,s)=m_1(\xi,\eta) m_3(\xi,\eta)= \frac{Q_5 \psi_j(\eta)\Psi_2(\xi,\eta)}{A\:2^{5j} \:\min(s,2^{-4j})) } \: q_0 \Psi_2(\eta,\sigma)  \label{m*92}
\end{align}
and the localizations $\psi_j(\eta)\varphi_k(\eta-\sigma)$ and $ \psi_j(\eta) \varphi_k(\sigma)$ respectively to obtain $\|F_{5,2,i}(\widehat{f},\widehat{f},\widehat{f})\|_2\lesssim  t^{0+}\|f\|_X^3 $ by \eqref{s22k}. We note that $\|\eqref{h2_1_2}\|_2$ can be bounded in a similar way. 

We note that one has $ \frac{Q_4(\xi,\eta)}{B} = q_0(\xi,\eta,\sigma)$ and therefore, studying \eqref{h2_1_4} reduces to estimating $\|F_{5,1}(s\partial_s\widehat{f},\widehat{f},\widehat{f})\|_2$, $\|F_{5,1}(s\widehat{f},\partial_s\widehat{f},\widehat{f})\|_2$,$\|F_{5,1}(s\widehat{f},\widehat{f},\widehat{\partial_sf})\|_2$ and $\|F_{5,2,i}(s\partial_s\widehat{f},\widehat{f},\widehat{f})\|_2$, $\|F_{5,2,i}(s\widehat{f},\partial_s\widehat{f},\widehat{f})\|_2$, $\|F_{5,2,i}(s\widehat{f},\widehat{f},\partial_s\widehat{f})\|_2 $. Using \eqref{dt_hat_f} and the multiplier \eqref{m*9} with localization $\psi_j(\xi-\eta)$ we can estimate 
\begin{align}
\|F_{5,1}(s\widehat{f},\partial_s\widehat{f},\widehat{f})\|_2 & \les \int_1^ts^3\sum_j 2^j \|T_{m^*}(P_ju, u^2, u)\|_2ds \nn \\ 
& \les \int_1^t \sum_j \min(2^{4j}, s^{-1}, s^{-\frac{5}{4}}2^{-3j}  )\|x^{2+}f\|_2\|f\|_X^3ds \lesssim t^{0+}\|f\|_X^4 \label{F91}
\end{align}
 with the same way as in \eqref{s22j}. One can bound $\|F_{5,1}(s\widehat{f},\widehat{f},\partial_s\widehat{f})\|_2$ similarly. To bound $\|F_{5,1}(s\partial_s\widehat{f},\widehat{f},\widehat{f})\|_2$, we need to first estimate $\|P_j(u^2)\|_p$. We compute 
\begin{equation*}
	\begin{split}
		\|P_j(u^2)\|_p&=\|P_j(\sum_kP_ku\cdot u)\|_p\leq \|P_j(P_{<j-5}u\cdot u)\|_p+\|P_j(\sum_{k\geq j-5}P_ku\cdot u)\|_p\\
		&=\|P_j(P_{<j-5}u\cdot P_{j-3\leq n\leq j+3}u)\|_p+\|P_j(\sum_{k\geq j-5}P_ku\cdot u)\|_p\\
		&\leq \|u\|_\infty (\|P_{j-3\leq n\leq j+3}u\|_p+\sum_{k\geq j-5}P_ku)
	\end{split}
\end{equation*}
and therefore,   
\begin{equation}\label{Pj_u2_10}
	\begin{split}
		\|P_j(u^2)\|_{10}&\leq s^{-\frac{5}{4}}\|f\|_X\left( \sum_{n=j-3}^{j+3}s^{-2}2^{-4n}\|f\|_\frac{10}{9}+\sum_{k\geq j-5}s^{-2} 2^{-4k}\|f\|_\frac{10}{9} \right)\\
		&\lesssim s^{-\frac{13}{4}}2^{-4j}\|f\|_X\|f\|_\frac{10}{9}
	\end{split}
\end{equation}
 Hence, we can compute the following. 
 
 \begin{align}\label{F91Pu2}
		\|F_{5,1} (\partial_s\widehat{f},\widehat{f},\widehat{f})\|_2& \les \int_1^ts^3\sum_j 2^j \|T_{m^*}(P_ju, u, u^2)\|_2ds \nn \\
        &\int_1^t s^3\sum_j 2^j\min(\|u\|_8^2\|u^2\|_4, \|P_j(u^2)\|_{10}\|u\|_5^2)ds \nn  \\
		&\lesssim \int_1^t \sum_j\min(s^{-\frac{3}{4}}2^j, s^{-\frac{7}{4}+\varepsilon}2^{-3j})\|f\|_X^4ds\lesssim t^\varepsilon\|f\|_X^4
\end{align}
 
Now,  considering the multiplier in \eqref{m*91} and localizing with \( \varphi_j(\eta)\varphi_k(\eta - \sigma) \)  we obtain:
\begin{multline*}
	\|F_{5,2,1} (s\partial_s\widehat{f},\widehat{f},\widehat{f})\|_2+\|F_{5,2,1} (s\widehat{f},\partial_s\widehat{f},\widehat{f})\|_2+\|F_{5,2,1} (s\widehat{f},\widehat{f},\partial_s\widehat{f})\|_2 \\ 
    \les \int_1^t s^3\sum_k\sum_{j\leq k+10} 2^{5j} \min(s,2^{-4k}) \Big( \|T_{m^*}(u, P_ku, u^2)\|_2+ \|T_{m^*}(u, P_k(u^2), u)\|_2) ds,  
\end{multline*}
and using \eqref{m*92}  and localizing with \( \varphi_j(\eta)\varphi_k(\sigma) \) we obtain 
\begin{multline*}
	\|F_{5,2,2} (s\partial_s\widehat{f},\widehat{f},\widehat{f})\|_2+\|F_{5,2,2} (s\widehat{f},\partial_s\widehat{f},\widehat{f})\|_2+\|F_{5,2,2} (s\widehat{f},\widehat{f},\partial_s\widehat{f})\|_2 \\ 
    \les \int_1^t s^3\sum_k\sum_{j\leq k+10} 2^{5j} \min(s,2^{-4k}) \Big( \|T_{m^*}(u, u, P_k(u^2))\|_2+ \|T_{m^*}(u, u^2, P_ku)\|_2) ds.
\end{multline*}
All of these terms can be bounded  by $t^{0+}\|f\|_X^3$ by \eqref{F91} and \eqref{F91Pu2} and the fact that $2^{5j} \min(s,2^{-4k})\leq 2^j$ in the domain of the sum. 

We finally focus on estimating  $\|\eqref{h2_1_5}\|_2$, where we use the same cut-off functions as in \eqref{defF1} and \eqref{defF2} and study $\|F_{6,2,j}(\widehat{xf},\widehat{f},\widehat{f})\|_2$,  $\|F_{6,2,j}(\widehat{xf},\widehat{f},\widehat{f})\|_2$,  $\|F_{6,1}(\widehat{f},\widehat{xf},\widehat{f})\|_2$, $\|F_{6,2,j}(\widehat{f},\widehat{xf},\widehat{f})\|_2$, $\|F_{6,1}(\widehat{f},\widehat{f},\widehat{xf})\|_2$ and $ \|F_{6,2,j}(\widehat{f},\widehat{f},\widehat{xf})\|_2$. 
Using the same Coifman-Meyer multipliers as in \eqref{m*9} with $Q_5$ exchanged with $Q_6$, we obtain 
\begin{align*}
 \|F_{6,1}(\widehat{xf},\widehat{f},\widehat{f})\|_2 \les \int_1^t s^2\sum_j 2^{2j} \|T_{m^*}(P_j(e^{-is\Delta^2}xf), u,u)\|_2 ds  \les t^{0+}\|f\|_X^3 
\end{align*}
by \eqref{F61}. On the other hand, we have 
\begin{align} \label{s22jxf}
 \|F_{6,1}(\widehat{f},\widehat{xf},\widehat{f})\|_2+ &\|F_{6,1}(\widehat{f},\widehat{f},\widehat{xf})\|_2 \les \int_1^t s^2\sum_j 2^{2j}  \|T_{m^*}(P_j(u),e^{-is\Delta^2}xf,u)\|_2 ds \nn \\ 
 & \les \int_1^t s^2\sum_j 2^{2j}\|P_ju\|_{10}\|u\|_{10}\|e^{-is\Delta^2}(xf)\|_\frac{10}{3}ds \nn \\
		&\lesssim \int_1^t \sum_j\min(s^{-\frac{3}{2}+\varepsilon}2^{-2j}, s^{\frac{1}{2}+\varepsilon}2^{6k})\|f\|_X^3ds \lesssim t^\varepsilon \|f\|_X^3
\end{align}
where we used $\|P_j u\|_{10} \les 2^{4k}\|f\|_{\frac{10}{9}}\les2^{4k} s^{0+}\|f\|_X$ in the second inequality together with

\begin{align}\label{Decay_xf_10_3}
		\|e^{-is\Delta^2}(xf)\|_\frac{10}{3}&\leq\sum_{2^j\leq t^{-\frac{1}{4}}}\|P_je^{-is\Delta^2}(xf)\|_\frac{10}{3}+\sum_{2^j> t^{-\frac{1}{4}}}\|P_je^{-is\Delta^2}(xf)\|_\frac{10}{3} \nn \\
		&\lesssim \sum_{2^j\leq t^{-\frac{1}{4}}} 2^{2j} \|xf\|_\frac{10}{7}+\sum_{2^j> t^{-\frac{1}{4}}}s^{-1}2^{-2j}\|xf\|_\frac{10}{7} \lesssim s^{-\frac{1}{2}+\varepsilon}\|f\|_X
\end{align}
Moreover, using the multiplier $\eqref{m*91}$ with $Q_5$ exchanged with $Q_6$, and localizing in $\varphi_l(\xi-\eta)\varphi_j(\eta) \varphi_k(\eta-\sigma)$, we reduce estimating $\|F_{6,2,1}(\widehat{f},\widehat{xf},\widehat{f})$, $\|F_{6,2,1}(\widehat{xf},\widehat{f},\widehat{f})$ and $\|F_{6,2,1}(\widehat{f},\widehat{f},\widehat{xf})$ to estimating the $L^2$ norm of the following terms
\begin{align*}
&\int_1^t s^2\sum_k\sum_{j\leq k+10} \sum_{l\leq j+10}2^{6j}\min(s,2^{-4k})\|T_{m^*}(P_l(e^{-is\Delta^2}xf),P_k u,u)\|_2ds  \\ 
&\int_1^t s^2\sum_k\sum_{j\leq k+10} \sum_{l\leq j+10}2^{6j}\min(s,2^{-4k})\|T_{m^*}(P_lu,P_k u,e^{-is\Delta^2}xf)\|_2ds   
\end{align*}

The first term has been estimated by $t^{0+} \|f\|_X^3$ in \eqref{F621} and the second term reduces to $\eqref{s22jxf}$ using Young's inequality and the fact that $2^{6j}\min(s,2^{-4k})\leq 2^{2k+}$ in the domain of the sum. 

The terms $\|F_{6,2,2}(\widehat{f},\widehat{xf},\widehat{f})$, $\|F_{6,2,2}(\widehat{xf},\widehat{f},\widehat{f})$ and $\|F_{6,2,2}(\widehat{f},\widehat{f},\widehat{xf})$ can be estimated similarly  using the multiplier $\eqref{m*92}$ with $Q_5$ exchanged with $Q_6$, and localizing in $\varphi_l(\xi-\eta)\varphi_j(\eta) \varphi_k(\sigma)$. 

\end{proof}

\begin{lemma}
	Assuming \eqref{Boot_est_g} and \eqref{Boot_est_h}, we have
	\begin{equation}
		\|x^3h_3\|_2\lesssim  t^{\frac{1}{24}}(\|f\|_X^2+\|f\|_X^3+\|f\|_X^4)
	\end{equation} \label{Decay_x3h3} 
\end{lemma}
\begin{proof} We start computing 
\begin{align}
\label{h3_1}\partial_\xi^3(\widehat{h_3}(\xi,t))&=\int_1^t s^2\int \frac{Q_9}{A}e^{is\varphi}\widehat{f}(\xi-\eta,s)\widehat{f}(\eta,s)d\eta ds\\
		\label{h3_2}&+\int_1^t s\int \frac{q_6}{A}e^{is\varphi}\widehat{xf}(\xi-\eta,s)\widehat{f}(\eta,s)d\eta ds\\
		\label{h3_3}&+\int_1^t s\int \frac{q_5}{A}e^{is\varphi}\widehat{f}(\xi-\eta,s)\widehat{f}(\eta,s)d\eta ds\\
		\label{h3_4}&+\int_1^t \int \frac{q_1}{A}e^{is\varphi}\widehat{f}(\xi-\eta,s)\widehat{f}(\eta,s)d\eta ds\\
		\label{h3_5}&+\int_1^t \int e^{is\varphi} \Big( \frac{q_2}{A}\widehat{xf}(\xi-\eta,s)+\frac{q_3}{A} \widehat{x^2f}(\xi-\eta,s)\Big)\widehat{f}(\eta,s)d\eta ds\\
		\label{h3_7}&+\int_1^t s^{-1}\int \frac{1}{A}e^{is\varphi}\widehat{x^3f}(\xi-\eta,s)\widehat{f}(\eta,s)d\eta ds
	\end{align}
First, note that using $\eqref{fracwuse}$, we can estimate $\| \eqref{h3_4}\|_2+ \|, \eqref{h3_5}\|_2\les t^{0+}\|f\|_X^2$. To estimate $\|\eqref{h3_3}\|_2$, we use \eqref{key2}. After this application, the resulting terms either take the form of \eqref{h3_2}, \eqref{h3_4}, or one of the following:
\begin{align} \label{sm1}
 t\int\frac{q_1}{A}e^{it\varphi}\widehat{f}(\xi-\eta,t)\widehat{f}(\eta,t)d\eta, \,\,\ \int_1^t s\int \frac{q_1}{A}e^{is\varphi}\widehat{f}(\xi-\eta,s)\partial_s\widehat{f}(\eta,s)d\eta ds. 
\end{align}
Using \eqref{fracwuse}, the \(L^2\) norm of both terms can be estimated as  $t^{0+}\|f\|_X^2+\|f\|_X^3$.

We also use \eqref{key2} to estimate $\|\eqref{h3_1}\|_2$. After we apply, the resulting terms are either in the form of \eqref{h3_3}, \eqref{h3_2} or one of the following 
\begin{align*}
t^2\int \frac{q_9}{A^2}  e^{it\varphi}\widehat{f}(\xi-\eta,t)\widehat{f}(\eta,t)d\eta, \,\,\ \int_1^t s^2\int \frac{q_5}{A^2}  e^{is\varphi}\widehat{f}(\xi-\eta,s)\partial_s\widehat{f}(\eta,s)d\eta ds
\end{align*}
The second term appeared in \eqref{h2_3} and its $L^2$ norm is estimated by $t^{0+}\|f\|_X^3$. To estimate the first term, we use the decomposition \eqref{Decom_Bi} of $f$. We observe that the resulted terms can be bounded by the $L^2$ norm of \eqref{h1_1_5} and \eqref{h1_2_1}. As a result, we obtain $\|\eqref{h3_1}+\eqref{h3_3}\|_2\lesssim \|f\|_X^2+t^{\frac{1}{48}+}\|f\|_X^3$ provided that we establish the required bound for $\|\eqref{h3_2}\|_2$.

We again use \eqref{key2} to decompose \eqref{h3_2}. After the decomposition, we obtain terms in the form of \eqref{h3_4}, \eqref{h3_5}, \eqref{x2h_2_f}, or one of the following:
\begin{align}
\int_1^t \int s \frac{q_2}{A} e^{is\varphi} \partial_s \partial_{\eta} \widehat{f}(\xi-\eta,s) \widehat{f}(\eta,s) \, d\eta ds, \quad
\int_1^t \int  \frac{Q_7}{A^2} e^{is\varphi} \widehat{xf}(\xi-\eta,s) \widehat{xf}(\eta,s) \, d\eta ds \label{xfxf}
\end{align}
Applying integration by parts in \(\eta\) to the first term, we reduce it to terms appearing in either \eqref{h2_3}, \eqref{h2_4}, or \eqref{h2_6}. To estimate the second term we apply \eqref{Decom_Bi} and reduce it to estimating:
	\begin{align}
		\label{h3_2_1}&\int_1^t \int \frac{Q_7}{A^2}e^{is\varphi}\widehat{xh}(\xi-\eta,s)\widehat{xh}(\eta,s)d\eta ds\\
		\label{h3_2_2}&\int_1^t \int \frac{Q_7}{A^2}e^{is\varphi}\widehat{xf}(\xi-\eta,s)\widehat{xg}(\eta,s)d\eta ds\\
		\label{h3_2_3}&\int_1^t \int \frac{Q_7}{A^2}e^{is\varphi}\widehat{xg}(\xi-\eta,s)\widehat{xh}(\eta,s)d\eta ds
	\end{align}
	We estimate \eqref{h3_2_1} using Lemma~\ref{lem_xf_4}:
	\begin{equation*}
		\|\eqref{h3_2_1}\|_2\lesssim \int_1^t s^{\frac{1}{4}}\|e^{-is\Delta^2}(xh)\|_4^2ds\lesssim \log(t)\|f\|_X^2
	\end{equation*}
We use \eqref{eq:xg}, the expression of $\widehat{xg}$ and write 
	\begin{align}
		\eqref{h3_2_2}&\label{102_1}=\int_1^t s\iint \frac{q_3(\xi,\eta,s)}{A}\frac{Q_3(\eta,\sigma)}{C}e^{is\psi}\widehat{xf}(\xi-\eta,s)\widehat{f}(\eta-\sigma,s)\widehat{f}(\sigma,s)d\sigma d\eta ds\\
		\label{102_2}&+\int_1^t\iint \frac{q_7}{A^2}\frac{1}{C}e^{is\psi}\widehat{xf}(\xi-\eta,s)\widehat{xf}(\eta-\sigma,s)\widehat{f}(\sigma,s)d\sigma d\eta ds
	\end{align}

The first term can be estimated using \eqref{Decay_xf} for $\alpha =1$ and \eqref{6.44est}. 
We apply Lemma \ref{lem_xf_4} and obtain
	\begin{equation*}
		\|\eqref{102_2}\|_2\lesssim\int_1^t s^\frac{5}{4}\|e^{-is\Delta^2}(xf)\|_4^2\|u\|_\infty ds\lesssim \|f\|_X^3
	\end{equation*}
	Notice that \eqref{h3_2_3} can be estimated as \eqref{h3_2_2}, so we finished estimating $\|\eqref{h3_2}\|_2. $
    
    We are left with the final term  \eqref{h3_7}. We first decompose it using \eqref{Decom_Bi} and study the following terms 
	\begin{align}
		\label{h3_7_1}&\int_1^t s^{-1}\int \frac{1}{A}e^{is\varphi}\widehat{x^3h}(\xi-\eta,s)\widehat{h}(\eta,s)d\eta ds, \,\,\, \int_1^t s^{-1}\int \frac{1}{A}e^{is\varphi}\widehat{x^3f}(\xi-\eta,s)\widehat{g}(\eta,s)d\eta ds\\
		\label{h3_7_3}& \int_1^t s^{-1}\int \frac{1}{A}e^{is\varphi}\widehat{x^3h}(\xi-\eta,s)\widehat{g}(\eta,s)d\eta ds, \,\ \int_1^t s^{-1}\int \frac{1}{A}e^{is\varphi}\widehat{x^3g}(\xi-\eta,s)\widehat{f}(\eta,s)d\eta ds 
	\end{align}
	\eqref{h3_7_1} has previously appeared in \eqref{h1_1_1} and \eqref{h1_2_5}, where we applied \eqref{eqn_g} to express the integral in a trilinear form. Moreover, using \eqref{eqn_g}, estimating  the \( L^2 \) norm of the first term in \eqref{h3_7_3} reduces to estimating the $L^2$ norm of \eqref{h1_2_5}, since the time growth of the weighted norm estimates of \( h \) is slower than that of \( f \).  To estimate the \( L^2 \) norm of  the second term in \eqref{h3_7_3}, we first reparametrize it:
	\begin{equation}
		\eqref{h3_7_3}=\int_1^t \int q_0(\xi,\eta,s) e^{is\varphi}\widehat{f}(\xi-\eta,s)\widehat{x^3g}(\eta,s)d\eta ds  \label{q0x3g}
	\end{equation}
Here,  $q_0(\xi,\eta,s)$ refers to $(s A^{\prime})^{-1} $  where $A^\prime =\frac{1}{s}+iZ^\prime$ with $Z^\prime\gtrsim\max(|\xi|^4+|\xi-\eta|^4, |\xi|^4+|\eta|^4) $. We next use the expression for $ \widehat{x^3g}$ computed in \eqref{x3g_1} to represent 
	\begin{align}
		\eqref{h3_7_3}&\label{108_1}=\int_1^t s^3\iint q_0(\xi,\eta,s) \frac{Q_9(\eta,\sigma)}{C}e^{is\psi} \widehat{f}(\xi-\eta,s)\widehat{f}(\eta-\sigma,s)\widehat{f}(\sigma,s)d\sigma d\eta ds\\
		\label{108_2}&+\int_1^t s^2\iint q_0(\xi,\eta,s) \frac{q_6(\eta,\sigma,s)}{C}e^{is\psi} \widehat{f}(\xi-\eta,s)\widehat{xf}(\eta-\sigma,s)\widehat{f}(\sigma,s)d\sigma d\eta ds\\
		\label{108_3}&+\int_1^t s^2 \iint q_0(\xi,\eta,s)\frac{q_5(\eta,\sigma,s)}{C}e^{is\psi} \widehat{f}(\xi-\eta,s)\widehat{f}(\eta-\sigma,s)\widehat{f}(\sigma,s)d\sigma d\eta ds\\
		\label{108_5}&+\int_1^t \iint sq_0(\xi,\eta,s)\frac{q_3(\eta,\sigma,s)}{\: C} e^{is\psi} \widehat{x^2f}(\eta-\sigma,s) \widehat{f}(\xi-\eta,s)\widehat{f}(\sigma,s)d\sigma d\eta ds\\
        \label{108_8}&+\int_1^t \iint s q_0(\xi,\eta,s) \frac{q_2 (\eta,\sigma,s)}{\: C}e^{is\psi}  \widehat{xf}(\eta-\sigma,s) \widehat{f}(\xi-\eta,s)\widehat{f}(\sigma,s)d\sigma d\eta ds\\
		\label{108_6}&+\int_1^t s\iint q_0(\xi,\eta,s)\frac{q_1(\sigma,\eta,s)}{C}e^{is\psi} \widehat{f}(\xi-\eta,s)\widehat{f}(\eta-\sigma,s)\widehat{f}(\sigma,s)d\sigma d\eta ds\\
		\label{108_7}&+\int_1^t  \iint  \frac{q_0(\xi,\eta,s)}{C}e^{is\psi} \widehat{f}(\xi-\eta,s)\widehat{x^3f}(\eta-\sigma,s)\widehat{f}(\sigma,s)d\sigma d\eta ds
	\end{align}
We note that the terms in \eqref{108_5}  and \eqref{108_8} appeared in \eqref{h1_2_4} and  \eqref{h1_2_2}.  The fastest growth among these terms is $t^\frac{1}{24}\|f\|_X^3$.Moreover, we can bound $\|\eqref{108_6}\|_2$ using
$$
\|\eqref{108_6}\|_2 \les s \|u\|_{20}^2\|\nabla_s^{-3}u\|_\frac{5}{2} \les s^{-1} \|f\|_\frac{5}{4}\|f\|_X^2. 
$$
We write \eqref{108_3} as the sum of the following terms. 

\begin{multline}\label{G51}
	G_{5,i}(\widehat{f},\widehat{f},\widehat{f})\\ =\int_1^t s^2\iint \frac{ q_0(\xi,\eta,s) Q_5(\sigma,\eta)}{C} \Psi_i(\eta,\sigma)e^{is\psi} \widehat{f}(\xi-\eta,s)\widehat{f}(\eta-\sigma,s)\widehat{f}(\sigma,s)d\sigma d\eta ds
\end{multline}
for $i=1,2$. 
To estimate $\|G_{5,i}(\widehat{f},\widehat{f},\widehat{f})\|2$, we use 
\begin{align}\label{m*G1}
m^*(\xi,\eta,\sigma,s)= m^*_1(\eta,\sigma) m^*_3(\xi,\eta, \sigma,s)  = q_0(\xi,\eta, s) \:\frac{Q_5(\sigma,\eta)\Psi_1(\eta,\sigma) \varphi_j(\eta-\sigma) }{C 2^{5j} \min(s,2^{-4j})}\:   
\end{align}
and observe that $m^*(\xi,\eta,\sigma,s)$ defines a Coifman-Meyer multiplier with flag singularity with norm independent on $s$. Hence, we have 
\begin{align*}
\|G_{5,1}(\widehat{f},\widehat{f},\widehat{f})\|_2 \les \int_1^t s^2 \sum_{j} 2^{j} \|T_{m^*}(u, P_ju,u)\|_2 \les t^{0+} \|f\|^3_X. 
\end{align*}
In the last inequality we used \eqref{s22j}. Similarly, using the multiplier \eqref{m*G1} with $\varphi_j(\eta-\sigma)$ exchanged with 
$\varphi_j(\sigma)$, we bound 
\begin{align*}
\|G_{5,2}(\widehat{f},\widehat{f},\widehat{f})\|_2 \les \int_1^t s^2 \sum_{j} 2^{j} \|T_{m^*}(u, u,P_ju)\|_2 \les t^{0+} \|f\|^3_X. 
\end{align*}

With a similar decomposition, we reduce estimating $\|\eqref{108_2}\|_2$ to estimating $\|G_{6,j}(\widehat{f}, \widehat{xf},\widehat{f})\|_2$ for $j=1,2$. Using the multiplier $m^*$ in \eqref{m*G1} with $Q_5$ exchanged with $Q_6$, we obtain 
\begin{align*}
\|G_{6,1}(\widehat{f}, \widehat{xf},\widehat{f})\|_2 \les \int_1^t s^2 \sum_{j}2^{2j}\|T_{m^*}(u, P_j(e^{-is\Delta^2}(xf)), u)\|_2 \les t^{0+} \|f\|_X^3
\end{align*}
where we used  \eqref{F61} in the last inequality. 

To estimate $\|G_{6,2}(\widehat{f}, \widehat{xf},\widehat{f})\|_2$, we use the multiplier \eqref{m*G1} with $Q_5$ exchanged with $Q_6$, and $\varphi_j(\eta-\sigma)$  with 
$\varphi_j(\sigma)$ to see 
\begin{align} 
\|G_{6,2}(\widehat{f}, \widehat{xf},\widehat{f})\|_2 &\les \int_1^t s^2 \sum_j  2^{2j}\|T_{m^*}(u,e^{-is\Delta^2}(xf), P_ju)\|_2 
\end{align}
which can be bounded by $t^{\varepsilon}\|f\|_X^3$ using \eqref{s22jxf}. 

To estimate \eqref{108_7}, we decompose it into two parts: 
\begin{equation*}
	G_{0,i}(\widehat{f},\widehat{x^3f}, \widehat{f})=\int_1^t \iint  q_0(\xi,\eta,s)\frac{1}{C}\Psi_i(\eta,\sigma)  e^{is\psi}\widehat{f}(\xi-\eta,s)\widehat{x^3f}(\eta-\sigma,s)\widehat{f}(\sigma,s)d\sigma ds
\end{equation*} 
For $G_{0,1}(\widehat{f},\widehat{x^3f}, \widehat{f})$, first observe that 
\begin{equation*}
	 m^*(\xi,\eta,\sigma,s)=q_0(\xi,\eta,s)\cdot\frac{1}{C\min(s, 2^{-4j})}\Psi_1(\eta,\sigma)\varphi_j(\eta-\sigma)
\end{equation*}
is a Coifman-Meyer multiplier with flag singularity. Using Remark \ref{middle_j}, we consider when $s\leq 2^j\leq s^{-1}$:
\begin{equation*}
	\begin{split}
		\|G_{0,1}(\widehat{f},\widehat{x^3f}, \widehat{f})\|_{L^2}&\leq\sum_j\int_1^t \min(s, 2^{-4j})\|T_{m^*}(u, P_j(e^{-is\Delta^2}(x^3f)), u)\|_{L^2}ds\\
		&\lesssim\sum_j\int_1^t s^{\frac{63}{64}}2^{-\frac{j}{16}}\|u\|_{160}^2\|P_j(e^{-is\Delta^2}(x^3f))\|_{\frac{80}{39}}ds\\
		&\lesssim \log(t)\|f\|_X^2 \int_1^t s^{\frac{63}{64}-\frac{79}{32}}\|x^3f\|_{L^2}ds\lesssim t^{\frac{1}{26}}\|f\|_X^3
	\end{split}
\end{equation*}
and using 
\begin{equation*}
	 m^*(\xi,\eta,\sigma,s)=q_0(\xi,\eta,s)\cdot\frac{1}{C\min(s, 2^{-4j})}\Psi_1(\eta,\sigma)\varphi_j(\sigma)
\end{equation*}
we have 
\begin{align*}
\|G_{0,2}(\widehat{f},\widehat{x^3f}, \widehat{f})\|_{L^2}&\leq\sum_j\sum_{k\leq j+10} \int_1^t \min(s, 2^{-4j})\|T_{m^*}(u, P_k(e^{-is\Delta^2}(x^3f)),P_j u)\|_{L^2}ds\\
		&\sum_j \sum_{k\leq j+10} \int_1^t s\|u\|_{160}\|P_k(e^{-is\Delta^2} (x^3f))\|_{\frac{80}{39}} \|P_ju\|_{160}ds\\
		&\lesssim \sum_j \int_1^t 2^{-\frac{39j}{16}+}  s^{-\frac{173}{64}} \|f\|_X^2 \|x^3f\|_2 \les  \|f\|_X^3 
\end{align*}

To estimate $\|\eqref{108_1}\|_2$, we apply \eqref{key3} to decompose it to the following terms:
\begin{align}
	\label{h3_1_1_n}&\int_1^t s^2 \iint 
 q_0(\xi,\eta,\sigma,s)\frac{Q_9(\sigma,\eta)}{CB}e^{is\psi} \widehat{f}(\xi-\eta,s)\widehat{f}(\eta-\sigma,s)\widehat{f}(\sigma,s)d\sigma d\eta ds\\
    \label{h3_1_3_n}&\int_1^t s^2 \iint (q_0(\xi,\eta,s)+q_0(\eta,\sigma,s)) \frac{Q_9(\sigma,\eta)}{CB}e^{is\psi} \widehat{f}(\xi-\eta,s)\widehat{f}(\eta-\sigma,s)\widehat{f}(\sigma,s)d\sigma d\eta ds\\
	\label{h3_1_2_n}&t^2\iint q_0(\xi,\eta,t)\frac{Q_9(\sigma,\eta)}{CB}e^{it\psi} \widehat{f}(\xi-\eta,t)\widehat{f}(\eta-\sigma,t)\widehat{f}(\sigma,t)d\sigma d\eta \\
	\label{h3_1_4_n}&\int_1^t s^3\iint q_0(\xi,\eta,s)\frac{Q_9(\sigma,\eta)}{CB}e^{is\psi} \partial_{s}\{\widehat{f}(\xi-\eta,s) \widehat{f}(\eta-\sigma,s)\widehat{f}(\sigma,s)\}d\sigma d\eta ds\\
	\label{h3_1_6_n}&\int_1^t s^2\iint q_0(\xi,\eta,s)\frac{Q_{9}(\sigma,\eta)q_1(\xi,\eta,\sigma)}{CB}e^{is\psi} \widehat{f}(\xi-\eta,s)\widehat{xf}(\eta-\sigma,s)\widehat{f}(\sigma,s)d\sigma d\eta ds
\end{align}

We note that \eqref{h3_1_1_n} and \eqref{h3_1_3_n} can be reduced to studying \( G_{5,j}(\widehat{f},\widehat{f},\widehat{f}) \), while \eqref{h3_1_6_n} reduces to studying \( G_{6,j}(\widehat{f},\widehat{xf},\widehat{f}) \), where \( G_{i,j} \) is defined in \eqref{G51}. Both of these terms are estimated by \( t^{0+} \|f\|_X^3 \). Moreover, it follows that \( \|\eqref{h3_1_2_n}\|_2 \) satisfies the same bound as \( \|G_{6,j}(\widehat{f},\widehat{f},\widehat{f})\|_2 \).

Therefore, we continue analyzing  \eqref{h3_1_4_n}. We decompose this term into 
\begin{align*}
&G_{9,i}(\widehat{f},\widehat{f},\widehat{f})=  \int_1^t s^3 \iint \frac{q_0(\xi,\eta) Q_9(\eta,\sigma)}{CB} \phi_j(\xi,\eta,\sigma)  e^{is \psi} \widehat{f}(\xi-\eta,s)\widehat{f}(\eta-\sigma,s)\widehat{f}(\sigma,s)d\sigma d\eta ds \\
&G_{9,1,\ell}(\widehat{f},\widehat{f},\widehat{f})= \int_1^t s^3 \iint \frac{q_0(\xi,\eta) Q_9(\eta,\sigma)}{CB} \phi_j(\xi,\eta,\sigma)  \Psi_{\ell} (\eta,\sigma) e^{is \psi} \widehat{f}(\xi-\eta,s)\widehat{f}(\eta-\sigma,s)\widehat{f}(\sigma,s)d\sigma d\eta ds
\end{align*}
for $i=2,3$ and $\ell=1,2$. We first estimate $ \|G_{9,2}(\partial_s\widehat{f},\widehat{f},\widehat{f})\|_2 $, $ \|G_{9,2}(\widehat{f},\partial_s\widehat{f},\widehat{f})\|_2 $, $ \|G_{9,2}(\widehat{f},\widehat{f},\partial_s\widehat{f})\|_2 $
To use Theorem~\ref{th:flagcf}, we let 
\begin{align} \label{m*G92}
m^*(\xi,\eta,\sigma,s)= q_0(\xi,\eta) \: q_0(\sigma,\eta) \: \frac{Q_5(\eta,\sigma) \psi_j(\eta-\sigma)\phi_2(\xi,\eta,\sigma)}{B \: \min(s,2^{-4j})}
\end{align}
as Coifman-Meyer multiplier with flag singularity operator norm independent of $s$ and estimate 
\begin{align*}
\|G_{9,2}(\partial_s\widehat{f},\widehat{f},\widehat{f})\|_2+& \|G_{9,2}(\widehat{f},\partial_s\widehat{f},\widehat{f})\|_2+\|G_{9,2}(\widehat{f},\widehat{f},\partial_s\widehat{f})\|_2 \\
&\les \int_1^t s^3 \sum_{j} 2^{j} \Big( \|T_{m*}(u,P_ju,u^2)\|_2+ \|T_{m*}(u,u,P_j(u^2)\|_2 \Big) \les t^{0+} \|f\|_X^3 
\end{align*}
using \eqref{F91} and \eqref{F91Pu2}. The same approach can also estimate $ \|G_{9,3}(\widehat{f},\partial_s\widehat{f},\widehat{f})\|_2 $, $ \|G_{9,3}(\widehat{f},\widehat{f},\widehat{f})\|_2 $, $ \|G_{9,3}(\widehat{f},\widehat{f},\partial_s\widehat{f})\|_2 $
using $\psi_j(\sigma)$ instead of $\psi_j(\eta-\sigma)$, and $\phi_3(\xi,\eta,\sigma)$ instead of $\phi_2(\xi,\eta,\sigma)$ in \eqref{m*G92}.

We estimate $ \|G_{9,1,1}(\widehat{f},\partial_s\widehat{f},\widehat{f})\|_2 $, $ \|G_{9,1,1}(\widehat{f},\widehat{f},\widehat{f})\|_2 $, $ \|G_{9,1,1}(\widehat{f},\widehat{f},\partial_s\widehat{f})\|_2 $
using the multiplier 
\begin{align*}
m^*(\xi,\eta,\sigma,s)= q_0(\xi,\eta) \: \frac{\psi_j(\eta-\sigma)Q_9(\eta,\sigma)}{C\: 2^{9j} \: \min(s,2^{-4j})} \: \frac{\psi_k(\xi-\eta)\phi_1(\xi,\eta,\sigma)}{B \: \min(s,2^{-4k})}
\end{align*}
which then allows us the following estimate 
\begin{align}
\|&G_{9,1,1}(\partial_s\widehat{f},\widehat{f},\widehat{f})\|_2+ \|G_{9,1,1}(\widehat{f},\partial_s\widehat{f},\widehat{f})\|_2+\|G_{9,1,1}(\widehat{f},\widehat{f},\partial_s\widehat{f})\|_2 \nn \\
&\les \int_1^t s^3  \sum_{k}\sum_{j\leq k+10} 2^j\Big( \|T_{m*}(P_k(u^2),P_ju,u^2)\|_2+ \|T_{m*}(P_ku,P_j(u^2),u)\|_2 + \|T_{m*}(P_ku,P_ju,u^2)\|_2  \Big) \nn  
\end{align}
We have 
\begin{align}\label{G911}
 s^3 \sum_{k}\sum_{j\leq k+10} 2^j & \Big(\|T_{m*}(P_ku,P_j(u^2),u)\|_2 +\|T_{m*}(P_ku,P_ju,u^2)\|_2 \Big) \nn
 \\
 & \les s^3 \sum_{k} 2^{k+} \min(\|u\|_{10}\|u^2\|_{10}\|P_ku\|_{\frac{10}{3}}, \|P_ku\|_{10}\|u^2\|_{5}\|u\|_5) \nn \\
 &\les \sum_j \min(2^{4j}, s^{-1}, s^{-\frac{5}{4}}2^{-3j}  )\|x^{2+}f\|_2\|f\|_X^2
\end{align}
similar to \eqref{s22j}. Moreover, 
\begin{align}\label{G9111}
s^3 \sum_{k}\sum_{j\leq k+10} 2^j \|T_{m*}(P_k(u^2),P_ju,u)\|_2 &\les 
		s^3\sum_k\sum_{j\leq k+10 }2^j \|P_ju\|_{10}\|u\|_{\frac{10}{3}}\|P_k(u^2)\|_{10} \nn \\ 
        &\les\sum_k \min(s^\frac{1}{4}2^{5k}, s^{-\frac{27}{16}}2^{-\frac{23}{8}k}  )\|f\|_\frac{10}{9}\|f\|_X^3  
\end{align}
where we used \eqref{Pj_u2_10} to estimate $\|P_k(u^2)\|_{10}$. Therefore, we obtain 
$$
\|G_{9,1,1}(\partial_s\widehat{f},\widehat{f},\widehat{f})\|_2+ \|G_{9,1,1}(\widehat{f},\partial_s\widehat{f},\widehat{f})\|_2+\|G_{9,1,1}(\widehat{f},\widehat{f},\partial_s\widehat{f})\|_2 \les t^{\frac{1}{31}}\|f\|_X^4.$$

Finally, we can use the multiplier 
\begin{align*}
m^*(\xi,\eta,\sigma,s)= q_0(\xi,\eta) \: \frac{\psi_j(\sigma)Q_9(\eta,\sigma)}{C\: 2^{9j} \: \min(s,2^{-4j})} \: \frac{\psi_k(\xi-\eta)\phi_1(\xi,\eta,\sigma)}{B \: \min(s,2^{-4k})}
\end{align*}
to estimate 
    \begin{align}
\|&G_{9,1,2}(\partial_s\widehat{f},\widehat{f},\widehat{f})\|_2+ \|G_{9,1,2}(\widehat{f},\partial_s\widehat{f},\widehat{f})\|_2+\|G_{9,1,2}(\widehat{f},\widehat{f},\partial_s\widehat{f})\|_2 \nn \\
&\les \int_1^t s^3  \sum_{k}\sum_{j\leq k+10} 2^j\Big( \|T_{m*}(P_k(u^2),u,P_ju)\|_2+ \|T_{m*}(P_ku,(u^2),P_ju\|_2 + \|T_{m*}(P_ku,u,P_j(u^2)\|_2  \Big) \nn  
\end{align}
The first term can be estimated using \eqref{G9111}
and the last two can be estimated  using  \eqref{G9111} by $t^{\frac{1}{31}}\|f\|_X^4$. This establishes the required estimate for $\|x^3h_3\|_2$. 
\end{proof}
Prove Lemma~\ref{Lemma Decay_x2h_4} for further calculations. 
\begin{lemma}\label{Lemma Decay_x2h_4}
Assuming \eqref{Boot_est_g} and \eqref{Boot_est_h}, we have
\begin{align}
    &\label{Decay_x2f_4} \|e^{-is\Delta^2}(x^2f)\|_4\lesssim s^{-\frac{19}{48}+\varepsilon}\|f\|_X^2 \\ 
&\label{Decay_x2h_4}\|e^{-is\Delta^2}(x^2h)\|_4\lesssim s^{-\frac{113}{140}+\varepsilon}\|f\|_X^2
\end{align}
\end{lemma}
\begin{proof}
	By Plancherel, we have:
\begin{align}
    \|e^{-is\Delta^2}&(x^2h)\|_4\lesssim \left\| \sum_j\int e^{-is|\xi-\eta|^4}\Delta \widehat{h}(\xi-\eta)e^{is|\eta|^4}  \Delta\widehat{\overline{h}}(\eta)\varphi_j(\xi-\eta)\varphi_{\leq j}(\eta)d\eta\right\|_{L^2_\xi} \nn \\
			&\leq \sum_{2^j<s^{-\beta}}\left\|P_j(e^{-is\Delta^2}(x^2h))P_{\leq j}(e^{-is\Delta^2}(x^2h)) \right\|_2+\sum_{2^j\geq s^{-\beta}}\left\|P_j(e^{-is\Delta^2}(x^2h))P_{\leq j}(e^{-is\Delta^2}(x^2h)) \right\|_2 \nn \\
			&=:I(x^2h)+II(x^2h) \label{I(x^2h)}
\end{align}
	
We use \eqref{bern} to estimate estimate $I(x^2h)$ as 
	\begin{align*}
			I(x^2h)\leq \sum_{2^j<s^{-\beta}} \|P_j(e^{-is\Delta^2}(x^2h))\|_2\|P_{\leq j}(e^{-is\Delta^2}(x^2h)) \|_\infty \les \sum_{2^j<s^{-\beta}} 2^{\frac{17}{4}j}\|x^2h\|_\frac{40}{27}^2\lesssim s^{\frac{1}{12}+\varepsilon-\frac{17}{4}\beta}\|f\|_X^2. 
	\end{align*}
Moreover, using Lemma \ref{lem:ben} we have 
	\begin{align*}
			II\leq & \sum_{2^j\geq s^{-\beta}} \|P_j(e^{-is\Delta^2}(x^2h))\|_\frac{40}{13} \|P_{\leq j}(e^{-is\Delta^2}(x^2h)) \|_\frac{40}{7} \\
            & \lesssim \sum_{2^j\geq s^{-\beta}} s^{-\frac{7}{8}}2^{-\frac{1}{8}j}\|x^2h\|_\frac{40}{27}\|x^2h\|_2\lesssim s^{-\frac{5}{6}+\varepsilon+\frac{1}{8}\beta}\|f\|_X^2
	\end{align*}
Letting $\beta=\frac{22}{105}$ we obtain the inequality \eqref{Decay_x2h_4}. 

To obtain \eqref{Decay_x2f_4} we use $f$ in \eqref{I(x^2h)} instead of $h$ and using \eqref{bern} we obtain
\begin{align*}
    I(x^2f)&\leq \sum_{2^j<s^{-\beta}} \|P_j(e^{-is\Delta^2}(x^2f))\|_2\|P_{\leq j}(e^{-is\Delta^2}(x^2f)) \|_\infty \lesssim \sum_{2^j<s^{-\beta}} 2^{\frac{5}{2}j}\|x^2f\|_2^2\lesssim s^{\varepsilon-\frac{5}{2}\beta}
\end{align*}
Similarly by Lemma~\ref{lem:ben}, we have 
	\begin{equation*}
		\begin{split}
			II(x^2f)&\leq \sum_{2^j\geq s^{-\beta}} \|P_j(e^{-is\Delta^2}(x^2f))\|_\frac{16}{5} \|P_{\leq j}(e^{-is\Delta^2}(x^2f)) \|_\frac{16}{3}\\
			&\lesssim \sum_{2^j\geq s^{-\beta}} s^{-\frac{15}{16}+\frac{63}{128}+\varepsilon}2^{-\frac{5}{16}j}\|f\|_X^2\lesssim s^{-\frac{57}{128}+\frac{5}{16}\beta+\varepsilon}
		\end{split}
	\end{equation*}
Letting  $\beta=\frac{19}{120}$ we finish the proof.

\end{proof}

\begin{lemma}
	Assuming \eqref{Boot_est_g} and \eqref{Boot_est_h}, we have
	\begin{equation}
		\|x^3h_4\|_2\lesssim t^{\frac{1}{24}}\|f\|_X^2 \label{Decay_x3h4} 
	\end{equation}
\end{lemma}

\begin{proof}
We start computing $\partial_\xi^3(\widehat{h_4})$:
	\begin{align}
		\partial_\xi^3\widehat{h_4}(\xi,t)&=\label{h4_1}\int_1^t s^3\int \frac{Q_1\cdot Q_9 \varphi_\eta}{A}e^{is\varphi}\widehat{f}(\xi-\eta,s)\widehat{f}(\eta,s)d\eta ds\\
		\label{h4_2}&+\int_1^t s^2\int \frac{Q_7\varphi_\eta}{A}e^{is\varphi}\widehat{xf}(\xi-\eta,s)\widehat{f}(\eta,s)d\eta ds\\
		\label{h4_4}&+\int_1^t s\int \frac{Q_4\varphi_\eta}{A}e^{is\varphi}\widehat{x^2f}(\xi-\eta,s)\widehat{f}(\eta,s)d\eta ds\\
        \label{h4_7}&+\int_1^t \int \frac{Q_4}{A}e^{is\varphi}\widehat{x^3f}(\xi-\eta,s)\widehat{f}(\eta,s)d\eta ds\\
		\label{h4_5}&+\int_1^ts\int e^{is\varphi} \Big( \frac{q_6(\xi,\eta,s)}{A}\widehat{xf}(\xi-\eta,s)+ \frac{Q_5}{A}\widehat{f}(\xi-\eta,s)\Big) \widehat{f}(\eta,s)d\eta ds \\
        \label{h4_3}&+\int_1^t \int \Big( \frac{s^2 Q_9}{A}e^{is\varphi}\widehat{f}(\xi-\eta,s)+\frac{q_3(\xi,\eta,s)}{A}\widehat{x^2f}(\xi-\eta,s)\Big)\widehat{f}(\eta,s)d\eta ds \\
		\label{h4_9}&+\int_1^t\int e^{is\varphi} \Big( \frac{q_2(\xi,\eta,s)}{A}\widehat{xf}(\xi-\eta,s)+ \frac{q_1}{A} \widehat{f}(\xi-\eta,s)\Big) \widehat{f}(\eta,s)d\eta ds
	\end{align}

Note that the terms in \eqref{h4_5} appeared in \eqref{h3_2} and \eqref{h3_3}, and the terms in   \eqref{h4_3} appeared in \eqref{h3_1} and \eqref{h3_5}. Moreover, the terms in \eqref{h4_9} can be trivially bounded using \eqref{fracwuse}. 

We next decompose \eqref{h4_7} using \eqref{Decom_Bi} into the following terms 
\begin{align}
		\label{h4_7_1}&\int_1^t \int \frac{Q_4}{A}e^{is\varphi}\widehat{x^3h}(\xi-\eta,s)\widehat{h}(\eta,s)d\eta ds, \qquad \int_1^t \int \frac{Q_4}{A}e^{is\varphi}\widehat{x^3g}(\xi-\eta,s)\widehat{f}(\eta,s)d\eta ds,\\
		\label{h4_7_2}&\int_1^t \int \frac{Q_4}{A}e^{is\varphi}\widehat{x^3f}(\xi-\eta,s)\widehat{g}(\eta,s)d\eta ds,\qquad \int_1^t \int \frac{Q_4}{A}e^{is\varphi}\widehat{x^3h}(\xi-\eta,s)\widehat{g}(\eta,s)d\eta ds.
\end{align}
Using  \eqref{fracwuse} we can estimate $L^2$ norm of the first term  in \eqref{h4_7_1} by $\|f\|_X^2$. Noting that $\frac{Q_4}{A}= m_0(\xi,\eta,s)$, we  estimated the $L^2$ norm of the second term as in \eqref{q0x3g}. The two terms in \eqref{h4_7_2} can be estimated 
as 
\begin{equation*}
	\begin{split}
		\|\eqref{h4_7_2}\|_2&\lesssim\int_1^t \|T_{\frac{Q_4}{A}}(e^{-is\Delta^2}(x^3f), e^{-is\Delta^2}g)\|_2+\|T_{\frac{Q_4}{A}}(e^{-is\Delta^2}(x^3h), e^{-is\Delta^2}g)\|_2ds\\
		&\lesssim \int_1^t (\|x^3f\|_2+\|x^3h\|_2)\|e^{-is\Delta^2}g\|_{\infty}ds\lesssim \int_1^t s^{\frac{1}{2}+\frac{1}{45} -\frac{3}{2}+}ds\|f\|_X^4\lesssim t^{\frac{1}{45}+}\|f\|_X^4
	\end{split}
\end{equation*}
where we used the fact that $\|e^{-it\Delta^2}g\|_{\infty} \les t^{-\f32}\|f\|_X^2$, see Lemma \ref{lem:decay_g}.

 Therefore, we focus on the terms  \eqref{h4_1}, \eqref{h4_2} and \eqref{h4_4}. We start with \eqref{h4_1}, where we apply \eqref{key2} and have:
\begin{align}
	\eqref{h4_1}&=\label{h4_1_1}\int_1^ts^2\int \frac{q_1(\xi,\eta,s) \cdot Q_9 \varphi_\eta}{A^2}e^{is\varphi}\widehat{f}(\xi-\eta,s)\widehat{f}(\eta,s)d\eta ds\\
	\label{h4_1_2}&+t^3 \int \frac{Q_1\cdot Q_9 \varphi_\eta}{A^2}e^{it\varphi}\widehat{f}(\xi-\eta,t)\widehat{f}(\eta,t)d\eta\\
	\label{h4_1_4}&+\int_1^t s^3\int \frac{q_{9}(\xi,\eta,s)}{A}e^{is\varphi}\partial_s\{\widehat{f}(\xi-\eta,s)\widehat{f}(\eta,s)\}d\eta ds \\
	\label{h4_1_5}&+\int_1^t s^3\int \frac{Q_{11}\varphi_\eta^2}{A^2}e^{is\varphi}\widehat{f}(\xi-\eta,s)\widehat{f}(\eta,s)d\eta ds
\end{align}
 Applying integration by parts in $\eta$ to \eqref{h4_1_1} we obtain the terms in the form of 
%\begin{align} \label{ibpterm5}
%	& \int_1^ts \int \frac{q_{6}(\xi,\eta,s)}{A}e^{is\varphi}\widehat{xf}(\xi-\eta,s)\widehat{f}(\eta,s)d\eta ds  \\ 
	%& \int_1^ts \int \frac{q_9(\xi,\eta,s)}{A}e^{is\varphi}\widehat{f}(\xi-\eta,s)\widehat{f}(\eta,s)d\eta ds \label{ibpterm9}
%\end{align}
%These terms appeared in \eqref{h3_3} and  \eqref{h3_2} and respectively. 
\eqref{h3_3} and \eqref{h3_2}, whose $L^2$ norms are estimated by $t^{\frac{1}{47}}\|f\|_X^2$. Similarly, applying integration by parts in $\eta$ we can estimate $\|\eqref{h4_1_2}\|_2$ by $t^{\frac{1}{47}}\|f\|_X^2$. We have seen the term \eqref{h4_1_4} in \eqref{h2_1} and we estimated its $L^2$ norm by $t^{\frac{1}{31}}\|f\|_X^3 $.

Applying \eqref{key2} to decompose \eqref{h4_1_5}, we obtain the following terms 
\begin{align}
	&\label{h4_1_5_1}\int_1^t s^2\int \frac{q_{11}\varphi_\eta^2}{A^3}e^{is\varphi}\widehat{f}(\xi-\eta,s)\widehat{f}(\eta,s)d\eta ds +t^3\int \frac{Q_{11}\varphi_\eta^2}{A^3}e^{it\varphi}\widehat{f}(\xi-\eta,t)\widehat{f}(\eta,t)d\eta, \\
	\label{h4_1_5_4}&\int_1^t s^3\int \frac{Q_{11}\varphi_\eta^2}{A^3}e^{is\varphi}\partial_s\{\widehat{f}(\xi-\eta,s)\widehat{f}(\eta,s)\}d\eta ds \\
	\label{h4_1_5_5}& \int_1^ts^3\int \frac{Q_{12}\varphi_\eta^3}{A^3}e^{is\varphi}\widehat{f}(\xi-\eta,s)\widehat{f}(\eta,s)d\eta ds.
\end{align}
By integration by  parts in $\eta$, the first term in  \eqref{h4_1_5_1}  can be reduced to the terms in \eqref{h3_3} and \eqref{h3_2}.  Moreover, $L^2$ bound of the second term in \eqref{h4_1_5_1} follows similarly. Using the definition of $q_9(\xi,\eta,s)$, the term in \eqref{h4_1_5_4} can be reduced to  \eqref{h4_1_4}. Thus, we focus on  \eqref{h4_1_5_5} to complete estimating $\| \eqref{h4_1}\|_2$. We first rewrite \eqref{h4_1_5_5} as 
\begin{equation*}
	\eqref{h4_1_5_5}=\int_1^t \int \frac{Q_{12}}{A^3}\partial_\eta^3(e^{is\varphi})\widehat{f}(\xi-\eta,s)\widehat{f}(\eta,s)d\eta ds
\end{equation*}
Using the decomposition \eqref{Decom_Bi} for $f$, we decompose \eqref{h4_1_5_5} into the following terms 
\begin{align}
	\label{h4_a}&\int_1^t \int \frac{Q_{12}}{A^3}\partial_\eta^3(e^{is\varphi})\widehat{h}(\xi-\eta,s)\widehat{h}(\eta,s)d\eta ds\\
	\label{h4_b}&\int_1^t \int \frac{Q_{12}}{A^3}\partial_\eta^3(e^{is\varphi})\widehat{f}(\xi-\eta,s)\widehat{g}(\eta,s)d\eta ds\\
	\label{h4_c}&\int_1^t \int \frac{Q_{12}}{A^3}\partial_\eta^3(e^{is\varphi})\widehat{g}(\xi-\eta,s)\widehat{f}(\eta,s)d\eta ds\\
	\label{h4_d}&\int_1^t \int \frac{Q_{12}}{A^3}\partial_\eta^3(e^{is\varphi})\widehat{h}(\xi-\eta,s)\widehat{g}(\eta,s)d\eta ds
\end{align}

Integration by parts three times in $\eta$, we reduce studying \eqref{h4_a} to studying te following type of terms:
\begin{align}
	&\int_1^t \int e^{is\varphi} \Big( \frac{q_{1}} {A}\widehat{h}(\xi-\eta,s) +\frac{q_{2}}{A}e^{is\varphi}\widehat{xh}(\xi-\eta,s)\Big)\widehat{h}(\eta,s)d\eta  \label{h4_a_2} \\ 
&\int_1^t \int e^{is\varphi}\Big( q_0 \widehat{x^3h}(\xi-\eta,s)+ \frac{q_{3}}{A} \widehat{x^2h}(\xi-\eta,s) \Big) \widehat{h}(\eta,s)  d\eta ds \label{h4_a_3}\\
 &\int_1^t \int \frac{q_3}{A} e^{is\varphi}\widehat{xh}(\xi-\eta,s)\widehat{xh}(\eta,s)d\eta ds \label{h4_a_4}\\
	\label{h4_a_1}& \int_1^t \int q_0(\xi,\eta,s) e^{is\varphi}\widehat{x^2h}(\xi-\eta,s)\widehat{xh}(\eta,s)d\eta ds 
	\end{align}
We note that $L^2$ norm of the terms in \eqref{h4_a_2} and \eqref{h4_a_3} can be estimated using \eqref{fracwuse} by $\|f\|_X^2$. Moreover, $L^2$ norm of the term in \eqref{h4_a_4} and  \eqref{h4_a_1} can be estimated using \eqref{Decay_xh_4} and \ref{Decay_x2h_4} respectively by $\|f\|_X^2$. 

We next apply integration by parts to \eqref{h4_b} and obtain the terms in the following form 
\begin{align}
&\label{h4_b_1}\int_1^t \int e^{is\varphi}  \Big( \frac{q_{9}(\xi,\eta,s)}{A^3} \widehat{f}(\xi-\eta,s) +\frac{q_{10}(\xi,\eta,s)}{A^3}\widehat{xf}(\xi-\eta,s)\Big) \widehat{g}(\eta,s)d\eta ds \\
	\label{h4_b_3}&+\int_1^t \int e^{is\varphi} \Big(\frac{q_{10}(\xi,\eta,s)}{A^3}\widehat{f}(\xi-\eta,s) + \frac{q_{11}}{A^3} \widehat{xf}(\xi-\eta)\Big)\widehat{xg}(\eta,s)d\eta ds \\
	\label{h4_b_4}&+\int_1^t \int \frac{q_{11}}{A^3} e^{is\varphi}\widehat{x^2f}(\xi-\eta,s) \widehat{g}(\eta,s)d\eta ds + \int_1^t \int \frac{q_{12}}{A^3} e^{is\varphi} \widehat{x^3f}(\xi-\eta,s) \Big) \widehat{g}(\eta,s)d\eta ds \\ 
	\label{h4_b_6}&+\int_1^t \int \frac{q_{11}}{A^3}e^{is\varphi}\widehat{f}(\xi-\eta,s)\widehat{x^2g}(\eta,s)d\eta ds +\int_1^t \int q_0 e^{is\varphi}\widehat{x^2f}(\xi-\eta,s)\widehat{xg}(\eta,s)d\eta ds\\
	\label{h4_b_9}&+\int_1^t \int \frac{q_{12}}{A^3}e^{is\varphi}\widehat{xf}(\xi-\eta,s)\widehat{x^2g}(\eta,s)d\eta ds\\
	\label{h4_b_10}&+\int_1^t \int q_0 e^{is\varphi}\widehat{f}(\xi-\eta,s)\widehat{x^3g}(\eta,s)d\eta ds
\end{align}
Using \eqref{eqn_g}, the definition of \( \widehat{g} \), we observe that \eqref{h4_b_1} reduces to \eqref{h1_2_11} and \eqref{h1_2_2}. The \( L^2 \) norms of these terms are estimated by \( t^{\f1{48}+}\|f\|_X^3 \). Moreover, \eqref{h4_b_10} appeared in \eqref{q0x3g}.  
 Similarly, using \eqref{eqn_g}, we can express 
\begin{align}
	\eqref{h4_b_3}&=\label{170_1}\int_1^t s \iint \frac{q_{2+\alpha}(\xi,\eta,s)}{A}\frac{q_3(\eta,\sigma,s)}{C}e^{is\psi}\widehat{x^{\alpha}f}(\xi-\eta,s)\widehat{f}(\eta-\sigma,s)\widehat{f}(\sigma)d\sigma d\eta \\
	\label{170_2}&+\int_1^t \iint \frac{q_{2+\alpha}}{A}\frac{1}{C}e^{is\psi}\widehat{x^{\alpha}f}(\xi-\eta,s)\widehat{xf}(\eta-\sigma,s)\widehat{f}(\sigma)d\sigma d\eta ds
\end{align}
for $\alpha=0,1$. When $\alpha=0$ \eqref{170_2} appeared in \eqref{h1_2_2} and its $L^2$ norm is estimated by $t^{\f1{48}+}\|f\|_X^3$. Moreover, when $\alpha=1$, then it reduces to $\eqref{102_2}$. We can estimate \eqref{170_2} for $\alpha=1$ using \eqref{Decay_xf} for $\alpha =1$ and \eqref{6.44est}, and we have the following for $\alpha=0$. 
\begin{align*}
\|\nabla_s^{-1}u\|_{\f52}\|u\|_{20} \les s^{-\f74}\|f\|_X^2 \|u\|_{5}\les s^{-\f{10}{4}} \|f\|_X^3 
\end{align*}

We continue with the terms in \eqref{h4_b_4}. The first one  reduces to \eqref{h1_2_4} upon using \eqref{eqn_g}, whose \( L^2 \) norm is estimated by \( t^{\frac{1}{48}+} \|f\|_X^3 \). The second term in \eqref{h4_b_4} can be estimated as

	 \begin{equation*}
	 	\int_1^t \|e^{-is\Delta^2}(x^3f)\|_2\|e^{-is\Delta^2}g\|_\infty ds\lesssim \int_1^t s^{-\f54} \|x^{\f52+}g\|_2\|f\|_2ds\lesssim t^{\frac{1}{90}+}\|f\|_X^3.
	 \end{equation*}
Moreover, using \eqref{Decay_x2g}, we can bound $L^2$ norm of the first term in \eqref{h4_b_6} by 
 \begin{align*} \int_1^t s^{\f14}\|u\|_{\infty} \|x^2g\|_2 ds \les t^{\varepsilon_0} (\|f\|_X^3 +\|f\|_X^2)
\end{align*}

We use \eqref{eqn_g} for the second term in  \eqref{h4_b_6} and obtain the following type of terms 
\begin{align}
	& \label{176_1}\int_1^t \iint q_0(\xi,\eta,s) \frac{1}{C}	e^{is\psi}\widehat{x^2f}(\xi-\eta,s)\widehat{xf}(\eta-\sigma,s)\widehat{f}(\sigma,s)d\sigma d\eta ds\\
	\label{176_2}& \int_1^t s \iint q_0(\xi,\eta,s) \frac{q_3(\eta,\sigma,s)}{C}	e^{is\psi}\widehat{x^2f}(\xi-\eta,s)\widehat{f}(\eta-\sigma,s)\widehat{f}(\sigma,s)d\sigma d\eta ds 
\end{align}
Using \eqref{decay_xf_4} we can estimate 
\begin{align*}
    \|\eqref{176_1}\|_2 \les \int_1^t \|e^{-is\Delta^2} (xf)\|_4 \|u\|_4  \: ds \|f\|_{X}\les  t^{-\frac{20}{128}}\|f\|^3_{X}
\end{align*}
and using Lemma~\ref{lem_xf_4} and \ref{Decay_x2f_4}, we obtain
\begin{equation}
	\begin{split}
		\|\eqref{176_1}\|_2\lesssim\int_1^t s\|e^{-is\Delta^2}(x^2f)\|_4\|e^{-is\Delta^2}(xf)\|_4\|u\|_\infty ds\lesssim t^{\frac{5}{384}+\varepsilon}\|f\|_X^3\lesssim t^\frac{1}{48}\|f\|_X^3.
	\end{split}\label{est6140}
\end{equation}
We again use \eqref{eqn_g} in \eqref{h4_b_9} to reduce the problem estimating the following terms 
\begin{align}
	& \label{191}\int_1^t s^2\iint q_0(\xi,\eta,s) \frac{q_6(\eta,\sigma,s)}{C} e^{is\psi}\widehat{xf}(\xi-\eta,s)\widehat{f}(\eta-\sigma,s)\widehat{f}(\sigma,s) d\sigma d\eta ds\\
	& \int_1^t \iint q_0(\xi,\eta,s)\frac{q_0(\eta,s)}{C} e^{is\psi}\widehat{xf}(\xi-\eta,s)\widehat{x^2f}(\eta-\sigma,s)\widehat{f}(\sigma,s) d\sigma d\eta ds \label{192}\\
	& \int_1^t s \iint q_0(\xi,\eta,s)\frac{q_3(\eta,\sigma,s)}{C} e^{is\psi}\widehat{xf}(\xi-\eta,s)\widehat{xf}(\eta-\sigma,s)\widehat{f}(\sigma,s) d\sigma d\eta ds \label{193}\\
	& \int_1^t s \iint  q_0(\xi,\eta,s)\frac{q_2(\eta,\sigma,s)}{C} e^{is\psi}\widehat{xf}(\xi-\eta,s)\widehat{f}(\eta-\sigma,s)\widehat{f}(\sigma,s) d\sigma d\eta ds \label{194}
\end{align}

We can bound $\|\eqref{192}\|_2$ using the same estimate as in \eqref{est6140}. Using \eqref{decay_xf_4} $\|\eqref{193}\|_2$ can be estimated by $ t^{-\f5{64}}\|f\|_X^{3}$, and $\|\eqref{194}\|_2$ can be bounded by $t^{\f1{48}} \|f\|_X^3$ using \eqref{Decay_xf} for $\alpha=1$ together with \eqref{6.44est}. Moreover,  we see that the term in \eqref{191} can be written as the sum of  $G_{6,i}(\widehat{xf},\widehat{f}, \widehat{f})$ which we defined in \eqref{G51}. By symmetry, we can estimate 
\begin{align*}
    \|G_{6,\ell}(\widehat{xf},\widehat{f}, \widehat{f}\|_2 \les \int_1^t s^2 \sum_{j}2^{2j}\|T_{m^*}(e^{-is\Delta^2}(xf), P_ju, u)\|_2 \les t^{0+}\|f\|^3_X
\end{align*}
where we used similar estimates as in \eqref{s22jxf} in the last inequality.

Applying a change of variables, the bound on \( \|\eqref{h4_c}\|_2 \) follows in the same way as the estimate for \( \|\eqref{h4_b}\|_2 \). Morever, by Lemma~\ref{dsh} and the fact that  the norm estimates for \( h \) exhibit slower time growth than those for \( f \) we can employ the same strategy used for \eqref{h4_b} to estimate $\|\eqref{h4_d}\|_2$, see Remark~\ref{rmkdsh2}. 

\begin{remark} \label{rmkdsh2} Similar to the reasoning in Remark~\ref{rmk:dsh}, strategy used for \eqref{h4_b} to estimate $\|\eqref{h4_d}\|_2$ particular care is needed in handling the term \eqref{h4_b_10} when $f$ is exchanged with $h$. Specifically, this term is \eqref{q0x3g} and we decompose it in terms including \eqref{108_1}. Moreover, the estimate $\|\eqref{108_1}\|_2$ relies on \eqref{key2} and requires the control  of $\|\eqref{h3_1_4_n}\|_2$. Noting \eqref{G911} and \eqref{G9111} to complete the required bound on $\|\eqref{h4_d}\|_2$ we prove Lemma~\ref{dsh}. 
\end{remark}

We continue with \eqref{h4_4}. Applying integration by parts in $\eta$
\begin{align}
\eqref{h4_4} =\eqref{h4_7}&+\int_1^t \int \frac{q_3}{A}e^{is\varphi}\widehat{x^2f}(\xi-\eta,s)\widehat{f}(\eta,s)d\eta ds \nn \\ &+\int_1^t \int q_0(\xi,\eta,s) e^{is\varphi}\widehat{x^2f}(\xi-\eta,s)\widehat{xf}(\eta,s)d\eta ds.  \label{611est} 
\end{align} 
The second term appeared in \eqref{h3_5}. We decompose the last term using \eqref{Decom_Bi} into:
\begin{align}
	\label{h4_4_1}&\int_1^t \int q_0(\xi,\eta,s) e^{is\varphi}\widehat{x^2h}(\xi-\eta,s)\widehat{xh}(\eta,s)d\eta ds\\
	\label{h4_4_2}&\int_1^t \int q_0(\xi,\eta,s) e^{is\varphi}\widehat{x^2f}(\xi-\eta,s)\widehat{xg}(\eta,s)d\eta ds\\
	\label{h4_4_3}&\int_1^t \int q_0(\xi,\eta,s) e^{is\varphi}\widehat{x^2g}(\xi-\eta,s)\widehat{xf}(\eta,s)d\eta ds\\
	\label{h4_4_4}&\int_1^t \int q_0(\xi,\eta,s) e^{is\varphi}\widehat{x^2h}(\xi-\eta,s)\widehat{xf}(\eta,s)d\eta ds
\end{align}
\eqref{h4_4_1} and \eqref{h4_4_2} are appeared in \eqref{h4_a_1} and \eqref{h4_b_6}. By change of variable, \eqref{h4_4_3} is  the same term as in \eqref{h4_b_9}. As $h$ behaves better than $f$, $\|\eqref{h4_4_4}\|_2$ can be estimated like $\|\eqref{h4_b_9}\|_2$. This completes the estimate of $\|\eqref{h4_4}\|_2$. 

We next focus on \eqref{h4_2}, the last term in $\partial_{\xi}^3(h_4)$. We apply integration by parts in $\eta$ to see 
\begin{align}
	\eqref{h4_2}&=\label{h4_2_1}\int_1^t s\int \frac{q_6}{A}e^{is\varphi} \widehat{xf}(\xi-\eta,s)\widehat{f}(\eta,s)d\eta ds	\\
	\label{h4_2_3}& + \int_1^t s\int \frac{Q_7}{A} e^{is\varphi}\widehat{xf}(\xi-\eta,s)\widehat{xf}(\eta,s)d\eta ds\\
	\label{h4_2_2}&+ \int_1^t s\int \frac{Q_7}{A}e^{is\varphi} \widehat{x^2f}(\xi-\eta,s)\widehat{f}(\eta,s)d\eta ds 
	\end{align}
\eqref{h4_2_1} has been estimated in \eqref{h3_2}. We apply \eqref{key2} to reduce estimating  $\|\eqref{h4_2_3}\|_2$  to estimating $L^2$ norm of the following terms: 
\begin{align}
 &\label{h4_2_3_1}\int_1^t \int \frac{Q_7}{A^2} 	e^{is\varphi}\widehat{xf}(\xi-\eta,s)\widehat{xf}(\eta,s)d\eta ds+ \int \frac{Q_7}{A^2}e^{it\varphi}\widehat{xf}(\xi-\eta,t)\widehat{xf}(\eta,t)d\eta, \\
\label{h4_2_3_3}&\int_1^t\int \frac{Q_8}{A^2}e^{is\varphi}\partial_\eta(\widehat{xf}(\xi-\eta,s)\widehat{xf}(\eta,s))d\eta ds\\
\label{h4_2_3_2}, &\int_1^t s\int \frac{Q_7}{A^2} e^{is\varphi}\partial_s(\widehat{xf}(\xi-\eta,s)\widehat{xf}(\eta,s))d\eta ds.
\end{align}
The $L^2$ norm of the first term in \eqref{h4_2_3_1} is estimated in \eqref{xfxf}. Note that we do not use the integration in $s$ for estimating \eqref{xfxf}, so the $L^2$ norm of the  second term in \eqref{h4_2_3_1} can be calculated similarly. The term \eqref{h4_2_3_3} is the last term in \eqref{611est}. To estimate the $L^2$ norm of  \eqref{h4_2_3_2}, we apply integration by parts in $\eta$ and obtain the terms in \eqref{x2h_2_p6} and \eqref{x2h_2_f} whose $L^2$ norms are estimated by $\log(t) ( \|f\|_X^2 +\|f\|_X^3+\|f\|_X^4) $. 
 
 We next estimate $\|\eqref{h4_2_2}\|_2$. Applying \eqref{key2}, we reduce estimating $\|\eqref{h4_2_2}\|_2$ to estimating the following terms: 
\begin{align}
	& \label{h4_2_2_1}\int_1^t \int \frac{q_3}{A} e^{is\varphi} \widehat{x^2f}(\xi-\eta,s)\widehat{f}(\eta,s)d\eta ds+t\int  \frac{q_3}{A}e^{it\varphi} \widehat{x^2f}(\xi-\eta,t)\widehat{f}(\eta,t)d\eta\\
	\label{h4_2_2_2}&\int_1^t s\int \frac{Q_7}{A^2}e^{is\varphi} \widehat{x^2f}(\xi-\eta,s)\partial_s\widehat{f}(\eta,s)d\eta ds\\
	\label{h4_2_2_3}&\int_1^t \int \frac{Q_8}{A^2}e^{is\varphi}\left(\widehat{x^3f}(\xi-\eta,s)\widehat{f}(\eta,s) +\widehat{x^2f}(\xi-\eta,s)\widehat{xf}(\eta,s)\right)d\eta ds\\
\label{h4_2_2_4}	&\int_1^t s\int \frac{Q_7}{A^2}e^{is\varphi} \partial_s(\widehat{x^2f}(\xi-\eta,s))\widehat{f}(\eta,s)d\eta ds
\end{align}
\eqref{h4_2_2_1} has appeared in \eqref{h3_5} and \eqref{x3g_5} while \eqref{h4_2_2_2} is also estimated in \eqref{h2_4} already. The first part of \eqref{h4_2_2_3} can be estimated using \eqref{weighted} and the second part is estimated in \eqref{h4_2_3_3}. 

Using change of variable and integration by parts in $\eta$, we rewrite \eqref{h4_2_2_4} as :
\begin{align*}
	\eqref{h4_2_2_4}&=\int_1^t s\int \frac{q_2}{A^\prime}e^{is\varphi}\widehat{f}(\xi-\eta,s)\partial_\eta\partial_s\widehat{f}(\eta,s)d\eta ds \stepcounter{equation}\tag{\theequation}\label{163_1}\\
	&+\int_1^t s\int \frac{Q_7}{(A^\prime)^2} e^{is\varphi}\widehat{xf}(\xi-\eta,s)\partial_\eta\partial_s\widehat{f}(\eta,s)d\eta ds \stepcounter{equation}\tag{\theequation}\label{163_2}\\
	&+\int_1^t s^2\int \frac{Q_{10}}{(A^\prime)^2}e^{is\varphi}\widehat{f}(\xi-\eta,s)\partial_\eta\partial_s\widehat{f}(\eta,s)d\eta ds \stepcounter{equation}\tag{\theequation}\label{163_3}
\end{align*}
\eqref{163_1} appeared in \eqref{xfxf}. We apply integration by parts for the second time to \eqref{163_2} to see
\begin{align}
	\eqref{163_2}&=\label{163_2_1}\int_1^t s\int \frac{q_2}{A^\prime}e^{is\varphi}\widehat{xf}(\xi-\eta,s)\partial_s\widehat{f}(\eta,s)d\eta ds\\
\label{163_2_2}	&+\int_1^t s^2\int \frac{Q_{10}}{(A^\prime)^2}e^{is\varphi}\widehat{xf}(\xi-\eta,s)\partial_s\widehat{f}(\eta,s)d\eta ds\\
	\label{163_2_3}&+\int_1^t s\int \frac{Q_7}{(A^\prime)^2} e^{is\varphi}\widehat{x^2f}(\xi-\eta,s)\partial_s\widehat{f}(\eta,s)d\eta ds
\end{align}
These terms appeared in \eqref{h2_4}, \eqref{h2_2} and \eqref{h2_2} respectively.

We finally focus on  \eqref{163_3}.
Observe that:
\begin{equation*}
\begin{split}
	\partial_\eta\partial_s\widehat{f}(\eta,s)&=i\partial_\eta\left(\int e^{is\varphi(\eta,\sigma)}\widehat{f}(\eta-\sigma,s)\widehat{f}(\sigma,s)d\sigma \right)\\
	&=s\int Q_3(\eta,\sigma)  e^{is\varphi(\eta,\sigma)}\widehat{f}(\eta-\sigma,s)\widehat{f}(\sigma,s)d\sigma+i \int e^{is\varphi(\eta,\sigma)}\widehat{f}(\eta-\sigma,s)\widehat{xf}(\sigma,s)d\sigma
\end{split}
\end{equation*}
where we applied also change of variable to the second term. Plugging in \eqref{163_3}, we have two terms:
\begin{align}
\label{163_3_1}	&\int_1^t s^2\iint \frac{q_{6}(\xi,\eta)}{A^\prime}e^{is\psi}\widehat{f}(\xi-\eta,s)\widehat{f}(\eta-\sigma,s)\widehat{xf}(\sigma,s)d\sigma\\
\label{163_3_2}&\int_1^t s^3\iint \frac{Q_{10}(\xi,\eta)}{(A^\prime)^2}Q_3(\eta,\sigma)e^{is\psi}\widehat{f}(\xi-\eta,s)\widehat{f}(\eta-\sigma,s)\widehat{f}(\sigma,s)d\sigma
\end{align}
We decompose \eqref{163_3_1} in to $F_{6,1}(\widehat{f},\widehat{f},\widehat{xf})$ and $F_{6,2,j}(\widehat{f},\widehat{f},\widehat{xf})$ where $F_{i,\ell}$ and $F_{i,\ell,k}$ are defined as in \eqref{defF1}, \eqref{defF2}. Then we can bound 
\begin{align*}
    \|F_{6,1}(\widehat{f},\widehat{f},\widehat{xf}) \|_2 & \les \int_1^t s^2\sum_j 2^{2j}\|T_{m^*}(P_j(u), u, e^{-is\Delta^2}(xf))\|_2ds \les t^{0+} \|f\|_X^3
\end{align*}
where we used \eqref{s22jxf} in the last line. 
Moreover, using Young's inequality and \eqref{s22jxf} we can estimate 
\begin{multline*}
 \|F_{6,2,1}(\widehat{f},\widehat{f},\widehat{xf})\|_2 \les \int_1^t s^2\sum_j\sum_{k\leq j+10} 2^{2k} \|T_{m^*}(P_k(u), P_j(u), e^{-is\Delta^2}(xf))\|_2  \: ds 
\les t^{0+} \|f\|_X^3.  
\end{multline*}
Finally, the following estimate
\begin{align*}
  \int_1^t s^2\sum_{s^{-1}\leq 2^j \leq s}& \sum_{k\leq j+10} 2^{2k}  \|T_{m^*}(P_k(u), u, P_j(e^{-is\Delta^2}(xf)))\|_2 \: ds  \\
&\les  \int_1^t s^2\sum_{s^{-1}\leq 2^j \leq s} 2^{2j+}\|P_je^{-is\Delta^2}(xf)\|_\frac{10}{3}\|u\|_{\infty}\|u\|_{5}ds \les  \int_1^t s^{-1+\varepsilon_0+} \|f\|_X^3 \: ds \les t^{\varepsilon_0+} \|f\|_X^3
\end{align*}
together with Remark~\ref{middle_j} gives $\|F_{6,2,2}(\widehat{f},\widehat{f},\widehat{xf})\|_2 \les t^{\varepsilon_0+} \|f\|_X^3 $. 
To estimate $\| \eqref{163_3_2}\|_2$, we apply \eqref{key3} and reduce the problem estimating $L^2$ norm of the following terms: 
\begin{align*}
	&\int_1^t s^2\iint \frac{q_5}{A^\prime}q_0(\xi,\eta,\sigma,s)e^{is\psi}\widehat{f}(\xi-\eta,s)\widehat{f}(\eta-\sigma,s)\widehat{f}(\sigma,s)d\sigma d\eta ds\\
&t^3\iint  \frac{q_5}{A^\prime}q_0(\xi,\eta,\sigma,s)e^{it\psi}\widehat{f}(\xi-\eta,t)\widehat{f}(\eta-\sigma,t)\widehat{f}(\sigma,t)d\sigma d\eta\\
&\int_1^t s\iint \frac{Q_{6}(\xi,\eta)}{A^\prime }\frac{Q_4(\xi,\eta,\sigma) }{B}e^{is\psi}\widehat{xf}(\xi-\eta,s)\widehat{f}(\eta-\sigma,s)\widehat{f}(\sigma,s)d\sigma d\eta ds\\
	&\int_1^t s\iint \frac{Q_{6}(\xi,\eta)}{A^\prime }\frac{Q_4(\xi,\eta,\sigma)}{B}e^{is\psi}\widehat{f}(\xi-\eta,s)\widehat{xf}(\eta-\sigma,s)\widehat{f}(\sigma,s)d\sigma d\eta ds 
\end{align*}
where we used $ Q_1(\xi,\eta) Q_3(\eta,\sigma)= Q_4(\xi,\eta,\sigma)$ in the first term.  Notice that $L^2$ norm of the last   term can be estimated similar to $\| \eqref{163_3_1}\|_2$. Similar terms ($A$ is exchanged with $A$) to the first two terms  appeared in \eqref{h2_1_1} and \eqref{h2_1_2} respectively. Moreover, the third term can be seen as \eqref{h2_1_5} since we consider $Q_9(\xi,\eta)Q(\eta,\sigma)= Q_6(\xi,\eta)Q_4(\xi,\eta,\sigma)$ in \eqref{h2_1_5}. 

\end{proof}

\begin{lemma} \label{dsh} Assuming \eqref{Boot_est_g} and \eqref{Boot_est_h} we have 
\begin{align*}
\|e^{is\Delta^2} (\partial_sh)\|_{p} \les s^{-\f52+ \f 5{p}}, \,\,\,\ \|e^{is\Delta^2} (\partial_sh)\|_{10}\les s^{-\f{13}{4}} 2^{-4j} \|f\|_X\|f\|_{\frac{10}{3}}
\end{align*}
Hence, $\|e^{is\Delta^2} (\partial_sh)\|_{p}$ satisfies the same bound that $\|u^2\|_p$ satisfies and $\|e^{is\Delta^2} (\partial_sh)\|_{10}$ holds \eqref{Pj_u2_10}. 

\end{lemma}
\begin{proof}
Using $h=f-g-f_*$ and noticing that $\partial_s\widehat{f_*}=0$, we have $$\partial_s(\widehat{h})=\partial_s(\widehat{f})- \partial_s(\widehat{g}) = e^{it|\xi|^4}\widehat{u^2} -\partial_s(\widehat{g}).  $$ Therefore, it suffices to prove the statement with \( \partial_s h \) replaced by \( \partial_s g \).
To do that, we first estimate $\|e^{-it\Delta^2}( \partial_s(g))\|_p$ for any $0<p<\infty$. Note that we have 

\begin{equation*}
	\begin{split}
		(e^{-is\Delta^2}(\partial_sg))^{\widehat{}}(\xi,s)&= \int m_0(\xi,\eta,s) \widehat{u}(\xi-\eta,s)\widehat{u}(\eta,s)d\eta \\
		&\ +\int \frac{1}{A}\widehat{u}(\xi-\eta,s)\widehat{u^2}(\eta,s)d\eta+\int \frac{1}{A}\widehat{u^2}(\xi-\eta,s)\widehat{u}(\eta,s)d\eta\\
		&=\widehat{I}+\widehat{II}+\widehat{III}
	\end{split}
\end{equation*}
Therefore, for any $0<p<\infty$
\begin{align} \label{dsgbound}
\|e^{-it\Delta^2}( \partial_s(g))\|_p & \les \|T_{m_0(\xi,\eta,s)}(u,u)\|_p+ s \|T_{m_0(\xi,\eta,s)}(u,u^2)\|_p  \nn \\
& \les \|u\|_{\infty} \|u\|_p + s \|u\|_{\infty}^2 \|u\|_{p}\les s^{-\f52 + \frac{5}{2p}} \|f\|_X^2 
\end{align}
Hence, $\|e^{-it\Delta^2}( \partial_s(g))\|_p$ satisfies the bounds that $\|u^2\|_p$ satisfies and so $\|e^{-it\Delta^2}( \partial_s(h))\|_p$. 

 Next we show that $\|P_je^{-is\Delta^2}(\partial_sg)\|_{10} \les s^{-\frac{13}{4}} 2^{-4j} \|f\|_X \|f\|_{\frac{10}{9}}$. We have 

\begin{equation*}
	\begin{split}
		P_j(I)&=P_jT_{m_0}(u,u)=P_jT_{m_0}(P_{< j-5}u,P_{j-3\leq n\leq j+3}u)+\sum_{k\geq j-5}P_jT_{\frac{P_4}{A}}(P_ku,u).
	\end{split}
\end{equation*}
Therefore,
\begin{equation*}
	\begin{split}
		\|P_jI\|_{10}&\leq \|T_{m_0}(P_{< j-5}u,P_{j-3\leq n\leq j+3}u)\|_{10}+\sum_{k\geq j-5}\|T_{m_0}(P_ku,u)\|_{10}\\
		&\les \|u\|_\infty (\sum_{n=j-3}^{j+3}\|P_nu\|_{10}+\sum_{k\geq j-5} \|P_ku\|_{10})\les s^{-\frac{13}{4}}2^{-4j}\|f\|_X\|f\|_\frac{10}{9}.
	\end{split}
\end{equation*}
Moreover, 
\begin{equation*}
	\begin{split}
		\|P_j(II+III)\|_{10}&\les \|T_\frac{1}{A}(P_{<j-5}u, P_{j-3\leq n\leq j+3}(u^2))\|_{10}+\sum_{k\geq j-5}\|T_\frac{1}{A}(P_ku, u^2)\|_{10}\\
		&\les s(\|u\|_\infty \sum_{n=j-3}^{j+3}\|P_n(u^2)\|_{10}+\sum_{k\geq j-5} \|P_ku\|_{10}\|u^2\|_\infty) \les s^{-\frac{7}{2}}2^{-4j}\|f\|_\frac{10}{9}\|f\|_X^2. 
		\end{split}
\end{equation*}
\end{proof}
Finally, we can give the following Corollary. 

\begin{corollary}\label{lem_x3h}
	Assuming \eqref{Boot_est_g} and \eqref{Boot_est_h}, we have
	\begin{equation*}
		\|x^3h\|_2\les t^{\frac{1}{24}}\|f\|_X^2
	\end{equation*}
\end{corollary} 
\begin{proof}
The proof follows from the estimates \eqref{Decay_x3h1}, \eqref{Decay_x3h2},\eqref{Decay_x3h3} and \eqref{Decay_x3h4}. 
\end{proof}

Therefore, we have closed the bootstrap and established the estimates in \eqref{est_g} and \eqref{est_h}. As a corollary, we now have the dispersive decay for $h$ and thus $f$:

\begin{lemma}\label{lem_decay_B}We have:
\begin{equation*}
\begin{split}
	\|e^{-it\Delta^2}h\|_{L^\infty}&\lesssim t^{-\frac{5}{4}}\|f\|_X^2\\
	\|e^{-it\Delta^2 }B(f,f)\|_{L^{\infty}} &\les t^{-\frac{5}{4}}\|f\|^2_X 
\end{split}
\end{equation*}
\end{lemma}
\begin{proof}
To bound $e^{-it\Delta^2}h$, we use Lemma~\ref{lem:linear} to estimate 
\begin{align*}
\|e^{-it\Delta^2}h\|_{L^{\infty}} & \les t^{-\f54} \|\hat{h}\|_{L^{\infty}} + t^{-\frac{11}{8}} \| \hat{h}\|_{H^3} \\ & \les t^{-\f54}[ \| \hat{B}(f,f) \|_{L^{\infty}}+ \|\hat{g}\|_{L^{\infty}} + \|f_*\|_{L^{\infty}}]
 + t^{-\frac{11}{8}} \| \hat{h}\|_{H^3}. 
\end{align*}
which follows from the decomposition in \eqref{decompB}. By Corollary~\ref{corsupB}, we have the required bound for the first term on the right side of the second inequality. Moreover, we have 
\begin{equation*}
		\begin{split}
			|\widehat{g}(\xi,t)|&\leq \int \frac{1}{|\eta|^4}|\widehat{f}(\xi-\eta,t) ||\widehat{f}(\eta,t) |d\eta\les \|\widehat{f}\|_{L^2_\xi\cap L^\infty_\xi}^2\leq \|f\|_X^2
		\end{split}
	\end{equation*}
for the second term. Finally in Section~\ref{sec:hbounds} we show that \( \| \hat{h}\|_{H^3} \les t^{\f1{24}} \| f\|_X^2 \), which proves 
\[ \|e^{it\Delta^2}h\|_{L^{\infty}}\les t^{-\f54} \| f\|_X^2.\] Combined with Lemma \ref{lem:decay_g}, we have
\begin{equation*}
	\begin{split}
		\|e^{-it\Delta^2}B(f,f)\|_{L^\infty}&\leq\|e^{-it\Delta^2}g\|_{L^\infty} +\|e^{-it\Delta^2}h\|_{L^\infty}+\|e^{-it\Delta^2}f_*\|_{L^\infty}
	\end{split}
\end{equation*}
The last term does not depend on $t$, and estimating it can be reduced to the smallness of the initial data. Therefore, summing up the terms, we proved:
\begin{equation*}
	\|e^{-it\Delta^2 }B(f,f)\|_{L^{\infty}} \les t^{-\frac{5}{4}}\|f\|^2_X 
\end{equation*}
\end{proof}

\section{Proof of the main theorem}\label{proof}
Using the estimates we achieved in sections \ref{sec:gbounds} and \ref{sec:hbounds}, we can close the bootstrap with assumptions \eqref{Boot_est_g} and \eqref{Boot_est_h}. In summary, we proved all the estimates in \eqref{est_g} and in \eqref{est_h}. Adding these together, we have:
\begin{equation}\label{est_f}
	\begin{split}
	\|f\|_{L^2}\lesssim \|f\|_{X}^2, \qquad \|xf\|_{L^2}&\lesssim \|f\|_{X}^2, \qquad \|x^2f\|_{L^2}\lesssim\|f\|_X^2+\|f\|_X^3, \\
		\qquad \|x^3f\|_{L^2}\lesssim t^{\frac{1}{2}+\frac{1}{47}} \|f\|_{X}^2, &\qquad \|e^{-it\Delta^2}f\|_{L^\infty}\lesssim t^{-\frac{3}{2}+}\|f\|_{X}^2
	\end{split}
\end{equation}
We are left only one term in the $X$-norm: $\|f\|_{L^\infty H^{10}}$.
\begin{equation*}
	\begin{split}
		\|f(t)\|_{H^{10}}&\leq \|f_*\|_{H^{10}}+\int_1^t \|e^{is\Delta^2}(u^2)\|_{H^{10}}ds \leq \|f_*\|_{H^{10}}+\int_1^t \|u\|_{H^{10}}\|u\|_\infty  ds\\
		&\leq \|f_*\|_{H^{10}}+\int_1^t s^{-\frac{5}{4}}  ds\|f\|_X^2 \lesssim \|f_*\|_{H^{10}}+\|f\|_X^2
	\end{split}
\end{equation*}
so we close the bootstrap for $\|f\|_{L^\infty H^{10}}$. Using the estimates \eqref{est_f} that we have established, we have:
\begin{equation*}
	\hat{f}(\xi,t)\rightarrow\widehat{f_1}(\xi)+\int_1^t\int e^{is\varphi}\widehat{f}(\xi-\eta,s)\hat{f}(\eta,s)d\eta ds
\end{equation*} 
maps a small neighborhood of the origin in $X$ into itself. It is easy to notice that our computation can be modified to show that the map is also a contraction on a small neighborhood in the $X$ space. Thus, the existence of solution is proved. Notice that the $L^\infty$ decay of $u$ is integrable, so scattering follows. 

\bibliographystyle{amsplain}

\bibliography{references}

\end{document}